\documentclass[12pt,letterpaper]{article}
\usepackage{epsfig}
\usepackage[]{natbib}
\usepackage{amsmath,amsthm}
\usepackage{amssymb,amsfonts}
\usepackage{nccfoots,bbm,bm}
\usepackage{graphicx,color,caption}
\usepackage{booktabs,multirow,enumerate}
\usepackage{tikz,subfigure}
\usepackage{microtype}
\usepackage[colorlinks,citecolor=blue]{hyperref}
\usepackage[active]{srcltx}
\usepackage{xr}
\usepackage{comment}
\usepackage{longtable}
\usepackage{verbatim}

\usepackage[hmargin=1in,vmargin=1in]{geometry}
\newtheorem{lem}{Lemma}[section]

\newtheorem{thm}{Theorem}[section]
\newtheorem{example}{Example}[section]

\newtheorem{rem}{Remark}[section]
\numberwithin{equation}{section}
\numberwithin{figure}{section}
\numberwithin{table}{section}
\theoremstyle{plain}

\newcommand{\beq}{\begin{equation}}
\newcommand{\eeq}{\end{equation}}

\newcommand{\var}{\text{\normalfont Var}}
\newcommand{\cov}{\text{\normalfont Cov}}
\newcommand{\ee}{\mathbb{E}}
\newcommand{\pp}{\mathbb{P}}
\newcommand{\rr}{\mathbb{R}}
\newcommand{\nn}{\mathbb{N}}
\newcommand{\zz}{\mathbb{Z}} 

\def\qed{\rule{2mm}{2mm}}

\usepackage[section]{placeins}

\usepackage{graphicx}
\graphicspath{ {./Plots}
}

\begin{document}

\title{Permutation Testing for Dependence in Time Series}

\author{Joseph P. Romano\footnote{Supported by NSF Grant MMS-1949845. We thank Kevin Guo and Benjamin Seiler for helpful comments and conversations.}
  \\
Departments of Statistics and Economics \\
Stanford University \\
\href{mailto:romano@stanford.edu}{romano@stanford.edu}
\and
Marius A. Tirlea \\
Department of Statistics \\
Stanford University \\
\href{mailto:mtirlea@stanford.edu}{mtirlea@stanford.edu} 
}

\date{\today}

\maketitle

\begin{abstract}

Given observations from a stationary time series, permutation tests allow one to construct exactly level $\alpha$ tests under the null hypothesis of an i.i.d. (or, more generally, exchangeable) distribution. On the other hand, when the null hypothesis of interest is that the underlying process is an uncorrelated sequence, permutation tests are not necessarily level $\alpha$, nor are they approximately level $\alpha$ in large samples. In addition, permutation tests may have large Type 3, or directional, errors, in which a two-sided test rejects the null hypothesis and concludes that the first-order autocorrelation is larger than 0, when in fact it is less than 0. In this paper, under weak assumptions on the mixing coefficients and moments of the sequence, we provide a test procedure for which the asymptotic validity of the permutation test holds, while retaining the exact rejection probability $\alpha$ in finite samples when the observations are independent and identically distributed. A Monte Carlo simulation study, comparing the permutation test to other tests of autocorrelation, is also performed, along with an empirical example of application to financial data.
\end{abstract}

\bigskip

\begin{tabbing}
\end{tabbing}


\newpage

\section{Introduction}

In this paper, we investigate the use of permutation tests for detecting dependence in a time series.
When testing the null hypothesis that the underlying time series consists of independent, identically distributed (i.i.d.) observations, permutation tests can be constructed that control the probability of a Type 1 error exactly, for any choice of test statistic. Typically, the choice of test statistic is the first-order sample autocorrelation, or some function of many of the sample autocorrelations.
However, significant problems of error control arise, stemming from the fact that zero autocorrelation and independence are actually quite different. It is crucial to carefully specify the null hypothesis of interest, whether it is the case that the observations are i.i.d. or that the observations have zero autocorrelation. For example, if the null hypothesis specifies that the autocorrelation is zero and one uses the sample first-order autocorrelation as a test statistic when applying a permutation test, then the Type 1 error can be shockingly different from the nominal level, even asymptotically.  Nevertheless, one might think it reasonable to reject based on such a permutation test, since the test statistic appears ``large", relative to the null reference permutation distribution. However, even  if one views the null hypothesis as specifying that the time series is i.i.d., a rejection of the null hypothesis based on the sample autocorrelation is inevitably accompanied by the interpretation then that the true underlying autocorrelation is nonzero. Indeed, one typically makes the further claim that this correlation is positive (negative), when the sample autocorrelation can be large and positive (large and negative).
When the true autocorrelation is zero and there is a large probability of a Type 1 error, then lack of Type 1 error control also implies lack of Type 3, or directional, error control. That is, there can be a large probability that one declares the underlying correlation to be positive when it is in fact negative.

Assume $X_1, \, \dots, \, X_n$ are jointly distributed according to some strictly stationary, infinite dimensional distribution $P$, where the distribution $P$ belongs to some family $\Omega$. Consider the problem of testing $H: P \in \Omega_0$,  where $\Omega_0$ is some subset of stationary proceseses.  For example, we might be interested in testing that the underlying $P$ is a product of its marginals, i.e. the underlying process is i.i.d. 

The problem of testing independence in time series and time series residuals is fundamental to understanding the stochastic process under study. A frequently used analogue for testing independence is that of testing the hypothesis 

\beq
H_r \! :{\rho} (1) = \dots = {\rho} (r) =0\, \, ,
\label{eq.joint.hyp}
\eeq 

\noindent for some fixed $r$, where $\rho (k)$ is the $k$th-order autocorrelation. Examples of such tests include those proposed by \cite{boxpierce} and \cite{ljungbox}, and the testing procedure proposed by \cite{breusch} and \cite{godfrey}. However, such tests assume that the data-generating model is parametric or semi-parametic and, in particular, follows an ARIMA model. This assumption is, in general, violated for arbitrary $P$, and so these tests will not be exact for finite samples or asymptotically valid, as will be shown later. 

We propose a nonparametric testing procedure for the hypothesis 

\beq
H^{(k)} \! \! : \rho (k) = 0 \,\, ,
\label{eq.piece.hyp}
\eeq

\noindent based on permutation testing, whence we may construct a testing procedure for the hypothesis $H_m$ using multiple testing procedures. Later, we will also consider this joint testing of many autocorrelations simultaneously in a multiple testing framework. 

To review the testing procedure in application to this problem: let $S_n$ be the symmetric, or permutation, group of order $n$. Then, given any test statistic $T_n (X) = T_n(X_1, \, \dots, \, X_n)$, for each element $\pi_n \in S_n$, let $\hat{T}_{\pi_n} = T_n \left( X_{\pi_n (1)}, \, \dots, \, X_{\pi_n (n) }\right)$. Let the ordered values of the $\hat{T}_{\pi_n}$ be 

\beq
\hat{T}_{n} ^{(1)} \leq \, \dots \, \leq \hat{T}_{n} ^{(n!)} \, \, .
\eeq 

\noindent Fix a nominal level $\alpha \in (0, \,1)$, and let $m = n! - \left[ \alpha n!\right]$, where $[x]$ denotes the largest integer less than or equal to $x$. Let $M^+ (x)$ and $M^0 (x)$ be the number of values $\hat{T}_{n} ^{(j)}(x)$ which are greater than and equal to $\hat{T}_{n} ^{(m)} (x)$, respectively. Let

\beq
a(x) = \frac{\alpha n! - M^+ (x)}{M^0 (x)} \, .
\eeq 

\noindent Define the permutation test $\phi (X)$ to be equal to $1$, $a(X)$, or $0$, according to whether $T_n (X)$ is greater than, equal to, or less than $T_n ^{(m)} (X)$, respectively. Additionally, define the permutation distribution 

\beq
\hat{R}_n ^{T_n} (t):= \frac{1}{n!} \sum_{\pi _n \in S_n} I\left\{\hat{T}_{\pi_n} \leq t \right\} \, \, .
\label{eq:perm:defn}
\eeq 

\noindent Let $[n] = \{1, \, \dots,\, n\}$. We observe that, for $\Pi_n \sim \text{Unif}(S_n)$, independent of the sequence $\{X_i, \, i \in [n]\}$, and $X_{\Pi_n} = \left(X_{\Pi_n(1)}, \, \dots, \, X_{\Pi_n (n)}\right)$, the permutation distribution is the distribution of $T_n \left( X_{\Pi_n} \right)$ conditional on the sequence $\{X_i, \, i \in [n] \}$. Also, accounting for discreteness, the permutation test rejects if the observed test statistic $T_n$ exceeds the $1-\alpha$ quantile of the permutation distribution $\hat{R}_n ^{T_n}$. 

Under the randomization hypothesis that the joint distribution of the $X_i$ is invariant under permutation, the permutation test $\phi$ is exact level $\alpha$ (see \cite{tsh}, Theorem 15.2.1), but problems may arise when the null hypothesis $H^{(k)} \! \!: \rho(k) = 0$ holds true, but the sequence $X$ is not independent and identically distributed. Indeed, the distribution of an uncorrelated sequence is not invariant under permutations, and the randomization hypothesis does not hold (the randomization hypothesis guarantees finite-sample validity of the permutation test; see \cite{tsh}, Section 15.2). Such issues may hinder the use of permutation testing for valid inference, but we will show how to restore asymptotic validity to the permutation test.

For instance, consider the problem of testing $H^{(1)} \! \!: \rho (1) =0$, for some sequence $\{X_i, \, i \in [n]\} \sim P$, where $\rho(1) = 0$. If the sequence is not i.i.d., the permutation test may have rejection probability significantly different from the nominal level, which leads to several issues. If the rejection probability is greater than the nominal level, we may reject the null hypothesis, and conclude that there is nonzero first order autocorrelation, whereas in fact we have autocorrelation of some higher order, or some other unobserved dependence structure. A further issue is that of Type 3, or directional, error, in a two-sided test of $H^{(1)}  \!$. In this situation, one runs the risk of rejecting the null and concluding, for instance, that the first-order autocorrelation is larger than 0, when in fact it is less than 0. To illustrate this, if there exists some distribution $P_n$ of the sequence $\{X_i ,\,i \in [n]\}$ with first-order autocorrelation $\rho (1) = 0$ but rejection probability equal to $\gamma \gg \alpha$, by continuity it follows that there exists some distribution $Q_n$ of the sequence with first-order autocorrelation $\rho (1) < 0$, but two-sided rejection probability almost as large as $\gamma$. Under such a distribution, with probability almost $\gamma/2$, not only would we reject the null, but we would also falsely conclude that the first-order sample autocorrelation is greater than 0, when in fact the opposite holds. We will later show that $\gamma$ may be arbitrarily close to 1; see Example \ref{ex:mdep:gaussian} and Remark \ref{rem.counterex}. 
There are also issues if the rejection probability under the null is much smaller than the nominal level. In this case, again by continuity, we would have power significantly less than the nominal level even if the alternative is true, i.e. the test would be biased. 
The strategy to overcome these issues is essentially as follows: assuming stationarity of the sequence $\{X_i,i \geq 1\} $, we wish to show that the permutation distribution based on some test statistic is asymptotically pivotal, i.e. does not depend on the distribution of the $X_i$, in order for the critical region of the associated hypothesis test to not depend on parameters of the distribution of the $X_i$. We then wish to show that the limiting distribution of the test statistic under $H^{(1)} $ is the same as the permutation distribution, so that we may perform (asymptotically) valid inference. Without this matching condition, we may not claim that a permutation test is asymptotically valid, despite being exact under the additional assumption of independence of the sequence.

Significant work has been done on these issues in the context of other problems. \cite{neuhaus1993} discovered the idea of studentizing test statistics to allow for asymptotically valid inference in the permutation testing setting, \cite{janssen1997} compares means by appropriate studentization in a permutation test, \cite{janssen.pauls} give general results about permutation testing, \cite{chung.perm} consider studentizing linear statistics in a two-sample setting, \cite{omelka.pauly} compare correlations by permutation testing, and \cite{diciccio2017} consider testing correlation structure and regression coefficients. In the context of time series data, \cite{perm.test.mri} discuss the application of permutation testing to neuroimaging data, and \cite{short.circadian} consider the application of permutation testing as a method of testing for periodicity in biological data. In a more theoretical setting, \cite{jentsch2015} use randomization methods to test equality of spectral densities, and \cite{ritzwoller} consider permutation testing in the setting of dependent Bernoulli sequences.

The goal of this paper is to provide a framework for the use of permutation testing as a valid method for testing the hypothesis $H^{(k)} \! \!: \rho(k) =0$, which retains the exactness property under the assumption of independence of the $X_i$, but is also asymptotically valid for a large class of weakly dependent stationary sequences. In particular, throughout this paper, we consider the problem of testing $H^{(1)} \! \!  : \rho(1) = 0$, for $\{X_i, \, i \in [n]\}$ a weakly dependent sequence, and with test statistic a possibly studentized version of the sample autocorrelation $\hat{\rho}_n = \hat{\rho}_n(1)$, where

\beq
\hat{\rho}_n (k) \equiv \hat{\rho}_n(X_1, \, \dots, \, X_n; \, k) = 
\frac{\frac{1}{n-k} \sum_{i=1} ^{n-k} \left( X_i - \bar{X}_n \right) \left( X_{i+k} - \bar{X}_n \right)}{\hat{\sigma}_n ^2} \, \, .
\eeq

\noindent $\hat{\sigma}_n ^2 $ is the sample variance, given by

\beq
\hat{\sigma}_n ^2 = \frac{1}{n}\sum_{i=1} ^n (X_i - \bar{X}_n)^2 \, \, ,
\eeq

\noindent and $\bar{X}_n = \frac{1}{n} \sum_{i=1} ^n X_i$. We note that $\bar{X}_n$ and $\hat{\sigma}_n ^2$ are permutation invariant. Unless otherwise stated, we consider the problem of testing $H^{(1)} \! \!: \rho_1 = 0$, where $\rho_1$ is the first-order autocorrelation, and $\hat{\rho}_n$ refers to the first-order sample autocorrelation.

There are several different notions of weak dependence (see \cite{bradley2005} for a discussion thereof). Throughout this paper, we focus on the notions of $m$-dependence and $\alpha$-mixing.

The main results are given in Section \ref{sec.mdep} and \ref{sec.alpha}. In Section \ref{sec.mdep}, we give conditions for the asymptotic validity of the permutation test when $\{X_i, \, i \in [n]\}$ is a stationary, $m$-dependent sequence. In Section \ref{sec.alpha}, under slightly stronger moment assumptions, we extend the result of Section \ref{sec.mdep} to a much larger class of $\alpha$-mixing sequences, which includes a class of stationary ARMA processes. The technical arguments for Sections \ref{sec.mdep} and \ref{sec.alpha} are rather distinct, though the results in both sections allow one to construct valid permutation tests of correlations by appropriate studentization. Section \ref{sec.multiple} provides a framework for using individual permutation tests for different order autocorrelations in a multiple testing setting. Section \ref{sec.simulations} provides simulations illustrating the results. Section \ref{sec.fin} gives an application of the testing procedure to financial data. Section \ref{sec.cov} provides analogous results for testing the equivalent null hypothesis that the first-order autocovariance is equal to zero. The proofs are quite lengthy due to the technical requirements needed to prove the results; consequently, all proofs are deferred to the supplement.

\section{Permutation distribution for $m$-dependent sequences}\label{sec.mdep}

In this section, we consider the problem of testing the null hypothesis

\beq
H^{(1)} \! \! : \rho_1 = 0 \, \, ,
\eeq 

\noindent where $\rho_1 = \rho(\mathbf{X}; \, 1)$ is the first-order autocorrelation, in the setting where the sequence $\{X_i, \, i \in [n]\}$ is stationary and $m$-dependent, i.e. there exists $m \in \nn$ such that, for all $j \in \nn$, the sequences $\{X_i, \, i \in [j]\}$ and $\{X_i, \, i \geq j + m + 1 \}$ are independent. A special case of the $m$-dependence condition is that of $m= 0$, which corresponds to independence of realizations. When the distribution of $(X_1, \, \dots, \, X_n)$ is invariant under permutation, i.e. the sequence is exchangeable, the randomization hypothesis holds, and so one may construct permutation tests of the hypothesis ${H}_0$ with exact level $\alpha$. Note that in the case of $m$-dependence, exchangeability and independence are equivalent conditions\footnote{A proof of this statement is given in Lemma S.3.1 of the supplement.}. 
However, if the realizations of the sequence are not independent, the test may not be valid even asymptotically, i.e. the rejection probability of such a test need not be $\alpha$ for finite samples or even near $\alpha$ in the limit as $n \to \infty$. Hence the goal is to construct a testing procedure, based on some appropriately chosen test statistic, which has asymptotic rejection probability equal to $\alpha$, but which also retains the finite sample exactness property under the assumption of independence of the $X_i$. It is therefore important to analyze the asymptotic properties of the permutation distribution. We assume that the sequence of random variables $\{X_i, \, i \in [n]\}$ is strictly stationary.

We wish to consider a permutation test based on the first-order sample autocorrelation, $\hat{\rho}_n$. Our strategy is as follows: in order to determine the limiting behavior of the permutation distribution, $\hat{R}_n$, we apply Hoeffding's condition (see \cite{tsh}, Theorem 15.2.3). This condition requires that we derive the joint limiting distribution of the normalized first-order sample autocorrelation of the sequence under the action of two independent random permutations. More precisely, we consider the first-order sample autocorrelations of $X_{\Pi_n(1)}, \, \dots, \, X_{\Pi_n(n)}$ and  $X_{\Pi_n ' (1)}, \, \dots, \, X_{\Pi_n '(n)}$, where $\Pi_n$ and $\Pi_n ' $ are independent random permutations of $\{1, \, \dots, \, n\}$, each of which is independent of the sequence $\{X_i, \, i \in [n]\}$. We aim to show that the limiting joint distribution is that of two i.i.d. random variables, each having the limiting distribution of the first-order sample autocorrelation when observations are i.i.d., with the same marginal distribution as the underlying sequence.

To this end, a natural approach in this problem is to use Stein's method. Indeed, we begin by specializing a result of \cite{stein1972} to the case of a sum of random variables whose dependency graph has uniformly bounded degree.

\begin{thm} Let $n \in \mathbb{N}$. Let $X_1, \, \dots, \, X_n$ be random variables such that $\ee X_i = 0$ for all $i$, and, uniformly in $i$, 

\beq
\begin{aligned}
\ee\left[ X_i ^2\right] &\leq M_2 \, \, , \\
\ee \left[X_i ^4 \right]&\leq M_4\,\,.
\end{aligned}
\label{eq:moments.stein}
\eeq

\noindent Let $S_i = \left\{ j \in [n] : X_i \text{ and } X_j \text{ are not independent}\right\}$. Suppose $\left| S_i \right| \leq D < \infty$ for all $i$. Let 

\beq
\sigma_n ^2 = \ee \left[ \sum_{i=1} ^n X_i \sum_{j \in S_i} X_j \right] = \var\left( \sum_{i=1} ^n X_i \right) \, \, ,
\eeq

\noindent and let $W_n = \sum_{i=1} ^n X_i/\sigma_n$. Then, for all $t \in \mathbb{R}$,

\beq
\left| \pp \left( W_n \leq t \right) - \Phi (t) \right| \leq \frac{1}{\sigma_n} \left( 4 \left( \frac{nD^3 M_4}{\sigma_n^2}\right)^{1/2} + 2^{3/4} \pi^{-1/4} \left(\frac{n}{\sigma_n} \right)^{1/2} \left( D^5 M_2 M_4\right)^{1/4} \right) \, \, .
\label{eq.clt}
\eeq

\label{thm.clt}
\end{thm}

\noindent We note several consequences of this result:

\begin{rem}\rm\label{rem:rate.clt.stein}
The bound on the right hand side of (\ref{eq.clt}) is independent of $t$. $\qed$
\label{rem.unif.clt}
\end{rem}
\begin{rem}\rm
If $X_i = \tilde{X}_i /\sqrt{n}$, for some $\tilde{X}_i$ with bounded 4th moments, and $\sigma_n ^2 \asymp 1$, then the right hand side of (\ref{eq.clt}) is $O(n^{-1/4})$, i.e. this result provides a CLT. Note also that if, instead of the moment condition (\ref{eq:moments.stein}), we have that the sequence $\{X_i, \, i \in [n] \}$ is uniformly bounded, we may apply a result of \cite{rinott} and instead replace the right hand side of (\ref{eq.clt}) with $O(n^{-1/2})$. $\qed$
\label{rem.clt.clt}
\end{rem}

We may now use the result of Theorem \ref{thm.clt} to exhibit the asymptotic properties of the permutation distribution based on $\sqrt{n} \hat{\rho}_n$.

\begin{thm} Let $\{X_i, \, i \geq 1\}$ be an $m$-dependent stationary time series, with finite 8th moment, or an i.i.d. sequence with finite 4th moment. The permutation distribution, $\hat{R}_n$, as defined in (\ref{eq:perm:defn}), of $\sqrt{n} \hat{\rho}_n$, based on the test statistic $\hat{\rho}_n = \hat{\rho}(X_1, \, \dots, \, X_n)$, with associated group of transformations $S_n$, the symmetric group of order $n$, satisfies

\beq
\sup_{t \in \rr} \left| \hat{R}_n (t)-\Phi(t)\right| \overset{p}{\to} 0 \, \, ,
\eeq

\noindent as $n \to \infty$, where $\Phi$ is the standard Gaussian c.d.f.

\label{rand.dist.mdep}
\end{thm}

\begin{rem}\rm
In this result, observe that, under independence, the moment conditions required for asymptotic normality of the permutation distribution are weaker than under general $m$-dependence, since, in the case of $m$-dependence for arbitrary $m$, we have the additional requirement of finiteness of the variance of products of the $X_i$. $\qed$
\end{rem}

\noindent We have shown that the permutation distribution is asymptotically Gaussian, with mean and variance not depending on the underlying process, and that the result holds irrespective of whether or not the null hypothesis $H^{(1)}$ holds. This result may be interpreted as follows. Under the action of a random permutation, for large values of $n$, one would expect that the first-order sample autocorrelation of the sequence $\{X_n, \, n\in \nn \}$ behaves similarly to the case of $X_i \overset{i.i.d.}\sim F$, where $F$ is the marginal distribution of the $X_i$, since the dependence between consecutive terms in the permuted sequence will be very weak, on account of the large sample size and the localized dependence structure of the original sequence. However, the same is not true of the asymptotic distribution of the test statistic. Indeed, under the null hypothesis, the asymptotic distribution of $\sqrt{n} \hat{\rho}_n$ is also Gaussian with mean 0, but with variance not necessarily equal to 1. Therefore it is not possible to claim asymptotic validity of the permutation test based on this test statistic.

\begin{thm} Let $X_1, \, \dots, \, X_n$ be a strictly stationary sequence, with variance $\sigma^2 >0$, and first-order autocorrelation $\rho_1$, such that one of the the following two conditions holds. 

\begin{enumerate}

\item[i)] $\{X_i, \, i \in [n]\}$ is $m$-dependent, for some $m \in \nn$, and $\ee\left[X_1 ^4\right] < \infty$.

\item[ii)] $\{X_i, \, i \in [n]\}$ is $\alpha$-mixing, and, for some $\delta > 0$, we have that 

\beq
\ee\left[\left|X_1 \right|^{4 + 2 \delta}\right] < \infty \, \, ,
\eeq

\noindent and the $\alpha$-mixing coefficients $\alpha_X(\cdot)$ satisfy

\beq
\sum_{n \geq 1} \alpha_X(n)^\frac{\delta}{2+ \delta} < \infty \, \, .
\eeq

\end{enumerate}

\noindent Let $\hat{\rho}_n$ be the sample first-order autocorrelation. Let 

\beq
\begin{aligned}
\kappa^2 &= \var \left(X_1 ^2 \right) + 2\sum_{k \geq 2} \cov \left( X_1 ^2, \, X_k ^2 \right) \\
\tau_1 ^2 &= \var \left(X_1X_2\right) + 2\sum_{k \geq 2} \cov (X_1X_2, \, X_k X_{k+1} )\\
\nu_1 &= \cov \left( X_1 X_2, \, X_1 ^2 \right) + \sum_{k \geq 2} \cov\left( X_1^2, \, X_{k} X_{k+1} \right) + \sum_{k \geq 2} \cov \left(X_1 X_2,\, X_k ^2 \right) \, \, .
\end{aligned}
\eeq

\noindent Let 

\beq
\gamma_1^2 = \frac{1}{\sigma^4} \left( \tau_1^2  - 2 \rho_1 \nu_1 + \rho_1 ^2 \kappa^2 \right) \, \, .
\label{eq.gam.defn}
\eeq

\noindent Suppose that $\kappa^2, \, \tau_1 ^2, \, \gamma_1 ^2 \in (0, \, \infty)$. Then, as $n \to \infty$, 

\beq
\sqrt{n} \left(\hat{\rho}_n - \rho_1 \right) \overset{d} \to N\left(0, \, \gamma_1^2 \right) \, \, .
\eeq

\label{thm.unperm.rho.lim}
\end{thm}

Since, clearly, $\gamma_1 ^2 = 1$ does not hold in general, a permutation test based on the test statistic $\sqrt{n} \hat{\rho}_n$ will not be asymptotically valid. However, note that, under the additional restriction of independence, $\gamma_1^2 =1$ always, hence, as is consistent with the test being exact under independence, the permutation test will also be asymptotically valid in this case. One could conclude that the permutation test based on the above test statistic is not asymptotically valid in general, and attempt to find a different test statistic, for which the permutation distribution and test statistic distribution are asymptotically the same. 

Alternatively, one could adapt the test statistic above in some fashion, in order to resolve the issue of (asymptotically) mismatched variances in the permutation distribution and distribution of the test statistic. In particular, a natural way to adapt the test statistic is to studentize it by some estimator of $\gamma_1 ^2$, motivated by the heuristics that, under permutations, all dependence structure in the sequence will be broken, and the estimator will be approximately equal to 1. Therefore, despite the limiting distribution of $\sqrt{n} \hat{\rho}_n$ being different, in general, from the case when $X_i\overset{i.i.d.}\sim F$, where $F$ is the marginal distribution of the $X_i$, under appropriate studentization, the limiting behaviors will be the same.

To this end, we now consider a permutation test based on some studentized version of the test statistic $\sqrt{n} \hat{\rho}_n$. Provided we can find a weakly consistent estimator $\hat{\gamma}_n ^2 = \hat{\gamma}_n ^2 \left( X_1, \, \dots, \, X_n \right)$ of $\gamma_1^2$, such that, for $\Pi_n$ a random permutation independent of the sequence $\{X_i, \, i \in [n]\}$, we also have that $\hat{\gamma}_n^2\left(X_{\Pi_n (1)},\, \dots, \, X_{\Pi_n (n)} \right) = \var\left(X_1 \right)^2 + o_p(1)$, we may apply Slutsky's theorem for randomization distributions (\cite{chung.perm}, Theorem 5.2) to studentize the test statistic and construct an asymptotically valid permutation test. Combining this with an application of Ibragimov's central limit theorem for $\alpha$-mixing random variables (\cite{ibragimov}), and noting that stationary $m$-dependent sequences necessarily satisfy the mixing conditions laid out therein, we have the following result.

\begin{thm} Let $m \in \nn$. Let $X_1, \, \dots, \, X_n$ be a strictly stationary, $m$-dependent sequence, with variance $\sigma^2 >0$, first-order autocorrelation $\rho_1$, and finite 8th moment. Let $\hat{\sigma}_n ^2$ be the sample variance. Let $\kappa^2, \, \tau_1^2$, $\nu_1$ and $\gamma_1^2$ be as in Theorem \ref{thm.unperm.rho.lim}. Suppose that $\kappa^2, \, \tau_1 ^2, \, \gamma_1 ^2 \in (0, \, \infty)$. For $i \in \nn$, let $Y_i = \left(X_i  - \bar{X}_n \right) \left( X_{i+1}  - \bar{X}_n\right)$, and let $Z_i = \left(X_i  - \bar{X}_n \right)^2$. Let $b_n = o\left(\sqrt{n}\right)$ be such that, for all $n$ sufficiently large, $b_n \geq m+ 2$. Let

\beq
\begin{aligned}
\hat{K}_n ^2 &= \frac{1}{n} \sum_{i=1} ^n \left(Z_i - \bar{Z}_n \right)^2 +  \frac{2}{n}\sum_{j=1}^{b_n} \sum_{i=1}^{n-j} \left( Z_i - \bar{Z}_n\right)\left(Z_{i+j} - \bar{Z}_n\right) \\
\hat{T}_n ^2 &= \frac{1}{n} \sum_{i=1} ^{n-1} \left(Y_i - \bar{Y}_n \right)^2 + \frac{2}{n}\sum_{j=1}^{b_n}\sum_{i=1} ^{n-j-1} \left(Y_i - \bar{Y}_n \right) \left(Y_{i+j} -\bar{Y}_n \right) \\
\hat{\nu}_n &=\frac{1}{n} \sum_{i=1}^{n-1} \left(Y_i - \bar{Y}_n\right)\left(Z_i -\bar{Z}_n \right) + \frac{1}{n} \sum_{j=1}^{b_n} \sum_{i=1} ^{n- j-1} \left(Z_i - \bar{Z}_n \right)\left(Y_{i + j} - \bar{Y}_n\right)  +\\ &+\frac{1}{n} \sum_{j=1} ^{b_n} \sum_{i=1} ^{n-j} \left(Y_i - \bar{Y}_n \right)\left(Z_{i+j} - \bar{Z}_n\right) \,\,.
\end{aligned}
\label{mdep.var.est}
\eeq  

\noindent Let 

\beq 
\hat{\gamma}_n ^2 = \frac{1}{\hat{\sigma}_n ^4 } \left[ \hat{T}_n^2 - \hat{\rho}_n \hat{\nu}_n + \hat{\rho}_n ^2 \hat{K}_n ^2 \right] \, \, .
\label{eq.gam.var.defn}
\eeq

\begin{enumerate} 

\item[i)] As $n \to \infty$, 

\beq
 \frac{\sqrt{n}\left( \hat{\rho}_n -  \rho \right)}{\hat{\gamma}_n}\overset{d} \to N \left(0, \, 1 \right) \, \, .
\label{clt.student}
\eeq

\item[ii)] Let $\hat{R}_n$ be the permutation distribution, with associated group of transformations $S_n$, the symmetric group of order $n$, based on the test statistic $\sqrt{n} \hat{\rho}_n/ \hat{\gamma}_n$. Then, as $n \to \infty$,

\beq
\sup_{t \in \rr} \left| \hat{R}_n (t) - \Phi (t) \right| \overset{p} \to 0 \, \, ,
\label{rand.dist.student}
\eeq

\noindent where $\Phi$ is the standard Gaussian c.d.f.

\end{enumerate}

\label{rand.dist.mdep.stud}
\end{thm}

\begin{rem} \rm
Under the assumptions set out in Theorem \ref{rand.dist.mdep.stud}, in particular as a result of (\ref{clt.student}), the level $\alpha$ permutation test of the null ${H}^{(1)}\! \!: \rho_1 = 0$ based on the test statistic $\sqrt{n} \hat{\rho}_n/\hat{\gamma}_n$ is asymptotically valid. $\qed$
\end{rem}

\begin{rem}\rm If the dependence parameter $m$ is known, one may replace the upper limit $b_n$ in (\ref{mdep.var.est}) with $m +2$, since, for all $k > m+2$,

\beq
\begin{aligned}
\cov\left(X_1 X_2, \, X_{k +1} X_{k +2} \right) &= 0\\
\cov\left(X_1 ^2 , \, X_k^2 \right) &= 0 \\
\cov\left(X_1 ^2 , \, X_i X_{k+1}\right) &= 0 \\
\cov\left(X_1 X_2 , \, X_k ^2 \right) &= 0 \,\,.
\end{aligned}
\eeq

\noindent However, the construction provided in Theorem \ref{rand.dist.mdep.stud} does not require knowledge of $m$. In general, we require that $b_n$ is sufficiently large, in order to guarantee convergence of the estimator of $\hat{\gamma}_n ^2$. $\qed$
\end{rem}

\begin{example}\rm\label{ex:mdep:gaussian} (\emph{Products of i.i.d. random variables}) Let $\{Z_n, \, n \in \nn\}$ be mean zero, i.i.d., non-constant random variables, such that 

\beq
\ee \left[ Z_1 ^8 \right] < \infty \, \, .
\eeq

\noindent Fix $m \in \nn$, and, for each $i$, let

\beq
X_i = \prod_{j = i} ^{i + m -1} Z_j \, \, .
\eeq

\noindent We observe that the sequence $\{X_i, \, i \geq 1\} $ is stationary and $m$-dependent, and that, by Fubini's theorem, the $X_i$ have uniformly bounded 8th moments. It now suffices to show that $\kappa^2$ and $\tau_1^2$ are finite and strictly greater than 0, and $\gamma_1^2$, as defined in (\ref{eq.gam.defn}), is finite and strictly greater than zero. Let $M_k$ be the $k$th moment of $Z_1$, $k \geq 1$. Simple calculations show that

\beq
\begin{aligned}
\var\left(X_1 X_2 \right) &= M_2 ^2 M_4 ^{m-1} \\
\cov\left(X_1 X_2, \, X_k X_{k+1} \right) &= 0 \, \, , \, k \geq 2 \, \, .
\end{aligned}
\eeq

\noindent Hence

$$
\tau_1 ^2 = \frac{M_4 ^{m-1}}{M_2 ^{2(m-1)}} \in  (0, \, \infty)\, \, .
$$

\noindent Additionally, we have that 

$$
\begin{aligned}
\kappa^2 &= M_4 ^m - M_2^{2m} \in (0, \, \infty) \\
\nu_1 &= 0 \, \, .
\end{aligned}
$$

\noindent Hence we have that $\gamma_1^2 = M_4 ^{m-1} \in (0,\, \infty)$. It follows that we may apply the result of Theorem \ref{rand.dist.mdep.stud}, and conclude that the rejection probability of the permutation test based on the test statistic $\sqrt{n} \hat{\rho}_n/\hat{\gamma}_n $ converges to $\alpha$ as $n \to \infty$.

\label{ex.mdep.prod}
\end{example}

\begin{rem} \rm Example \ref{ex.mdep.prod} also provides an illustrative example of the need for studentization in the permutation test. Indeed, in the setting of Example \ref{ex.mdep.prod}, for $r \in \nn$ odd, let

\beq
Z_i = G_i ^r \, \, ,
\label{mdep.counterex}
\eeq

\noindent where $\{G_i, \, i \in \zz\}$ are independent standard Gaussian random variables. Hence

$$
\gamma_1 ^2 = \left(\frac{\ee\left[ G_i ^{4r} \right]}{\ee\left[ G_i ^{4r} \right]^2}\right) ^{m-1} = \left(\frac{(4r - 1)!! }{((2r-1)!!)^2}\right)^{m-1}\, \, .
$$

\noindent By Theorems and \ref{rand.dist.mdep} and \ref{thm.unperm.rho.lim}, it follows that the asymptotic rejection probability of the level $\alpha$ two-sided permutation test, based on the test statistic $\sqrt{n}\hat{\rho}_n$, converges to 

$$
2\left(1- \Phi\left( \left(\frac{((2r-1)!!)^2}{(4r- 1)!!}\right)^\frac{m-1}{2} z _{1- \alpha/2} \right) \right) \, \, ,
$$

\noindent as $n \to \infty$, where $z_{1- \alpha/2}$ is the $\alpha/2$ quantile of the standard normal distribution. It follows that this rejection probability can be arbitrarily close to 1 for large values of $n$, $m$ and $r$. Therefore, by continuity, there exists a distribution $Q_n$ of $(X_1, \, \dots, \, X_n)$ such that $\rho(1) < 0$, but the two-sided permutation test based on $\sqrt{n}\hat{\rho}_n$ would reject $H^{(1)}$\! \!\!, with probability arbitrarily close to $1/2$, and conclude that the first-order sample autocorrelation is greater than zero, when, in fact, the opposite is true. \qed
 
\label{rem.counterex}
\end{rem}

We have shown that a permutation test of the hypothesis $H^{(1)} \! \!: \rho_1 = 0$ is asymptotically valid under assumptions of $m$-dependence, with the permutation distribution converging in probability to the standard Gaussian distribution. 

In Section \ref{sec.alpha}, we extend the results of this section to a much richer class of time series, such as ARMA processes, and processes for which there can be dependence between $X_i$ and $X_j$ for arbitrarily large values of $\left| i - j \right|$. We extend these results by imposing a small additional constraint on the moments of the sequence, and by imposing fairly standard assumptions on the mixing coefficients of the underlying process.

\section{Permutation distribution for $\alpha$-mixing sequences}\label{sec.alpha}

In order to extend the results of Section \ref{sec.mdep} to the broader setting of $\alpha$-mixing sequences, we will show that an appropriately studentized version of the first-order sample autocorrelation has permutation distribution asymptotically not depending on the underlying process $\{X_n, \, n \in \nn \}$, and that, under the null hypothesis $H^{(1)} \! \!: \rho_1 = 0$, the test statistic has asymptotic distribution equal to that of the permutation distribution. As a review, let $\{X_n, \, n \in \zz \}$ be a stationary sequence of random variables, adapted to the filtration $\{\mathcal{F}_n\}$. Let

\beq
\mathcal{G}_n = \sigma\left(X_r : \, r \geq n\right) \, \, .
\eeq 

\noindent For $n \in \nn$, let $\alpha_X(n)$ be Rosenblatt's $\alpha$-mixing coefficient, defined as 

\beq
\alpha_X (n) = \sup_{A \in \mathcal{F}_0, \, B \in \mathcal{G}_n} \left| \pp(A \cap B) - \pp(A) \pp(B)  \right | \, \, .
\eeq

\noindent We say that $\{X_n\}$ is $\alpha$-mixing if $\alpha_X (n) \to 0$ as $n \to \infty$. 

Note that, analogously to the discussion in Section \ref{sec.mdep}, in the setting of $\alpha$-mixing sequences, all exchangeable sequences are also independent\footnote{A proof of this statement is given in Lemma S.3.1 of the supplement.}, i.e. any such testing procedure will retain the exactness property under the additional assumption of independence of the $X_i$, and this is the only condition under which the randomization hypothesis holds.

Unfortunately, however, the method of proof used in Section \ref{sec.mdep} can no longer apply, since the dependency graph of an arbitrary $\alpha$-mixing sequence $\{X_n, \, n \in \nn\}$ has infinite degree, and so we cannot apply Theorem \ref{thm.clt}. We proceed instead as follows, in the spirit of \cite{noether}. Suppose, for now, that the sequence $\{X_n, \, n \in \nn\}$ is uniformly bounded. For $\Pi_n \sim \text{Unif}(S_n)$, observing that the permutation distribution based on some test statistic $T_n \left(X_1, \, \dots, \, X_n \right)$ is the empirical distribution of $T_n \left( X_{\Pi_n(1)} ,\dots, \, X_{\Pi_n(n)} \right)$ {conditional} on the data $\{X_i, \, i \in [n]\}$, we may condition on the data and apply the central limit theorem of \cite{wald1943}, checking that appropriate conditions on the sample variance are satisfied. This allows us to obtain a convergence result for a distribution very closely related to that of the permutation distribution, but with additional centering and scaling factors.

We are now in a position to use a double application of Slutsky's theorem for randomization distributions, in order to remove the centering and scaling factors, thus obtaining the following result.

\begin{thm} Let $\{X_n, \, n\in \nn\}$ be a stationary, bounded, $\alpha$-mixing sequence. Suppose that 

\beq
\sum_{n \geq 1} \alpha_X (n) < \infty \, \,.
\eeq

\noindent The permutation distribution of $\sqrt{n} \hat{\rho}_n$ based on the test statistic $\hat{\rho}_n = \hat{\rho}(X_1, \, \dots, \, X_n)$, with associated group of transformations $S_n$, the symmetric group of order $n$, satisfies

\beq
\sup_{t \in \rr} \left| \hat{R}_n (t)-\Phi(t)\right| \overset{p}{\to} 0 \, \, ,
\eeq

\noindent as n $\to \infty$, where $\Phi(t)$ is the distribution of a standard Gaussian random variable.

\label{alpha.bd.perm}
\end{thm}

We now wish to remove the boundedness constraint of Theorem \ref{alpha.bd.perm} and extend its result to the setting of stationary, $\alpha$-mixing sequences with uniformly bounded moments of some order. In order to do this, we require the following lemma.

\begin{lem} For each $N, \, n \in \nn$, let $G_{N, \, n}: \rr \to \rr $ be an nondecreasing random function such that, for each $N$, and for all $t \in \rr$, as $n \to \infty$,

\beq
G_{N, \, n} (t) \overset{p}\to g_N(t) \, \, ,
\eeq

\noindent where, for each $N$, $g_N : \rr \to \rr$ is a function. Suppose further that, as $N \to \infty$, for each $t \in \rr$,

\beq
g_N (t) \to g(t) \, \, ,
\eeq

\noindent where $g: \rr \to \rr$ is continuous. Then, there exists a sequence $N_n \to \infty$ such that, for all $t \in \rr$, as $n \to \infty$, 

\beq
G_{N_n, \, n} (t) \overset{p} \to g(t) \, \, .
\eeq

\label{trunc.diag.lem}
\end{lem}

We may now proceed to extend Theorem \ref{alpha.bd.perm} as follows. Let $\{X_i, \, i \geq 1\}$ be an $\alpha$-mixing sequence, with summable $\alpha$-mixing coefficients, and let $G_{N, \, n} (t)$ be the permutation distribution, evaluated at $t$, of the truncated sequence $Y_i = \left(X_i \wedge N \right) \vee (-N)$, based on the test statistic $\sqrt{n}\hat{\rho}_n$. Let $g_N = g = \Phi$, where $\Phi$ is the distribution of a standard Gaussian random variable. By Theorem \ref{alpha.bd.perm}, the conditions of Lemma \ref{trunc.diag.lem} are satisfied, so we apply Lemma \ref{trunc.diag.lem} in order to find an appropriate sequence of truncation parameters $N_n$.

Then, for $\Pi_n \sim \text{Unif}(S_n)$, and $Y_i = \left(X_i \wedge N_n \right) \vee (-N_n)$, we relate the first-order sample autocorrelation $\hat{\rho}_n \left( X_{\Pi_n(1)}, \, \dots, \, X_{\Pi_n(n)} \right)$ to the truncated first-order sample autocorrelation $\hat{\rho}_n \left( Y_{\Pi_n(1)}, \, \dots, \, Y_{\Pi_n(n)} \right)$. Bounding the difference of these two autocorrelations in probability using \cite{doukhan}, Section 1.2.2, Theorem 3, and applying Slutsky's theorem for randomization distributions once more, we obtain the following result.

\begin{thm} Let $\{X_n, \, n \in \nn\}$ be a stationary, $\alpha$-mixing sequence, with mean $0$ and variance $1$. Suppose that, for some $\delta > 0$, 

\beq
\sum_{n \geq 1} \alpha_X (n) ^{\frac{\delta}{2 + \delta}} < \infty \, \, ,
\eeq

\noindent and 

\beq
\ee \left[ \left|X_1\right| ^{8 + 4 \delta} \right] < \infty \, \,.
\eeq

\noindent Then, the permutation distribution of $\sqrt{n} \hat{\rho}_n$ based on the test statistic $\hat{\rho}_n = \hat{\rho}_n(X_1, \, \dots, \, X_n)$, with associated group of transformations $S_n$, the symmetric group of order $n$, satisfies

\beq
\sup_{t \in \rr} \left| \hat{R}_n (t) - \Phi(t) \right|  \overset{p} \to 0 \, \, ,
\eeq

\noindent where $\Phi(t)$ is the distribution of a standard Gaussian random variable.

\label{alpha.perm}
\end{thm}

As in the case of $m$-dependence, despite the asymptotic normality of the permutation distribution, we may still not, in general, use a permutation test in this setting, since the asymptotic distribution of the test statistic under the null may not be the same as the permutation distribution. To that end, we again consider studentizing the test statistic $\sqrt{n} \hat{\rho}_n$ by an appropriate estimator of its standard deviation.

\begin{lem} Let $\{X_n, \, n \in \nn\}$ be a stationary, $\alpha$-mixing sequence such that, for some $\delta >0$, 

\beq
\ee \left[ \left| X_1\right| ^{8 + 4\delta}\right] < \infty \, \, ,
\eeq

\noindent and 

\beq
\sum_{n \geq 1} \alpha_X (n) ^{\frac{\delta}{2 + \delta} }< \infty  \, \, .
\eeq

\noindent Let 

\beq
\begin{aligned}
\kappa^2 &= \var \left(X_1 ^2 \right) + 2\sum_{k \geq 2} \cov \left( X_1 ^2, \, X_k ^2 \right) \\
\tau_1 ^2 &= \var \left(X_1X_2\right) + 2\sum_{k \geq 2} \cov (X_1X_2, \, X_k X_{k+1} )\\
\nu_1 &= \cov \left( X_1 X_2, \, X_1 ^2 \right) + \sum_{k \geq 2} \cov\left( X_1^2, \, X_{k} X_{k+1} \right) + \sum_{k \geq 2} \cov \left(X_1 X_2,\, X_k ^2 \right) \, \, .
\end{aligned}
\eeq

\noindent Let

\beq
\gamma_1^2 = \frac{1}{\sigma^4} \left( \tau_1^2  - 2 \rho_1 \nu_1 + \rho_1 ^2 \kappa^2 \right) \, \, .
\eeq

\noindent  Suppose that $\kappa^2, \, \tau_1 ^2, \, \gamma_1 ^2 \in (0, \, \infty)$. Let $b_n = o\left(\sqrt{n}\right)$ be such that $b_n \to \infty$ as $n \to \infty$. Let $\hat{K}_n ^2$, $\hat{T}_n ^2$, and $\hat{\nu}_n$ be as in (\ref{mdep.var.est}), and let $\hat{\gamma}_n^2$ be as in (\ref{eq.gam.var.defn}). Then, as $n \to \infty$, 

\beq
\hat{\gamma}_n ^2 \overset{p} \to \gamma_1^2 \, \, .
\eeq

\label{lem.stud.unperm.pap}
\end{lem}

\begin{lem} In the setting of Lemma \ref{lem.stud.unperm.pap}, let $\Pi_n \sim \text{\normalfont Unif}(S_n)$, independent of the sequence $\{X_n, \, n\in \nn\}$. For $i \in \nn$, let

\beq
\begin{aligned}
\tilde{Y}_i &= \left(X_{\Pi_n (i)} - \bar{X}_n \right) \left( X_{\Pi_n(i+1)} -\bar{X}_n \right) \\
\tilde{Z}_i &= \left(X_{\Pi_n (i)} - \bar{X}_n \right)^2 \, \, .
\end{aligned}
\eeq

\noindent Let

\beq
\begin{aligned}
\hat{K}_n ^2 &= \frac{1}{n} \sum_{i=1} ^n \left(\tilde{Z}_i - \bar{\tilde{Z}}_n \right)^2 +  \frac{2}{n}\sum_{j=1}^{b_n} \sum_{i=1}^{n-j} \left( \tilde{Z}_i - \bar{\tilde{Z}}_n\right)\left(\tilde{Z}_{i+j} - \bar{\tilde{Z}}_n\right) \\
\hat{T}_n ^2 &= \frac{1}{n} \sum_{i=1} ^{n-1} \left(\tilde{Y}_i - \bar{\tilde{Y}}_n \right)^2 + \frac{2}{n}\sum_{j=1}^{b_n}\sum_{i=1} ^{n-j-1} \left(\tilde{Y}_i - \bar{\tilde{Y}}_n \right) \left(\tilde{Y}_{i+j} -\bar{\tilde{Y}}_n \right) \\
\hat{\nu}_n &=\frac{1}{n} \sum_{i=1}^{n-1} \left(\tilde{Y}_i - \bar{\tilde{Y}}_n\right)\left(\tilde{Z}_i -\bar{\tilde{Z}}_n \right) + \frac{1}{n} \sum_{j=1}^{b_n} \sum_{i=1} ^{n- j-1} \left(\tilde{Z}_i - \bar{\tilde{Z}}_n \right)\left(\tilde{Y}_{i + j} - \bar{\tilde{Y}}_n\right)  +\\ &+\frac{1}{n} \sum_{j=1} ^{b_n} \sum_{i=1} ^{n-j} \left(\tilde{Y}_i - \bar{\tilde{Y}}_n \right)\left(\tilde{Z}_{i+j} - \bar{\tilde{Z}}_n\right) \,\,.
\end{aligned}
\eeq  

\noindent Let 

\beq 
\hat{\gamma}_n ^2 = \frac{1}{\hat{\sigma}_n ^4 } \left[ \hat{T}_n^2 - \hat{\rho}_n \hat{\nu}_n + \hat{\rho}_n ^2 \hat{K}_n ^2 \right] \, \, .
\eeq

\noindent We have that, as $n \to \infty$, 

\beq
\hat{\gamma}_n ^2 \overset{p} \to 1 \, \, .
\eeq

\label{pap.var.perm.stud}
\end{lem}

With the results of Lemmas \ref{lem.stud.unperm.pap} and \ref{pap.var.perm.stud}, we may once again apply Slutsky's theorem for randomization distributions and conclude the following.

\begin{thm} Let $\{X_n,\, n \in \nn\}$ be a strictly stationary, $\alpha$-mixing sequence, with variance $\sigma^2$ and  first-order autocorrelation $\rho_1$, such that, for some $\delta > 0$,

\beq
\ee \left[ \left| X_1 \right| ^{8 + 4 \delta}\right] < \infty \, \, ,
\label{1d.moment.cond}
\eeq

\noindent and 

\beq
\sum_{n \geq 1} \alpha_X (n) ^\frac{\delta}{2 + \delta} < \infty \, \, .
\label{1d.mixing.cond}
\eeq

\noindent Let $\kappa^2, \, \tau_1^2, \, \nu_1$ and $\gamma_1^2$ be as in Theorem \ref{thm.unperm.rho.lim}. Suppose that $\kappa^2, \, \tau_1 ^2, \, \gamma_1 ^2 \in (0, \, \infty)$. Let $b_n = o\left(\sqrt{n}\right)$ be such that $b_n \to \infty$ as $n \to \infty$. Let $\hat{K}_n ^2$, $\hat{T}_n ^2$, and $\hat{\nu}_n$ be as in (\ref{mdep.var.est}), and let $\hat{\gamma}_n^2$ be as in (\ref{eq.gam.var.defn}).

\noindent \begin{enumerate}

\item[i)] We have that, as $n \to \infty$,

\beq
\frac{\sqrt{n}\left( \hat{\rho}_n - \rho_1 \right)}{\hat{\gamma}_n} \overset{d} \to N \left( 0, \, 1 \right) \, \, .
\label{main.test.stat.1d}
\eeq

\item[ii)] Let $\hat{R}_n$ be the permutation distribution, with associated group of transformations $S_n$, the symmetric group of order $n$, based on the test statistic $\sqrt{n} \hat{\rho}_n/\hat{\gamma}_n$. Then, as $n \to \infty$,

\beq
\sup_{t \in \rr} \left| \hat{R}_n (t)-\Phi(t)\right| \overset{p}{\to} 0 \, \, ,
\label{stud.perm.test.1d}
\eeq

where $\Phi$ is the standard Gaussian c.d.f.

\end{enumerate}

\label{perm.test.big.thm}
\end{thm}

We now illustrate the application of Theorem \ref{perm.test.big.thm} to the class of stationary ARMA processes.

\begin{example}\rm\label{ex:arma}(\emph{ARMA process}) Let $\{X_i,\, i \in \zz \}$ satisfy the equation 

\beq
\sum_{i = 0} ^p B_i X_{t - i} = \sum_{k = 0} ^q A_k \epsilon_k \, \, ,
\label{eq.arma.ex}
\eeq

\noindent where the $\epsilon_k$ are independent and identically distributed, and $\ee \epsilon_k =0$, i.e. $X$ is an ARMA$(p, \, q)$ process. Let $X$ have first-order autocorrelation $\rho = 0$. Let 

\beq
P(z) := \sum_{i=0} ^p B_i z^i \, \, .
\eeq

\noindent If the equation $P(z) = 0$ has no solutions inside the unit circle $\{ z \in \mathbb{C}: \left | z \right| \leq 1\}$, there exists a unique stationary solution to (\ref{eq.arma.ex}). By \cite{mokkadem}, Theorem 1, if the distribution of the $\epsilon_k$ is absolutely continuous with respect to Lebesgue measure on $\rr$, and also that, for some $\delta >0$, 

\beq
\ee\left[ \left |\epsilon_1\right|^ {8 + 4 \delta}\right] > 0\, \,,
\eeq

\noindent we have that the sequence $\{X_i,\, i \in \nn \}$ satisfies the conditions of Theorem \ref{perm.test.big.thm}, as long as $\gamma_1^2$, as defined in (\ref{eq.gam.defn}), is finite and positive. Therefore, asymptotically, the rejection probability of the permutation test applied to such a sequence will be equal to the nominal level $\alpha$.

\end{example}

\begin{example}\rm\label{ex:ar2}(\emph{AR(2) process}) We specialize Example \ref{ex:arma} to the case of an AR(2) process with first-order autocorrelation equal to 0. Suppose that the strictly stationary sequence $\{X_i, \, i \geq 1\}$ satisfies, for all $t >2$, the equation

\beq
X_t = \phi X_{t-1} + \rho X_{t-2} + \epsilon_t \, \,,
\label{eq:ar2}
\eeq

\noindent where the $\epsilon_t$ are as in Example \ref{ex:arma}. The first-order autocorrelation of $X$ is given by

\beq
\rho(1) = \frac{\phi}{1-\rho} \, \, .
\eeq

\noindent Hence, for $X$ to be such that $\rho(1) = 0$, we must have that $\phi = 0$. In particular, it follows that $\{X_{2i}, \, i\geq 1\}$ and $\{X_{2i - 1}, \, i \geq 1\}$ are independent and identically distributed stationary AR(1) processes with parameter $\rho$. By the same argument as in Example \ref{ex:arma}, the requisite $\alpha$-mixing condition is satisfied, and so, in order for the result of Theorem \ref{perm.test.big.thm} to apply, it suffices to show that $\tau_1^2, \, \kappa^2$, and $\nu_1$, as defined in (\ref{mdep.var.est}), are finite and nonzero. If so, since $\rho(1) = 0$, we have that $\gamma_1 ^2 = \tau_1 ^2/\sigma^4$, and so the variance condition on $\gamma_1^2$ is automatically satisfied. We begin by noting that, for $i$ odd,

\beq
\cov \left(X_1, \, X_i\right)  = \rho^\frac{i-1}{2} \var(X_1) \, \,, 
\eeq

\noindent and similarly for the covariance between $X_2$ and $X_i$, for $i$ even. Also, note that $\ee\left[X_1 \right] =0$. Simple calculations show that 

\beq
\begin{aligned}
\var\left(X_1 X_2\right) &= \var \left(X_1 \right)^2 \\
\cov \left( X_1 X_2, \, X_k X_{k+1} \right) &= \rho^{k-1} \var \left(X_1 \right)^2 \, \,, \, k \geq 2 \, \, .
\end{aligned}
\eeq

\noindent It follows that $\tau_1^2 \in (0, \, \infty)$. Similarly, we have that 

\beq
\begin{aligned}
\nu_1 &= 0 \, \, ,
\end{aligned}
\eeq

\noindent and we have that 

\beq
\cov\left(X_1 ^2, \, X_k ^2 \right) = \begin{cases}
\rho^{k-1} \var\left(X_1 ^4 \right) \, \, , \, \text{ if } k \in 2 \nn \, \, ,\\
0 \, \, , \, \text{ otherwise.}
\end{cases}
\eeq

\noindent It follows that $\gamma_1 ^2 \in (0, \, \infty)$, and so the result of Theorem \ref{perm.test.big.thm} holds in this case.

\end{example}

\begin{rem}\rm In this section, we have only considered a permutation test of the hypothesis $H^{(1)} \! \! :\rho(1) = 0$. However, analogously, one may prove a similar result for a permutation testing procedure for the hypothesis $H^{(k)} \! \!: \rho_k = \rho(k) = 0$, where $k \in \nn$ is fixed. Indeed, note that the sequence $\{\tilde{Y}_i :i \geq 1\}$, where $\tilde{Y}_i = X_i X_{i + k} $ is $\alpha$-mixing, with $\alpha$-mixing coefficients given by 

\beq
\alpha_\xi(n) = \alpha_X ( n - k) \,\, .
\eeq

\noindent Furthermore, for $\Pi_n$ a random permutation independent of the $X_i$, under appropriate moment conditions for the $X_i$, we also have that 

\beq
\frac{1}{\sqrt{n}} \sum_{i=1} ^{n-1}X_{\Pi_n(i)} X_{\Pi_n(i+1)} \overset{d} = \frac{1}{\sqrt{n}} \sum_{i=1} ^{n-k}X_{\Pi_n(i)} X_{\Pi_n(i+k)} + o_p (1) \, \, ,
\eeq

\noindent since, for any fixed element $\sigma \in S_n$, $\Pi_n \sigma \overset{d}  = \Pi_n$. Hence, defining an appropriate estimator of the variance of $\hat{\rho}_n (k)$, we may similarly construct an asymptotically valid permutation test under the hypothesis $H^{(k)} \! $. To be precise, let $Y_i  = \left(X_i - \bar{X}_n\right)\left(X_{i+k} - \bar{X}_n\right)$, and let $Z_i = \left(X_i - \bar{X}_n\right)^2$. Let

\beq
\begin{aligned}
\hat{T}_{n, \, k} ^2 &= \frac{1}{n} \sum_{i=1} ^{n-1} \left(Y_i - \bar{Y}_n \right)^2 + \frac{2}{n}\sum_{j=1}^{b_n}\sum_{i=1} ^{n-j-k} \left(Y_i - \bar{Y}_n \right) \left(Y_{i+j} -\bar{Y}_n \right) \\
\hat{\nu}_{n, \, k} &=\frac{1}{n} \sum_{i=1}^{n-1} \left(Y_i - \bar{Y}_n\right)\left(Z_i -\bar{Z}_n \right) + \frac{1}{n} \sum_{j=1}^{b_n} \sum_{i=1} ^{n- j-k} \left(Z_i - \bar{Z}_n \right)\left(Y_{i + j} - \bar{Y}_n\right)  +\\ &+\frac{1}{n} \sum_{j=1} ^{b_n} \sum_{i=1} ^{\min\{n-j, \, n-k\}} \left(Y_i - \bar{Y}_n \right)\left(Z_{i+j} - \bar{Z}_n\right) \,\,.
\end{aligned} 
\eeq

\noindent Let $\hat{K}_{n} ^2$ and $\kappa^2$ be defined as in Lemma \ref{lem.stud.unperm.pap}. Let

\beq
\hat{\gamma}_{n, \, k} ^2 = \frac{1}{\hat{\sigma}_n ^4} \left( \hat{T}_{n, \, k} ^2 - 2 \hat{\rho}_{n}(k) \hat{\nu}_k + \hat{\rho}_{n}(k) ^2 \hat{K}_n^2 \right) \, \, ,
\label{gamma.k.est}
\eeq

\noindent where $\hat{\rho}_{n}(k)$ is the sample $k$-th order autocorrelation. Let 

\beq
\begin{aligned}
\tau_k ^2 &= \var\left(X_1 X_{k +1} \right) + 2 \sum_{j \geq 2} \cov \left( X_1 X_{k+1}, \, X_j X_{j+k} \right) \\
\nu_k &:= \cov \left( X_1 X_{k+1}, \, X_1 ^2 \right) + \sum_{j \geq 2} \cov\left( X_1^2, \, X_{j} X_{j+k} \right) + \sum_{j \geq 2} \cov \left(X_1 X_{k+1},\, X_j ^2 \right) \, \, ,
\end{aligned}
\eeq

and let

\beq
{\gamma}_k ^2 = \frac{1}{\sigma ^4} \left( \tau_k ^2 - 2 \rho_k \nu_k + \rho_k ^2 \kappa^2 \right) \, \, .
\eeq 

\noindent Assume that $\tau_k ^2 ,\,  \kappa^2  \in \rr_+$ and $\gamma_k ^2 \in \rr_+$. Under the same conditions on the sequence $\{X_n,\, n \in \nn\}$ as in Theorem \ref{perm.test.big.thm}, by an identical argument to the one given in the case $k = 1$, we will have that the permutation distribution based on the test statistic $\sqrt{n} \hat{\rho}_n (k)/\hat{\gamma}_{n, \, k}$, with associated group of transformations $S_n$, will satisfy (\ref{stud.perm.test.1d}), and the test statistic will satisfy a central limit theorem analogous to (\ref{main.test.stat.1d}). $\qed$

\label{rem.multiple}
\end{rem}

Having developed a permutation testing framework, we further derive an array version of Theorem \ref{perm.test.big.thm}, in order to provide a procedure under which one may compute the limiting power of the permutation test under local alternatives.

\begin{thm} For each $n \in \nn$, let $\left\{X_i^{(n)}, \, i \in [n]\right\}$, be stationary sequences of random variables. Suppose that the $X_i ^{(n)}$ are bounded, uniformly in $i$ and $n$, and that 

\beq
\sum_{n \geq 1} \sup_{r \geq n+1} \alpha_{X^{(r)}} (n) < \infty \, \, .
\eeq

\noindent The permutation distribution of $\sqrt{n} \hat{\rho}_n$, based on the test statistic $\hat{\rho}_n = \hat{\rho}_n \left(X_1 ^{(n)}, \, \dots, \, X_n ^{(n)} \right)$, with associated group of transformations $S_n$, the symmetric group of order $n$, satisfies, as $n \to \infty$,

\beq
\sup_{t \in \rr} \left| \hat{R}_n (t)-\Phi(t)\right| \overset{p}{\to} 0 \, \, .
\eeq

\label{1d.local.perm}

\end{thm}

We may view Theorem \ref{1d.local.perm} as an extension of Theorem \ref{alpha.bd.perm}. Analogously, we may extend Theorem \ref{alpha.perm}, and Lemmas \ref{lem.stud.unperm.pap} and \ref{pap.var.perm.stud}, and hence the result of Theorem \ref{perm.test.big.thm} holds for triangular arrays of stationary, $\alpha$-mixing sequences, replacing the condition (\ref{1d.moment.cond}) with

\beq
\sup_{n \geq 1} \ee \left[ \left|X_1 ^{(n)}\right| ^{8 + 4 \delta}\right] < C \, \, , 
\label{array.moment.cond}
\eeq

\noindent and the condition (\ref{1d.mixing.cond}) with 

\beq
\sum_{n \geq 1} \max_{r \geq n+1} \alpha_{X^{(r)}} (n)^\frac{\delta}{2 + \delta} < \infty \, \, .
\label{array.mixing.cond}
\eeq

\noindent In particular, it follows that we may apply the result of Theorem \ref{1d.local.perm} to triangular arrays of stationary sequences, where instead of (\ref{main.test.stat.1d}), we have the result

\beq
\frac{\sqrt{n} \left(\hat{\rho}_n -{\rho_n }\right)}{\hat{\gamma}_n} \overset{d} \to N(0, \, 1) \, \, ,
\eeq

\noindent where $\rho_n$, and $\hat{\gamma}_n$ are defined analogously to Theorem \ref{perm.test.big.thm}.

It follows that we may compute the power function of the permutation test under appropriate sequences of limiting local alternatives.

\begin{example}\rm\label{ex:ar1} (\emph{AR(1) process}) Consider a triangular array of AR$(1)$ processes, given by

\beq
X_i ^{(n)} = \rho_n X_{i-1}^{(n)}  + \epsilon_i ^{(n)} \, \, , \, \, \, i \in \{2, \, \dots, \, n \} \, \, ,
\label{ar.1.ex}
\eeq

\noindent where $\rho_n = h/\sqrt{n}$, for some fixed constant $h \in  (0, \, 1)$, and the $\epsilon_i ^{(n)}$ form a triangular array of independent standard Gaussian random variables. For each $n$, the autoregressive process defined in \ref{ar.1.ex} has a unique stationary solution, in which $\left( X_1 ^{(n)} , \, \dots, \, X_n ^{(n)} \right)$ follows a multivariate Gaussian distribution, with mean 0 and covariance matrix given by 

\beq
\cov\left(X_i ^{(n)}, \, X_j^{(n)} \right) = \frac{\rho_n^{\left| i - j \right| }}{1 - \rho_n ^2}  \, \, .
\eeq

\noindent Consider the problem of testing $H^{(1)} \! \!: \rho_1 = 0$ against the alternative $\rho_1 > 0$, using the permutation test described in Theorem \ref{1d.local.perm}.

By Theorem 1 of \cite{mokkadem}, we have that condition (\ref{array.mixing.cond}) is satisfied for e.g. $\delta = 1/2$, and, since the $X_i ^{(n)}$ have uniformly bounded second moment and are normally distributed, we also have that condition (\ref{array.moment.cond}) is satisfied. 

Hence, letting $\phi_n$ denote the permutation test conducted on the sequence $\left\{X_i ^{(n)},\,i \in [n] \right\}$,   we may apply the analogous result of Theorem \ref{1d.local.perm} to the triangular array of AR processes, whence, by an application of Slutsky's theorem, we obtain that 

\beq
\frac{\sqrt{n} \hat{\rho}_n}{\hat{\gamma}_n} - \frac{h}{\gamma^{(n)}} \overset{d} \to N(0, \, 1) \, \, ,
\eeq

\noindent and that the local limiting power function satisfies, for $z_{1-\alpha}$ the upper $\alpha$ quantile of the standard Gaussian distribution, 

\beq
\ee_{\rho_n} \phi_n \to 1- \Phi\left( z_{1- \alpha} - \lim_{n\to \infty} \frac{h   }{\gamma^{(n)} }\right) \, \, ,
\eeq

\noindent and

\beq
\left(\gamma^{(n)} \right)^2 = \frac{1}{\sigma_n ^4} \left[ \left(\tau^{(n)} \right)^2 - 2\rho_n \nu_1^ {(n)} + \rho_n ^2 \left( \kappa^{(n)} \right)^2 \right] \, \, ,
\eeq

\noindent where 

\beq
\begin{aligned}
\left( \kappa^{(n)} \right)^2 &=  \var \left(\left(X_1^{(n)}\right) ^2 \right) + 2\sum_{k \geq 2} \cov \left( \left(X_1^{(n)}\right) ^2, \, \left(X_k^{(n)}\right) ^2 \right)\\
\left(\tau^{(n)} \right)^2 &=\var\left(X^{(n)}_1 X_2 ^{(n)} \right) + 2 \sum_{k \geq 2} \cov \left( X_1^{(n)} X_2^{(n)}, \, X_k ^{(n)}X_{k+1} ^{(n)}\right) \\
\nu_1^{(n)} &:= \cov \left( X_1^{(n)} X_2^{(n)}, \, \left(X_1^{(n)}\right) ^2 \right) + \sum_{k \geq 2} \cov\left( \left(X_1^{(n)}\right)^2, \, X_{k} ^{(n)} X^{(n)}_{k+1} \right)+\\&+ \sum_{k \geq 2} \cov \left(X_1^{(n)} X_2^{(n)},\, \left(X_k^{(n)}\right)^2 \right) \, \, .
\end{aligned}
\eeq

\noindent Example 7.16 of \cite{vdv} establishes local asymptotic normality of the local alternative sequence to the null model corresponding to $h = 0$. Hence, by contiguity of the sequence of alternatives, it follows that, as $n \to \infty$,

$$
\gamma^{(n)} \to 1 \, \, ,
$$

\noindent and so, as $n \to \infty$,

\beq
\ee_{\rho_n} \phi_n \to 1- \Phi\left( z_{1- \alpha} - h \right) \, \, .
\label{power.eq}
\eeq
 








\label{local.power.ex}
\end{example}

\begin{rem} \rm We observe that, in the setting of Example \ref{local.power.ex}, the one-sided studentized permutation test is LAUMP (see \cite{tsh}, Definition 13.3.3). Indeed, by Example 7.16 and Theorem 15.4 in \cite{vdv}, coupled with the result of Lemma 13.3.2 in \cite{tsh}, we observe that the optimal local power of a one-sided test, against the alternatives $\rho_n = h/\sqrt{n}$, in the setting of Example \ref{local.power.ex}, is 

$$
\beta^* = 1 - \Phi\left( z_{1- \alpha}  - h  \right) \, \, .
$$

\noindent Since this is exactly the power in (\ref{power.eq}), it follows that the studentized permutation test is LAUMP.

\end{rem}

\begin{rem} \rm Note that, more generally, the same argument applies when computing the limiting local power of the studentized permutation test with respect to contiguous alternatives. Indeed, if contiguity can be established for some sequence of alternatives $\{P_n, \, n \in \nn\}$ with first-order autocorrelations $\rho_{1, \, n}= h/\sqrt{n}$, by a similar argument to the one presented in Example \ref{local.power.ex}, we will have that the convergence of $\hat{\gamma}_n^2$ to $\gamma_1^2$ (in probability) also holds under the contiguous sequence of alternatives. Hence the limiting power of the one-sided level $\alpha$ studentized permutation test, under the contiguous sequence of alternatives, will also be given by

$$
\ee_{P_n} \phi_n \to 1- \Phi\left(z_{1-\alpha} - \frac{h}{\gamma_1} \right) \, \, , 
$$

\noindent as $n \to \infty$. 
\end{rem}

\section{Multiple and joint hypothesis testing}\label{sec.multiple}

In this section, we outline multiple testing procedures which may be applied to test the hypotheses $H^{(k)}$, as defined in (\ref{eq.piece.hyp}), simultaneously. While we make use of the standard Bonferroni method of combining $p$-values, we argue such an approach is not overly conservative. We develop a method for testing joint null hypotheses of the form

$$
H_r: \rho(1) = \rho(2) = \dots = \rho(r) = 0 \, \, .
$$

\noindent It is desirable to perform such a test in a multiple testing framework, i.e. in the case of rejection of the null hypothesis, we often wish to accompany this rejection with inference on which of the individual hypotheses $H^{(k)}$ do not hold. To this end, it is necessary to construct a procedure allowing for valid inference, in the sense that the familywise error rate (FWER) is controlled at the nominal level $\alpha$. In general, we may apply the canonical Bonferroni correction; that is, given marginal $p$-values $\hat{p}_1, \, \dots, \, \hat{p}_r$ and a nominal level $\alpha$, we reject the null hypothesis $H_r$ if 

$$
\min_i \hat{p}_i \leq \frac{\alpha}{r} \, \, ,
$$

\noindent and assert that the hypothesis $H^{(k)}$ does not hold for any $k$ such that $\hat{p}_k \leq \alpha/r$. 
Under the further assumption of independence of $p$-values, we may use multiple testing procedures such as the \v{S}id\'ak correction, which rejects any $H^{(k)}$ for which 

$$
\hat{p}_k \leq 1 - (1- \alpha)^\frac{1}{r} \, \, .
$$

\noindent This procedure is marginally more powerful that the canonical Bonferroni procedure, but may not control FWER at the nominal level $\alpha$ if there is negative dependence between the $\hat{p}_k$. In order to understand the dependence structure between sample autocorrelations, and their corresponding permutation $p$-values, we have the following result.

\begin{thm} In the setting of Theorem \ref{perm.test.big.thm}, let $r \in \nn, \, r > 1$. For $k \in [r]$, let $\rho_k$ be the $k$th-order autocorrelation, and let $\hat{\rho}_k$ be the $k$th-order sample autocorrelation. Let $\Sigma \in \rr^{(r+1) \times (r+1)} = (\sigma_{ij})_{i, \, j = 0} ^r$ be such that 

$$
\sigma_{ij} = \begin{cases} 
&\var(X_1 X_{1 + i}) + 2 \sum_{l > 1} \cov(X_1 X_{1 + i}, \, X_l X_{l + i} )\, \, , \, i = j \\
&\cov(X_1 X_{1+i}, \, X_1 X_{1+j}) + \sum_{l > 1} \left[ \cov(X_1 X_{1 + i}, \, X_l X_{l+j} ) + \cov(X_1 X_{1 + j}, \, X_l X_{l+i} )\right] \, \, , \, i \neq j \, \, .
\end{cases}
\label{sigma.multiple}
$$

\noindent Let $A \in \rr^{(r+1) \times r }$ be given by 

$$
A = \begin{pmatrix} 
-\frac{\rho_1}{\sigma^4} & \dots & \dots &\dots&- \frac{ \rho_r}{\sigma^4} \\
\frac{1}{\sigma^2} & 0 & \dots&\dots  & 0 \\
0 & \frac{1}{\sigma^2} & 0 &\dots & 0 \\
\vdots &\vdots &\vdots & \vdots & \vdots \\
0 & \dots &\dots &\dots & \frac{1}{\sigma^2} 
\end{pmatrix} \, \, .
$$

\noindent Then, as $n \to \infty$, 

$$
\sqrt{n}  \left [ \begin{pmatrix}
\hat{\rho}_1\\
\vdots \\
\hat{\rho}_r \end{pmatrix}  - \begin{pmatrix}
{\rho}_1\\
\vdots \\
{\rho}_r \end{pmatrix} \right] \overset{d} \to N\left(0, \, A^T \Sigma A\right) \, \, .
$$

\label{thm.multiple}
\end{thm}

\begin{rem} \rm We observe that, in the i.i.d. setting, the sample autocorrelations are asymptotically independent. Indeed, in this case, we have that $\Sigma$, as defined in (\ref{sigma.multiple}), is diagonal, and, for $i \neq j$, for all $l \in \{0, \, \dots, \, r\}$, 

$$
A_{li} A_{lj} = 0 \, \, .
$$

\noindent Therefore, for $i, \, j \in [r]$, $i \neq j$, 

$$
\begin{aligned}
\left(A^T \Sigma A\right)_{ij} &= \sum_{l, \, s = 0} ^r A_{li} A_{sj} \sigma_{ls} \\
&= \sum_{l = 0} ^r A_{li}A_{lj} \sigma_{ll} \\
&= 0 \, \, . \, \, \, \, \,\qed
\end{aligned}
$$

\label{rem.multi.indep}
\end{rem}

By Remark \ref{rem.multi.indep} and the uniform convergence of the permutation distribution $\hat{R}_n$ in Theorem \ref{perm.test.big.thm}, we have that, for $1 \leq k \leq r$, leaving the dependence of $\hat{p}_k$ on $n$ implicit,

$$
\hat{p}_k = 1 - \Phi\left(\frac{\hat{\rho}_{n}(k)}{\hat{\gamma}_{n, \, k}} \right) + o_p (1) \, \, ,
$$

\noindent where $\hat{\gamma}_{n, \, k}$ is as defined in (\ref{gamma.k.est}). It follows that, in some settings, such as the i.i.d. setting, the marginal $p$-values are asymptotically independent. Therefore we may use the \v{S}id\'ak correction if, for instance, we use the null hypothesis $H_r$ as a portmanteau test of independence of realizations. However, more generally, using the Bonferroni cutoff of $\alpha /r$ is only marginally larger than the Bonferroni-\v{S}id\'ak correction, and it applies irrespective of the asymptotic dependence structure of the marginal $p$-values.
Indeed, since, by Theorem \ref{thm.multiple} and Remark \ref{rem.multi.indep}, any method must at least account for possibility of asymptotic independence, it follows that the cutoff should be at least as large as the \v{S}id\'ak correction.  But since the Bonferroni correction is not much larger than the \v{S}id\'ak correction, there does not appear to be much gain, in terms of power, in devising a method that precisely accounts for the joint dependence among the marginal $p$-values. Despite this, we may use a step-down procedure, such as that of \cite{holm}, to obtain a larger power.


We illustrate the application of the canonical Bonferroni procedure in Section \ref{sec.fin}, in application to historical log-return data.

\section{Simulation results}\label{sec.simulations} 

Monte Carlo simulations illustrating our results are given in this section. Tables \ref{tab:mdep} and \ref{tab:alpha} tabulate the rejection probabilities of one-sided tests for the permutation tests, in addition to those of the Ljung-Box and Box-Pierce tests. The nominal level considered is $\alpha =0.05$. The simulation results confirm that the permutation test is valid, in that, in large samples, it approximately attains level $\alpha$. The simulation results also confirm that, by contrast, both the Ljung-Box and Box-Pierce tests perform extremely poorly in non-i.i.d. settings. 

As a review, the Ljung-Box and Box-Pierce tests are used to test for independence of residuals in fitting ARMA models. This is done by a portmanteau test of the null hypothesis $H_r$, as defined in (\ref{eq.joint.hyp}), which is tested under the assumption that the residuals follow a Gaussian white noise process. For each $k \in \nn$, let $\hat{\rho}_k$ be the sample $k$th order autocorrelation. The one-sided Ljung-Box test compares the test statistic 

\beq
\hat{Q}_{LB, \, n} = n(n+2) \sum_{k=1} ^r \frac{\hat{\rho}_k ^2 }{n-k} 
\label{eq.ljung.box}
\eeq

\noindent to the quantiles of a $\chi_r ^2$ distribution, with rejection occurring for large values of $\hat{Q}_{LB,\, n}$. Similarly, the one-sided Box-Pierce test compares the test statistic 

\beq
\hat{Q}_{BP, \, n} = n \sum_{k=1} ^r {\hat{\rho}_k ^2 }
\label{eq.ljung.box}
\eeq

\noindent to the quantiles of a $\chi_r ^2$ distribution, with rejection occurring for large values of $\hat{Q}_{BP, \, n}$. In the case of $r = 1$, both the one-sided Ljung-Box and Box-Pierce tests compare the test statistic

\beq
\hat{Q}_n = C_n \hat{\rho}_1 ^2 
\label{eq.box}
\eeq

\noindent to the quantiles of a $\chi_1 ^2$ distribution, with rejection in both tests occurring for large values of $\hat{Q}_n$. The Ljung-Box and Box-Pierce tests primarily differ in their scaling in this case; namely, the Ljung-Box test takes $C_n = n(n+2)/(n-1)$ in (\ref{eq.box}), while the Box-Pierce test uses $C_n = n$.

In this simulation, we consider both $m$-dependent (in Table \ref{tab:mdep}) and $\alpha$-mixing (in Table \ref{tab:alpha}) processes. Table \ref{tab:mdep} gives the null rejection probabilities for sampling distributions of the form described in Example \ref{ex:mdep:gaussian}, in the case of Gaussian products, where the values of $m$ are listed in the first column. Note that $m = 0$ corresponds to the setting where the $X_i$ are independent standard Gaussian random variables. 

Table \ref{tab:alpha} gives the null rejection probabilities for processes of the form described in Example \ref{ex:ar2}, with $\rho = 0.5$. We include one additional example in the second row of Table \ref{tab:alpha}. The sample distribution in this row is as follows. $\{X_{2i}, \, i \in [n/2]\}$ and $\{X_{2i - 1},\, i \in [n/2]\}$ are independent and identically distributed sequences, with 

\beq
X_{2i} = Y_{2i} Y_{2(i+1)} \, \, ,
\eeq

\noindent for $Y$ as in Example \ref{ex:ar2}, with $\phi = 0$ and $\rho = 0.5$, and standard Gaussian innovations. 

For each situation, 10,000 simulations were performed. Within each simulation, the permutation test was calculated by randomly sampling 2,000 permutations. 

\begin{table}[]
\centering 
\begin{tabular}{cc rrrrrrr} 
\hline\hline 
$m$&$n$&10  &20&    50  &  80 &  100   &500 & 1000 \\ [0.5ex]
\hline 
\multirow{4}{*}{0}	&Stud. Perm.&0.0511 &0.0489&0.0465&0.0452&0.0500&0.0488&0.0525\\
				&Unst. Perm. &0.0503&0.0511&0.0493&0.0470&0.0494&0.0480&0.0527\\
   				 &Ljung-Box&0.0534&0.0544&0.0474&0.0488&0.0516&0.0488&0.0482 \\
    				 &Box-Pierce&0.0198&0.0365&0.0407&0.0448&0.0484&0.0480&0.0478\\
				 \hline
\multirow{4}{*}{1}&Stud. Perm.& 0.0654&0.0578&0.0611&0.0586&0.0582&0.0532&0.0534\\
				&Unst. Perm. &0.1010&0.1249&0.1388&0.1512&0.1553&0.1692&0.1749 \\
   				 &Ljung-Box&0.0888&0.1390&0.1873&0.2084&0.2102&0.2359&0.2651\\
    				 &Box-Pierce&	0.0342&0.1057&0.1737&0.2001&0.2039&0.2341&0.2645\\
				 \hline
\multirow{4}{*}{2}&Stud. Perm.&0.0718&0.0638&0.0661&0.0615&0.0683&0.0608&0.0582\\
				&Unst. Perm. &0.1288&0.1588&0.2041&0.2189&0.2327&0.2580&0.2721 \\
   				 &Ljung-Box&0.0999&0.1912&0.2975&0.3420&0.3494&0.4425&0.4645\\
    				 &Box-Pierce&0.0455&0.1555&0.2841&0.3319&0.3410&0.4414&0.4638\\
				 \hline
\multirow{4}{*}{3}	&Stud. Perm.&0.0714&0.0708&0.0638&0.0713&0.0706&0.0647&0.0566\\
				&Unst. Perm. &0.1411&0.1716&0.2332&0.2582&0.2748&0.3269&0.3364 \\
   				 &Ljung-Box&0.1000&0.2056&0.3404&0.4026&0.4310&0.5644&0.6034\\
    				 &Box-Pierce&0.0451&0.1693&0.3252&0.3946&0.4233&0.5634&0.6033
\\
\hline \hline
\end{tabular}
\caption{Monte Carlo simulation results for null rejection probabilities for tests of $\rho(1) = 0$, in an $m$-dependent Gaussian product setting.} \label{tab:mdep}
\end{table}

\begin{table}[] 
\hspace{-1.4cm}
\begin{tabular}{cc rrrrrrr} 
\hline\hline 
Distribution&$n$&10  &20&    50  &  80 &  100   &500 & 1000 \\ [0.5ex]
\hline 
\multirow{4}{*}{AR(2), $N(0, \, 1)$ innov.}&Stud. Perm.&0.0418&0.0215&0.0399&0.0420&0.0448&0.0464&0.0480\\
				&Unst. Perm. &0.0594&0.0901&0.1212&0.1403&0.1372&0.1541&0.1658
\\
   				 &Ljung-Box&0.2101&0.2360&0.2570&0.2492&0.2537&0.2527&0.2607 \\
    				 &Box-Pierce&0.1283&0.1972&0.2411&0.2399&0.2468&0.2516&0.2603\\
				 \hline
\multirow{4}{*}{AR(2) Prod., $N(0, \, 1)$ innov.}&Stud. Perm.&0.0385&0.0370&0.0332&0.0375&0.0382&0.0350&0.0366\\				
&Unst. Perm. &0.0555&0.0648&0.0836&0.0923&0.0939&0.1136&0.1181 \\
   				 &Ljung-Box&0.0938&0.0878&0.1012&0.1151&0.1189&0.1674&0.1708 \\
    				 &Box-Pierce&0.0493&0.0658&0.0913&0.1081&0.1138&0.1653&0.1705\\
				 \hline
\multirow{4}{*}{AR(2), $U[-1, \, 1]$ innov.}&Stud. Perm.&0.0496&0.0243&0.0390&0.0433&0.0444&0.0470&0.0464 \\
				&Unst. Perm. &0.0628&0.0940&0.1276&0.1412&0.1395&0.1569&0.1572 \\
   				 &Ljung-Box&0.2272&0.2464&0.2545&0.2539&0.2522&0.2530&0.2522\\
    				 &Box-Pierce&0.1403&0.2051&0.2394&0.2434&0.2445&0.2517&0.2517\\
				 \hline
\multirow{4}{*}{AR(2), $t_{9.5}$ innov.}&Stud. Perm.&0.0432&0.0218&0.0385&0.0423&0.0437&0.0531&0.0459\\
				&Unst. Perm. &0.0582&0.0902&0.1182&0.1346&0.1362&0.1581&0.1634\\
   				 &Ljung-Box&0.2050&0.2316&0.2567&0.2522&0.2524&0.2654&0.2590\\
    				 &Box-Pierce&0.1206&0.1941&0.2416&0.2436&0.2455&0.2638&0.2579\\
\hline \hline
\end{tabular}
\caption{Monte Carlo simulation results for null rejection probabilities for tests of $\rho(1) = 0$, in an $\alpha$-mixing setting.} 
\label{tab:alpha}
\end{table}

\begin{figure}[]
\begin{center}
\begin{tabular}{cc}
{\mbox{\epsfxsize=85mm\epsfbox{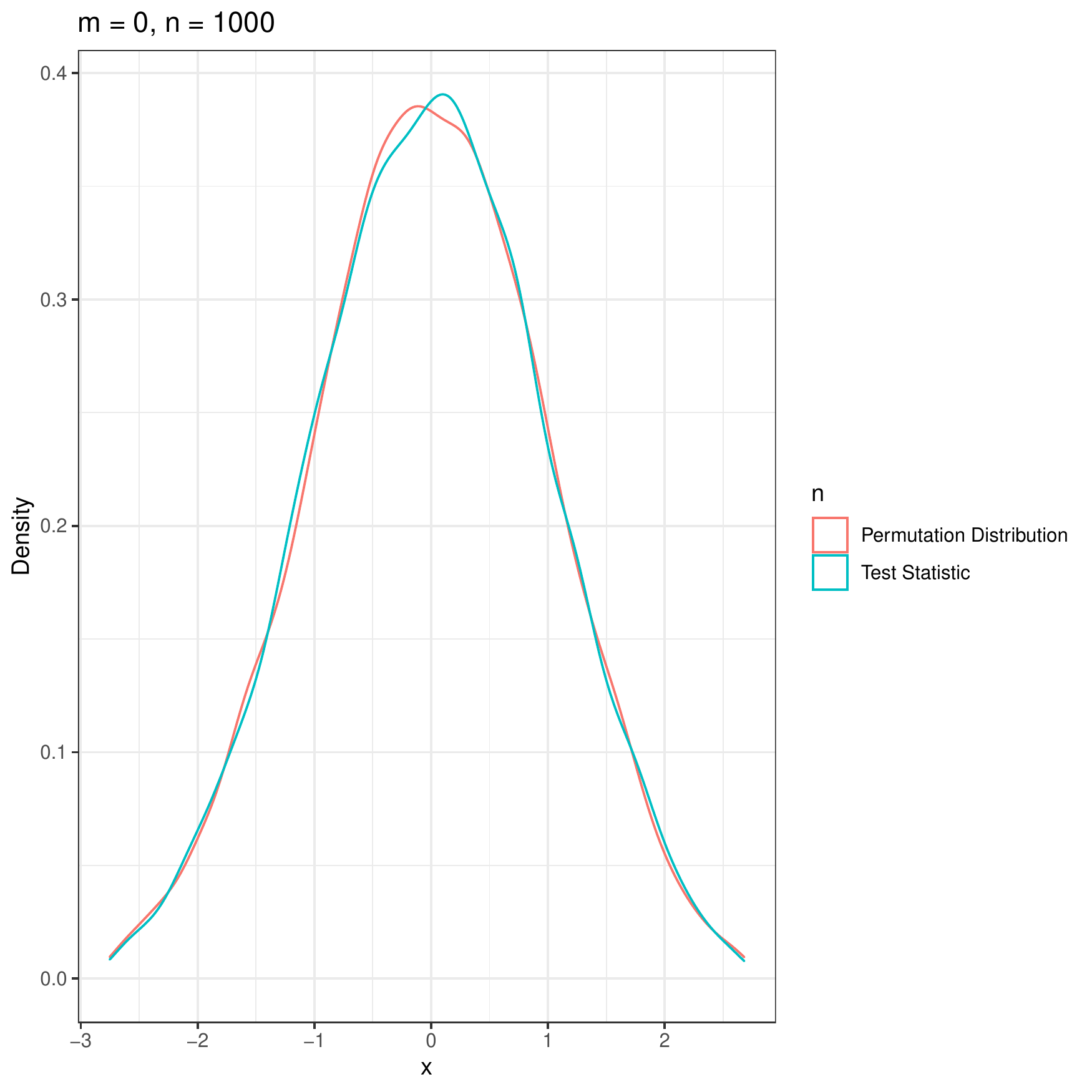}}} 
& {\mbox{\epsfxsize=85mm\epsfbox{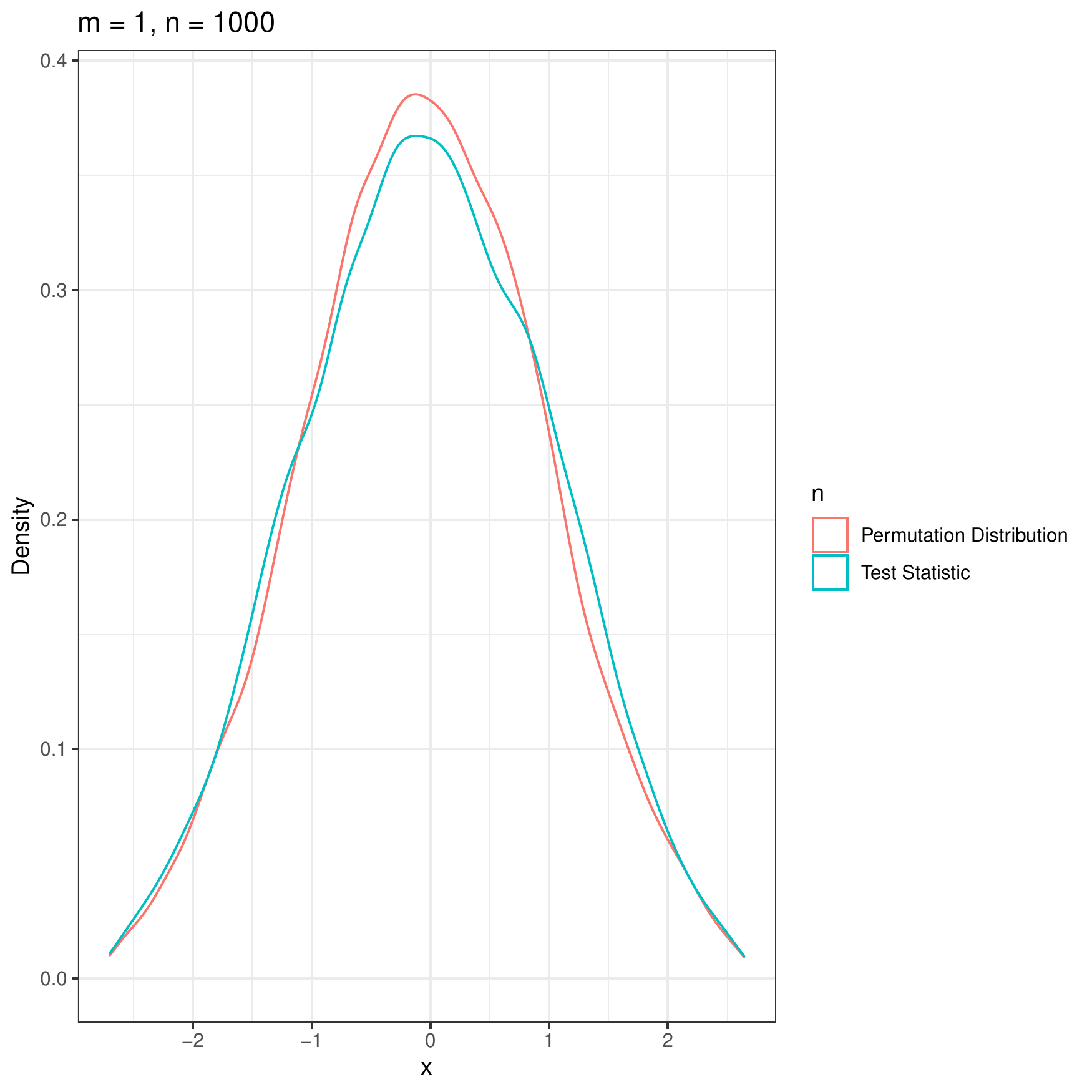}}} \\
{\mbox{\epsfxsize=85mm\epsfbox{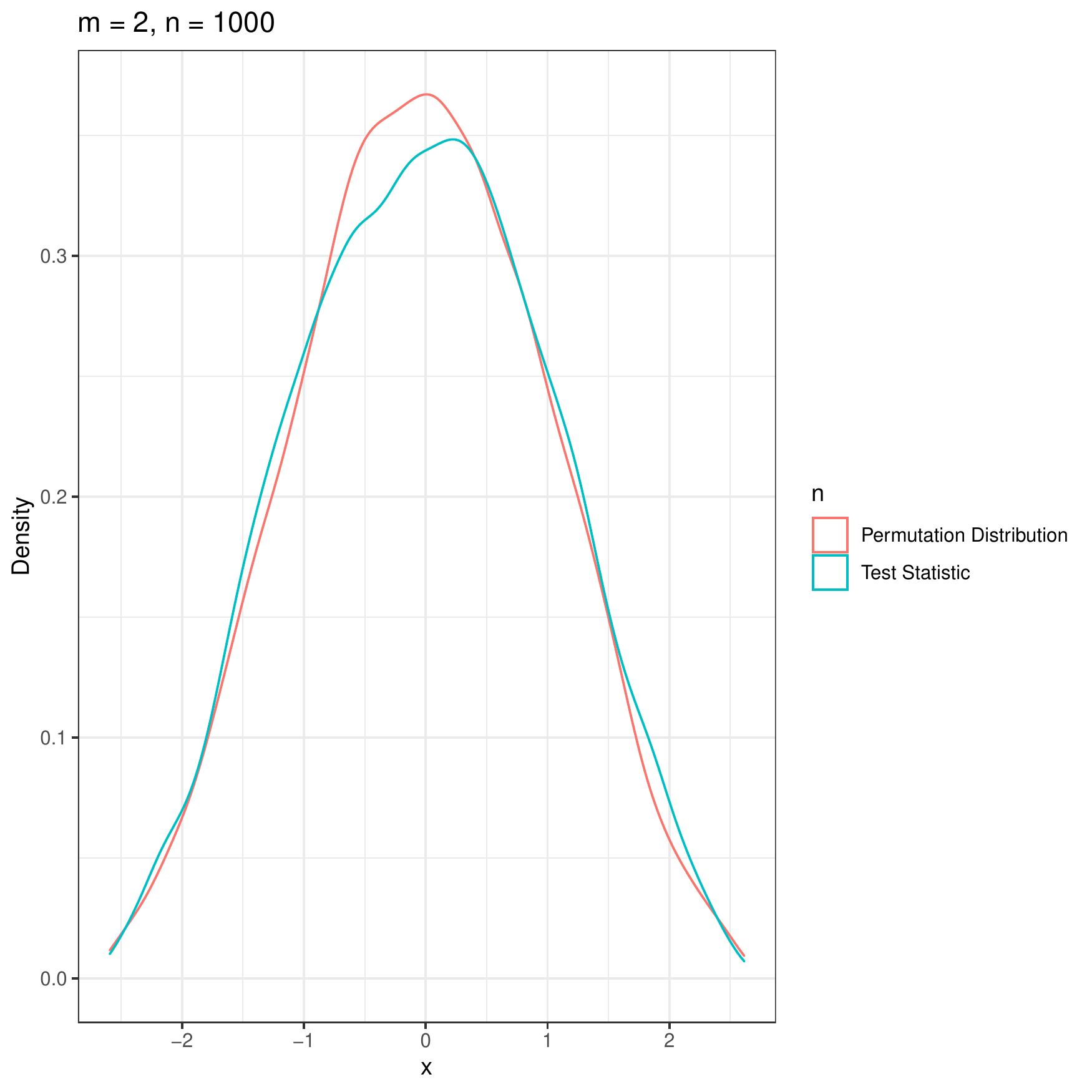}}} &
{\mbox{\epsfxsize=85mm\epsfbox{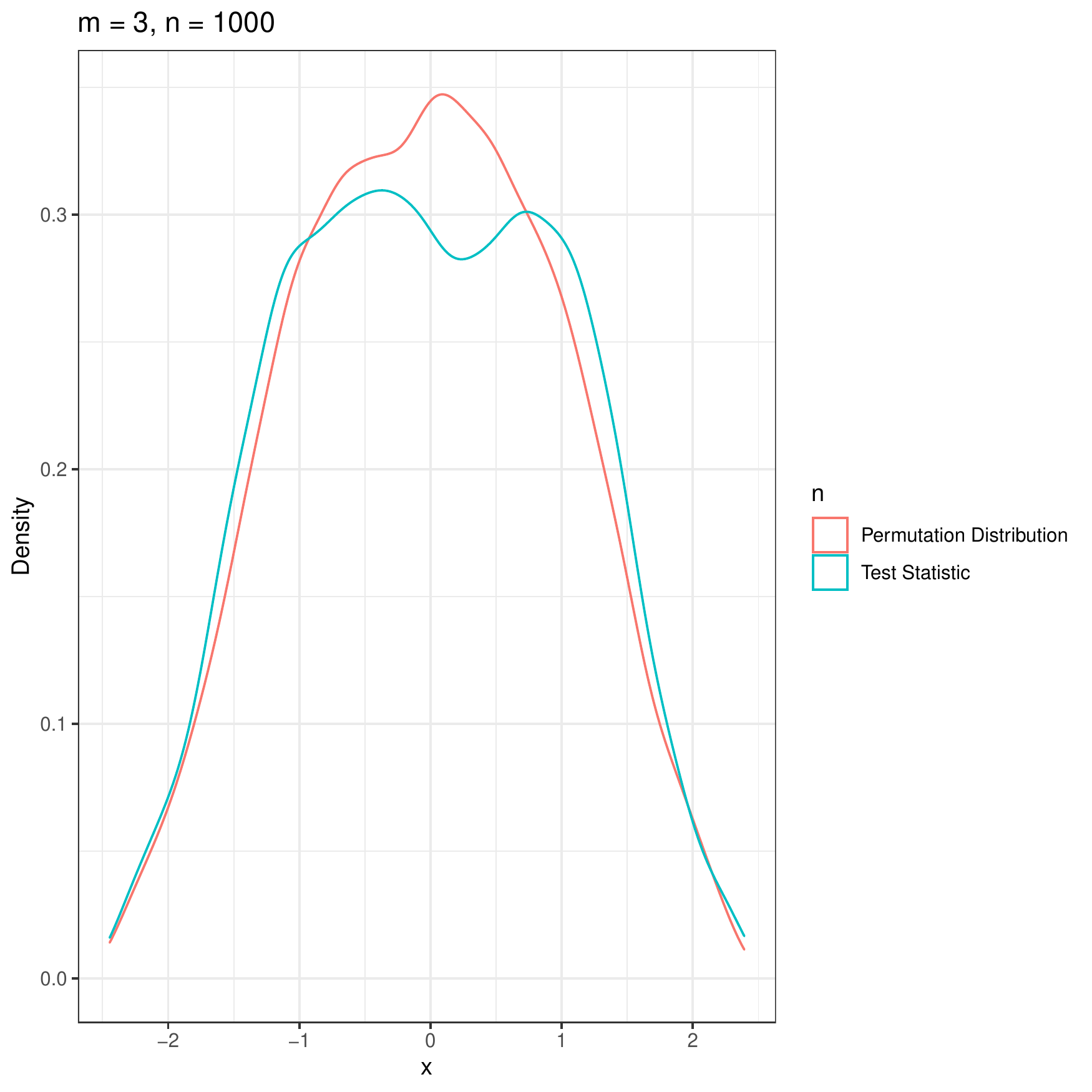}}}
\end{tabular}
\end{center}
\caption[a]{Figure showing kernel density estimates of the densities of the test statistic and permutation distribution in the $m$-dependent case. In the case of the permutation distribution, the KDE of the permutation distribution, pooled across simulations, is provided. The kernel used for the KDE is a Gaussian kernel.
}
\label{fig:density.mdep}
\end{figure}

\begin{figure}[]
\begin{center}
\begin{tabular}{cc}
{\mbox{\epsfxsize=51mm\epsfbox{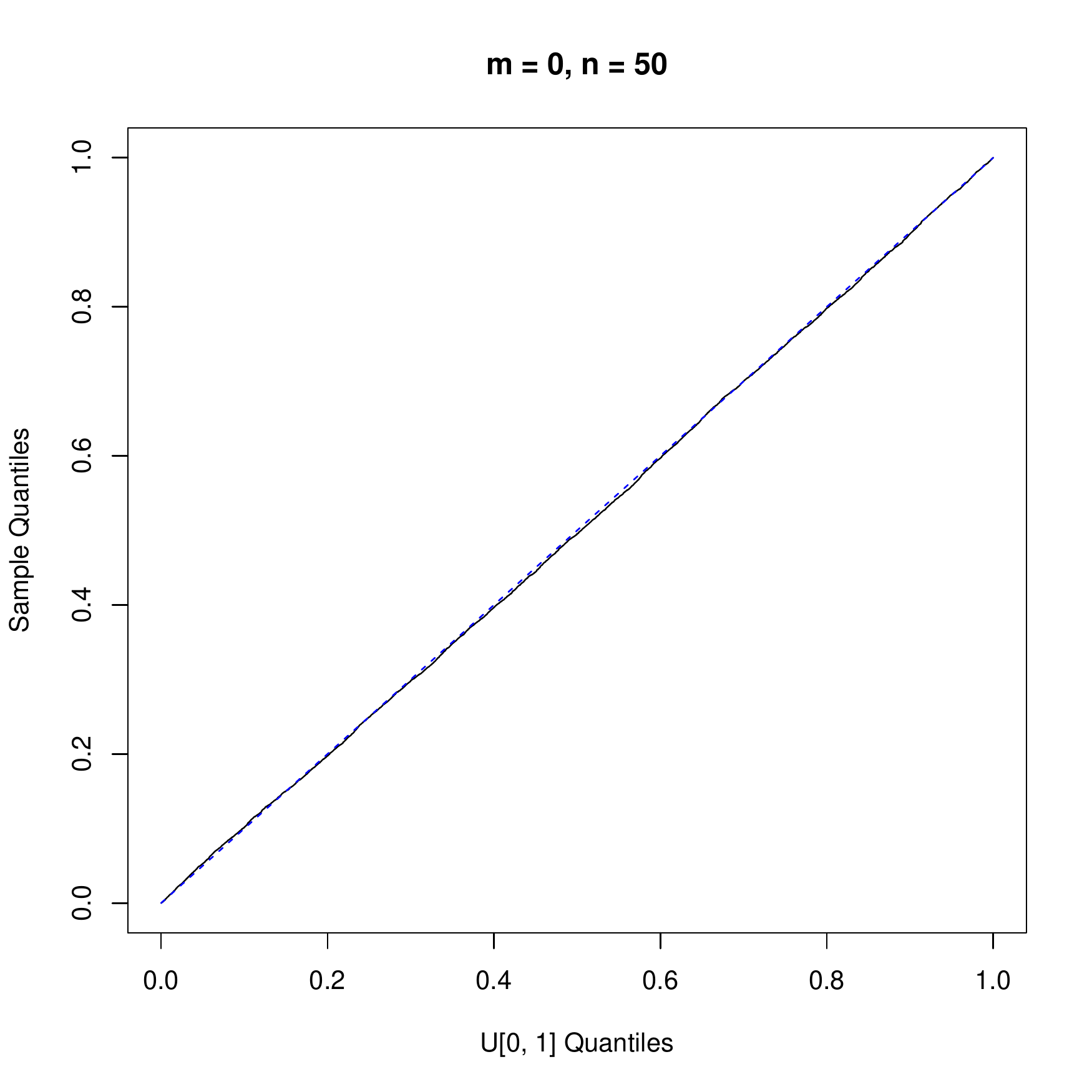}}} 
& {\mbox{\epsfxsize=51mm\epsfbox{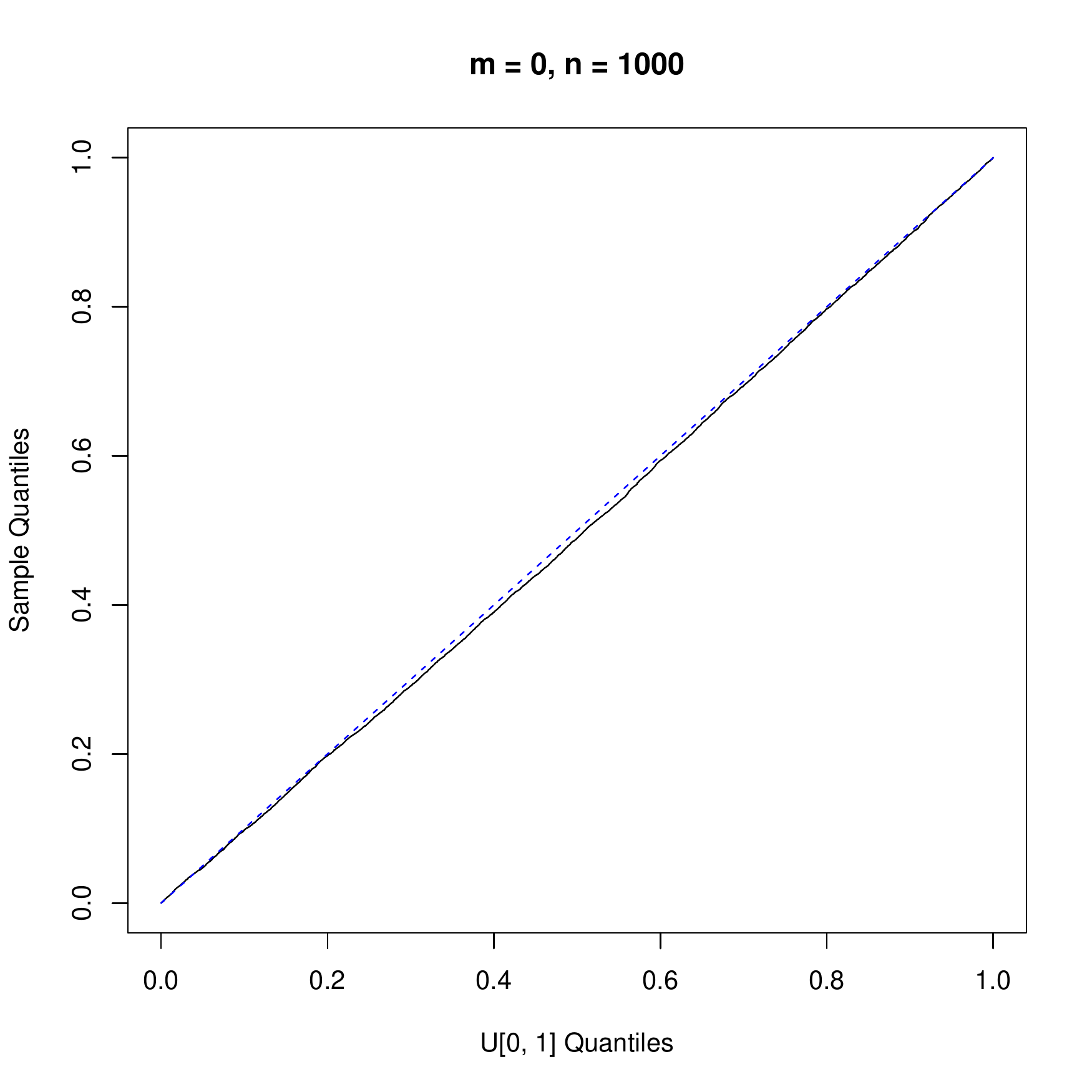}}} \\
{\mbox{\epsfxsize=51mm\epsfbox{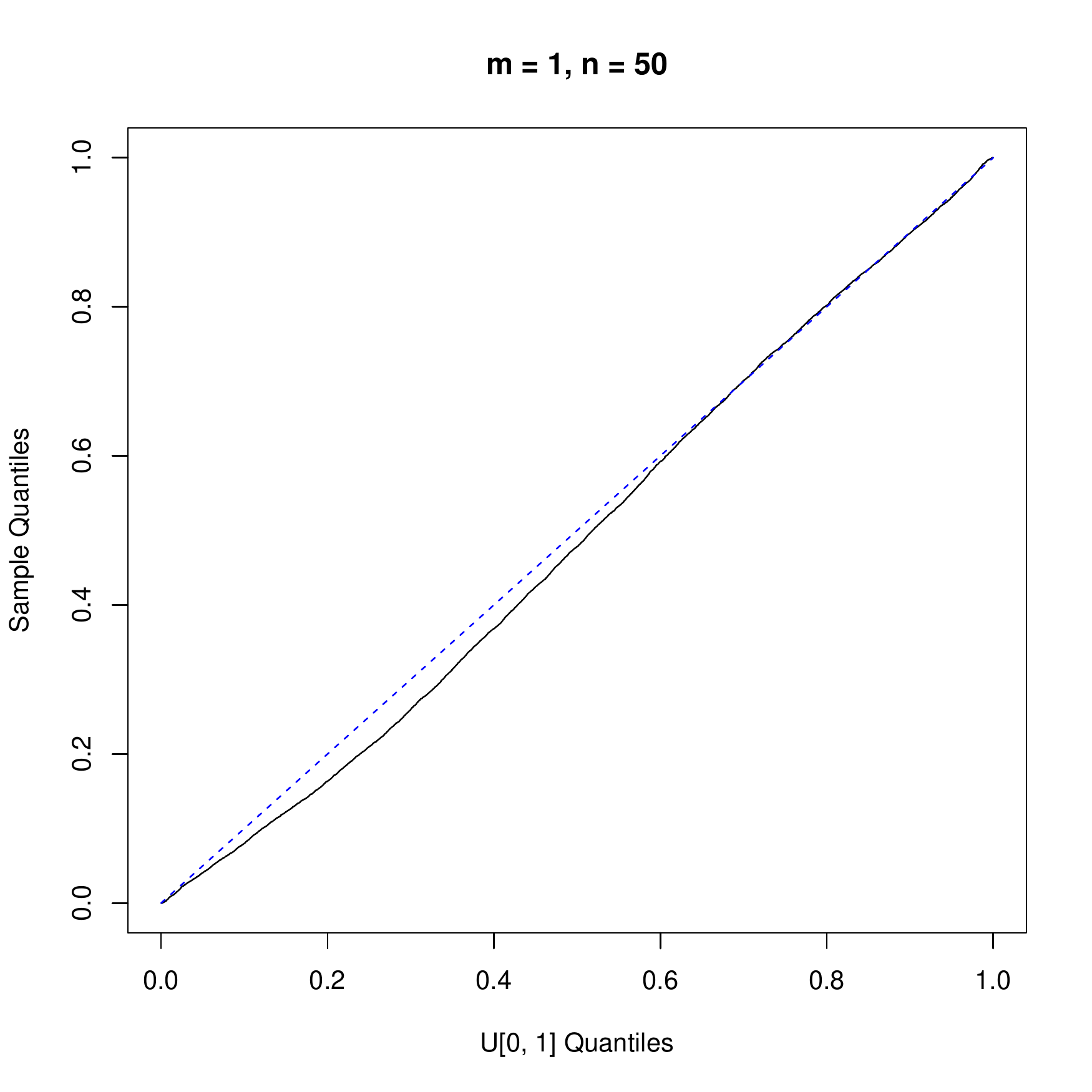}}} &
{\mbox{\epsfxsize=51mm\epsfbox{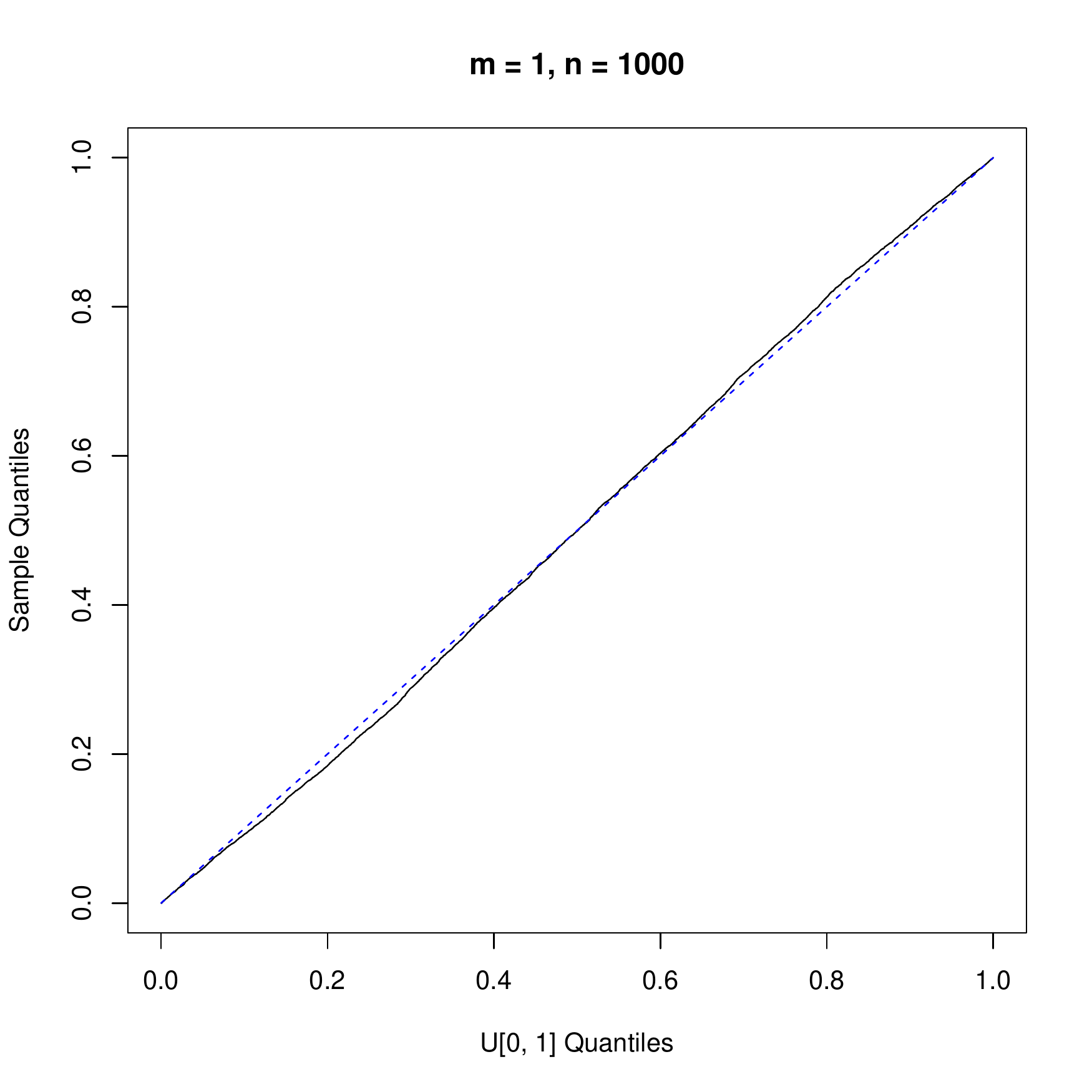}}} \\
{\mbox{\epsfxsize=51mm\epsfbox{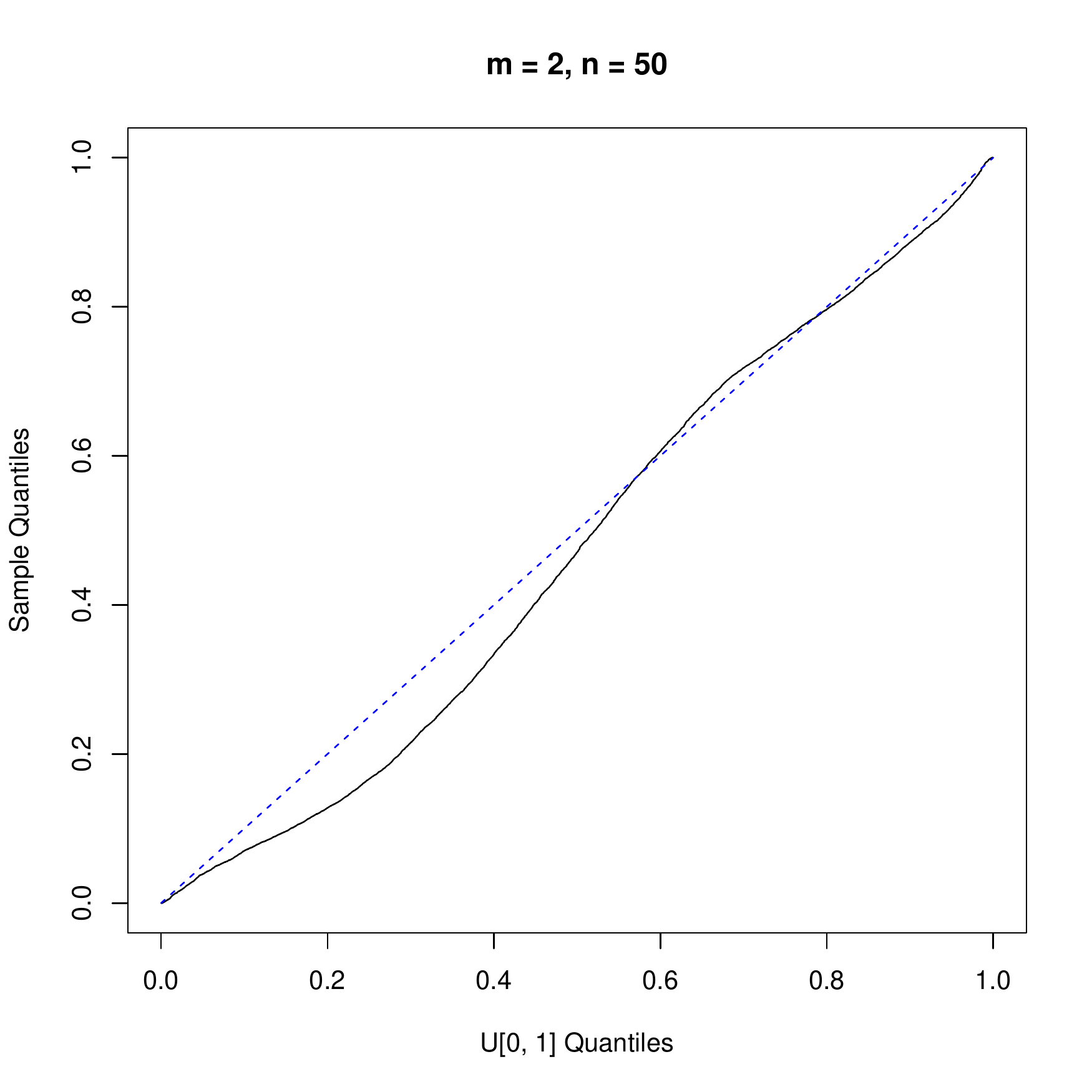}}} 
& {\mbox{\epsfxsize=51mm\epsfbox{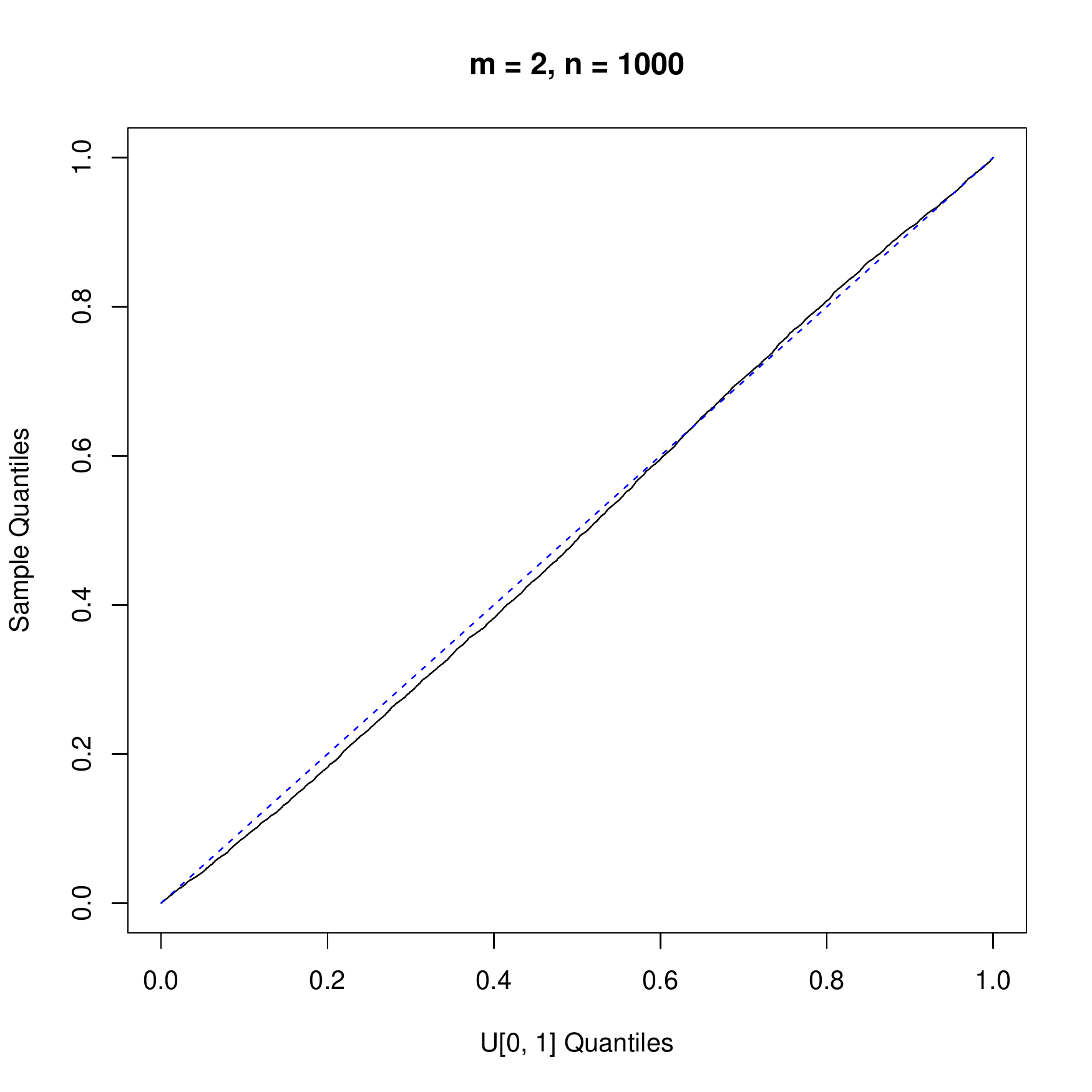}}} \\
{\mbox{\epsfxsize=51mm\epsfbox{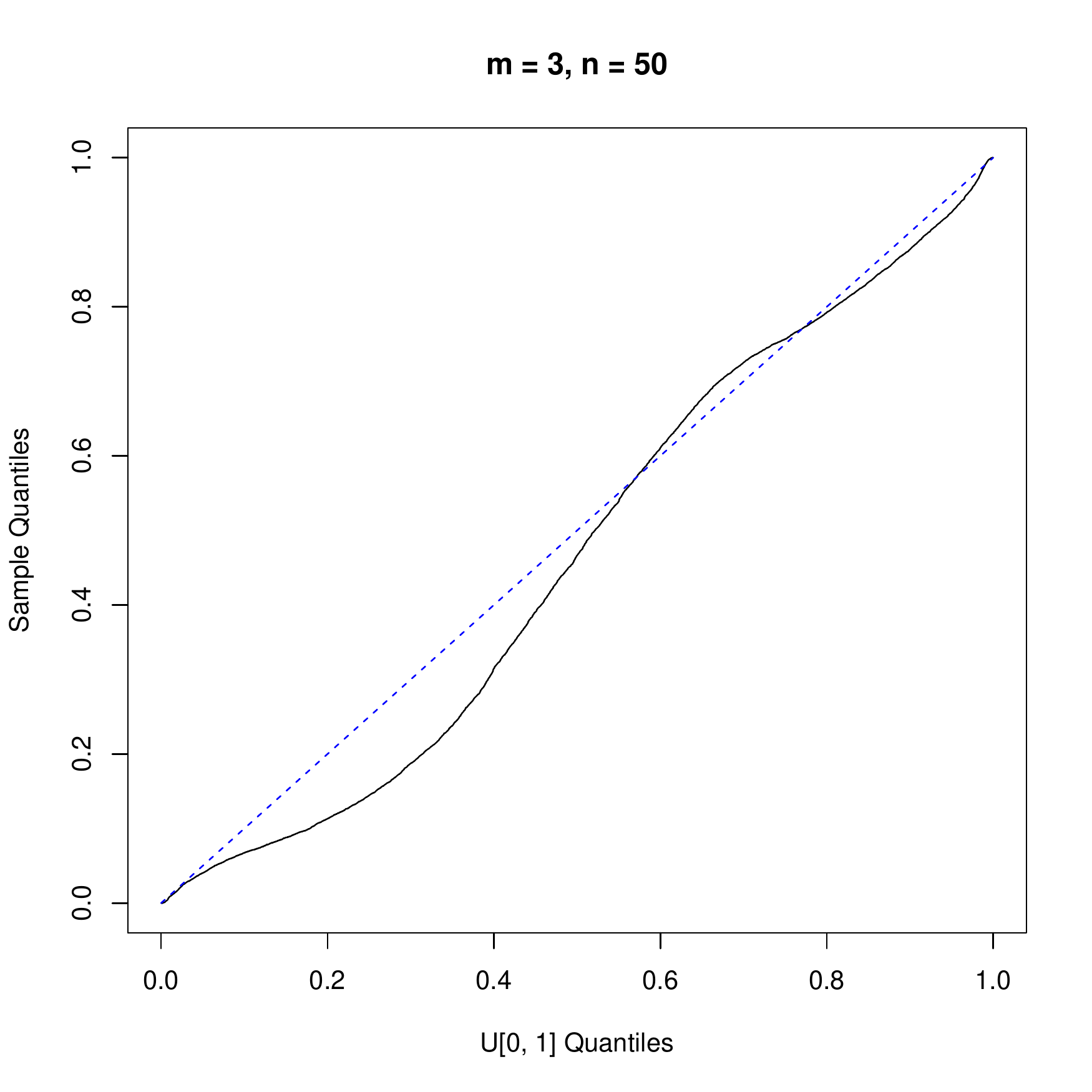}}} &
{\mbox{\epsfxsize=51mm\epsfbox{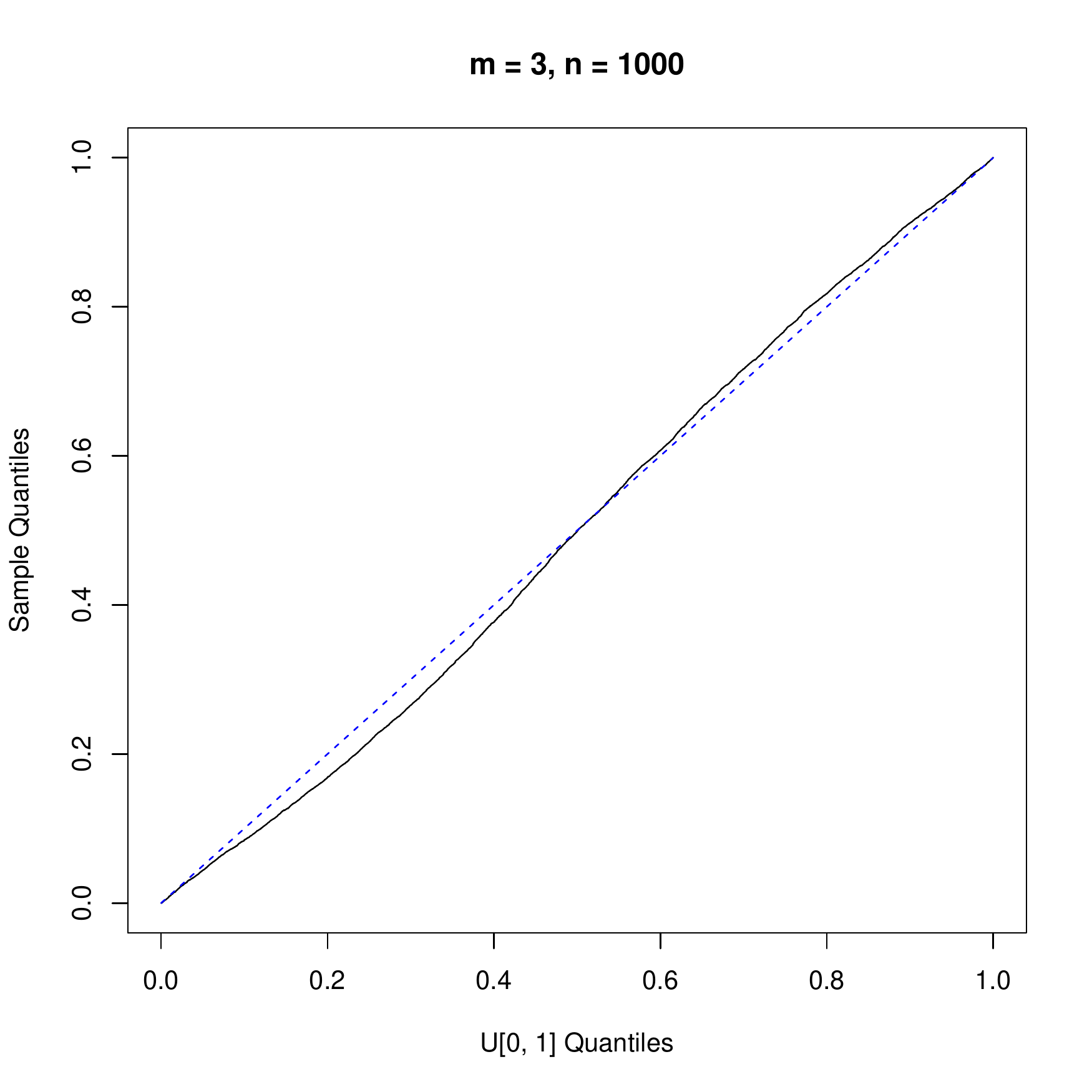}}} \\
\end{tabular}
\end{center}
\caption[a]{Figure showing QQ plots of the sample $p$-values obtained from one-sided permutation tests, in the $m$-dependent Gaussian product setting.
}
\label{fig:qq.mdep}
\end{figure}

The results of the simulation are further illustrated in Figures \ref{fig:density.mdep} and \ref{fig:qq.mdep}. Figure \ref{fig:density.mdep} shows kernel density estimates\footnote{These were obtained using the \texttt{density} function in R, using the default Gaussian kernel.}  of the distributions of the test statistic and the permutation distribution in the $m$-dependent setting described above. Figure \ref{fig:qq.mdep} provides QQ plots of the simulated $p$-values against the theoretical quantiles of a $U[0, \, 1]$ distribution, also in the $m$-dependent setting. These figures further confirm the asymptotic validity of the permutation test procedure in the $m$-dependent setting.

We observe several computational choices to be made when applying the permutation testing framework in practice. By the results of Lemmas \ref{lem.stud.unperm.pap} and \ref{pap.var.perm.stud}, for large values of $n$, the estimate $\hat{\gamma}_n^2 $ will be be strictly positive with high probability. However, for smaller values of $n$, it may be the case that a numerically negative value of $\hat{\gamma}_n ^2$ is observed, either when computing the test statistic or the permutation distribution. A trivial solution to this issue is the truncate the estimate at some sufficiently small fixed lower bound $\epsilon >0$. Note that, for appropriately small choices of $\epsilon$, i.e. $\epsilon < \gamma_1^2$, the results of Lemmas \ref{lem.stud.unperm.pap} and \ref{pap.var.perm.stud} still hold, i.e. inference based on this choice of studentization is still asymptotically valid. In practice, however, the suitability of a choice of $\epsilon$ for a particular numerical application is affected by the distribution of the $X_i$. For the above simulation, a constant value of $\epsilon = 10^{-6}$ was used.

A further choice is that of the truncation sequence $\{b_n,\, n \in \nn\}$ used in the definition of $\hat{\gamma}_n ^2$. Any sequence $  \{b_n\}$ such that, as $n \to \infty$, $b_n \to \infty$ and $b_n = o\left(\sqrt{n} \right)$ will be appropriate, although, in a specific setting, some choices of $\{b_n\}$ will lead to more numerical stability than others. In the simulations above, $\{b_n\}$ was taken to be $b_n = [n^{1/3}] + 1$, where $[x]$ denotes the integer part of $x$.

\begin{figure}[]
\begin{center}
\begin{tabular}{c}
{\mbox{\epsfxsize=100mm\epsfbox{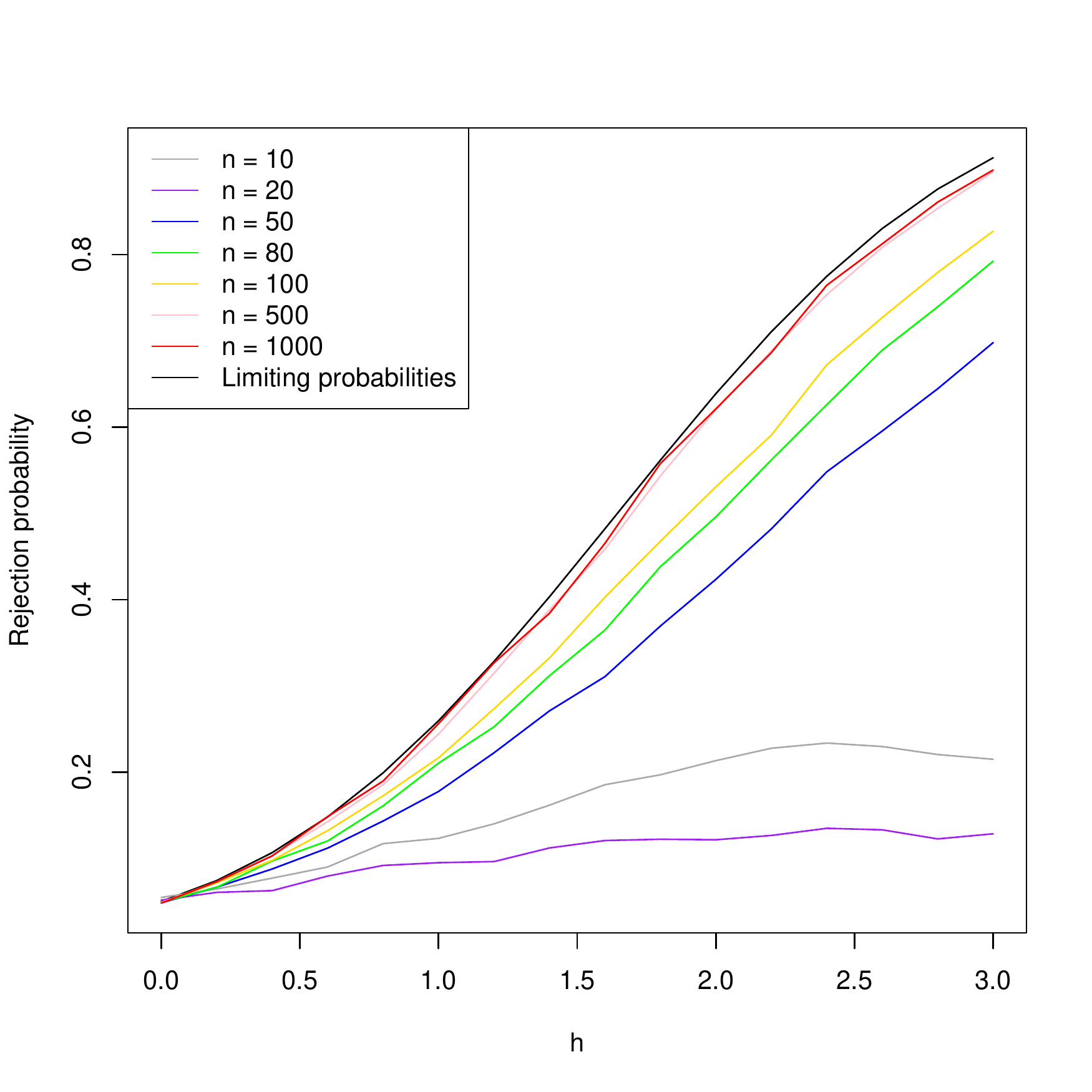}}} 
\end{tabular}
\end{center}
\caption{Monte Carlo simulation results for rejection probabilities of the tests $\rho(1) = 0$, in the setting of Example \ref{local.power.ex}, for different values of $h$ and $n$.} \label{tab.local}
\label{fig.local}
\end{figure}

We also provide Monte Carlo simulation results for the local limiting power of the one-sided studentized permutation test with local alternatives of the form described in Example \ref{local.power.ex}. The nominal level considered is $\alpha = 0.05$. For each situation, 10,000 simulations were performed. Within each simulation, the permutation test was calculated by randomly sampling 2,000 permutations. Figure \ref{fig.local} shows the null rejection probabilities for different values of $h$.


We observe that, for large values of $n$, the sample rejection probabilities are very close to the theoretical rejection probabilities computed in Example \ref{local.power.ex}.

\section{Application to financial data}\label{sec.fin}

In this section, we describe an application of the permutation test to financial stock data. 

Under the assumption that a certain version of the Efficient Market Hypothesis holds true (see \cite{fama_emh} and \cite{malkiel} for details), we have that the daily log-returns of a stock $R_t$, i.e. for $S_t$ the stock price at time $t$,

$$
R_t = \log(S_t) - \log(S_{t-1})  \, \, ,
$$

\noindent are serially uncorrelated. Stronger versions of the Efficient Market Hypothesis assert that the daily log-returns either form a martingale sequence, or are independent. It follows that a test of the lack of correlation of observed daily log-returns can provide evidence for, or against, these versions of the Efficient Market Hypothesis.

\begin{figure}[]
\begin{center}
\begin{tabular}{c}
{\mbox{\epsfxsize=100mm\epsfbox{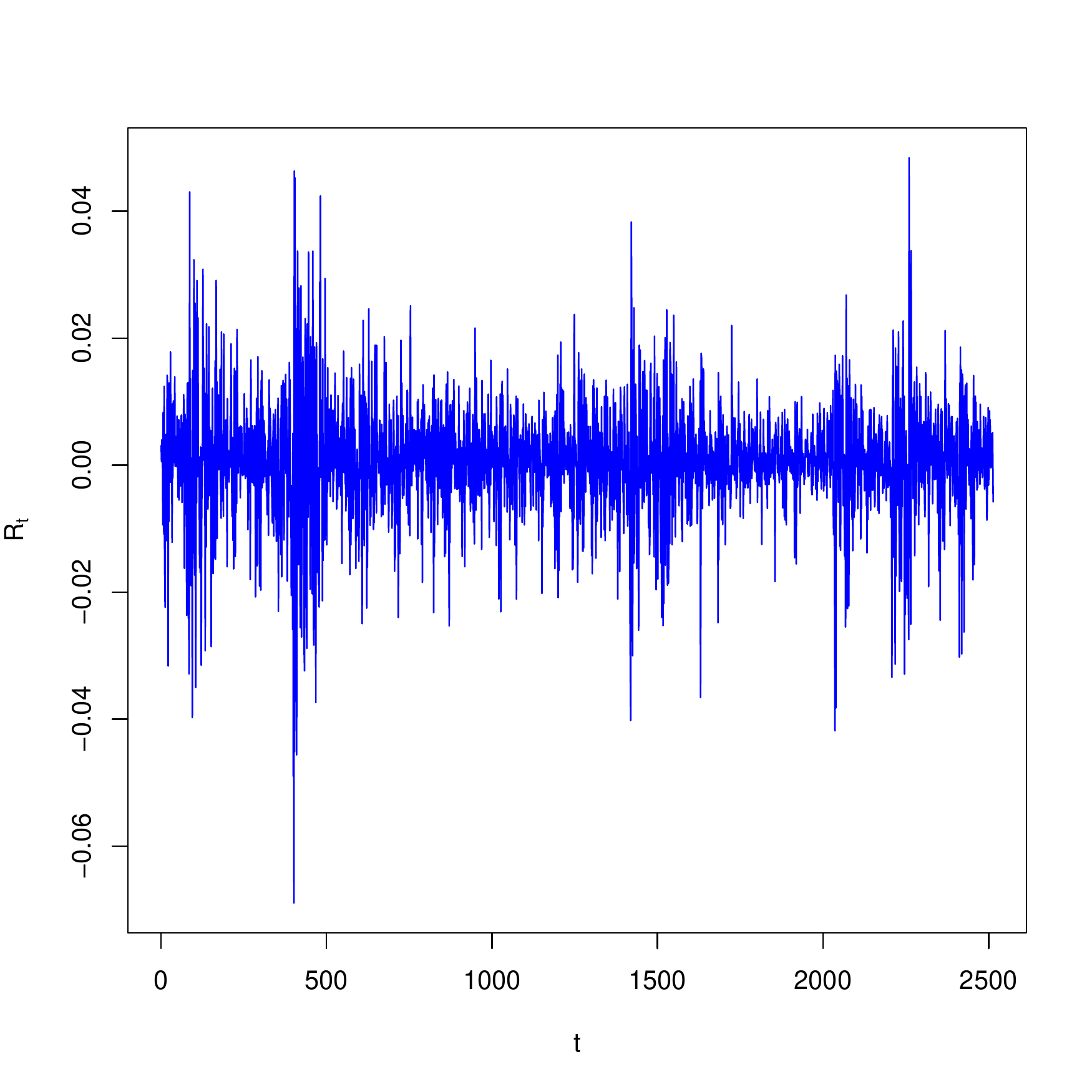}}} \\
{\mbox{\epsfxsize=100mm\epsfbox{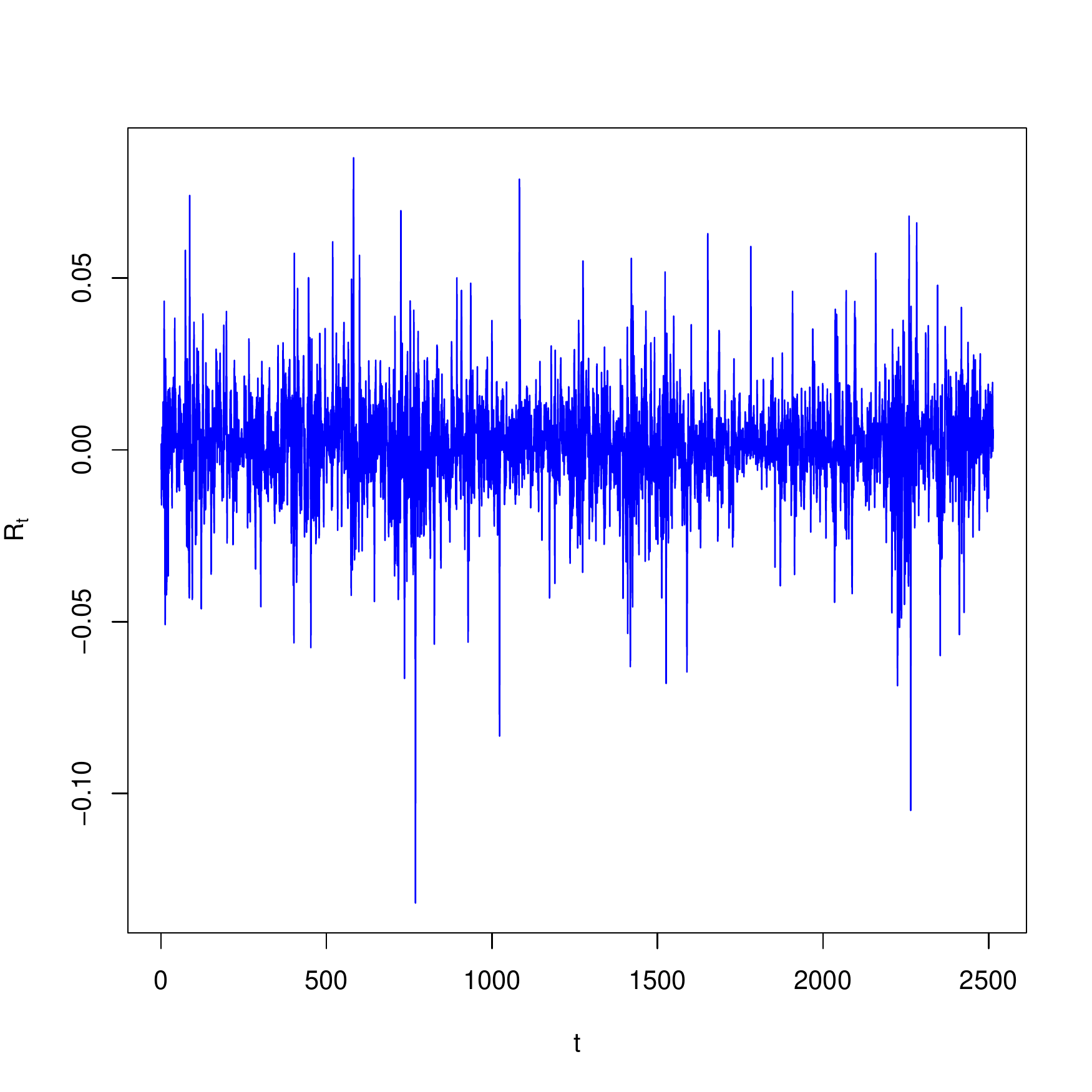}}}
\end{tabular}
\end{center}
\caption[a]{Figure showing plots of the daily log-returns $R_t$ against time. The top plot shows log-returns for the S\&P 500 index, and the bottom plot shows log-returns for Apple stock. In both plots, $t= 0$ corresponds to the first day of trading after 01/01/2010, i.e. 01/04/2010. 
}
\label{fig.ret}
\end{figure}

\begin{figure}[]
\begin{center}
\begin{tabular}{c}
{\mbox{\epsfxsize=95mm\epsfbox{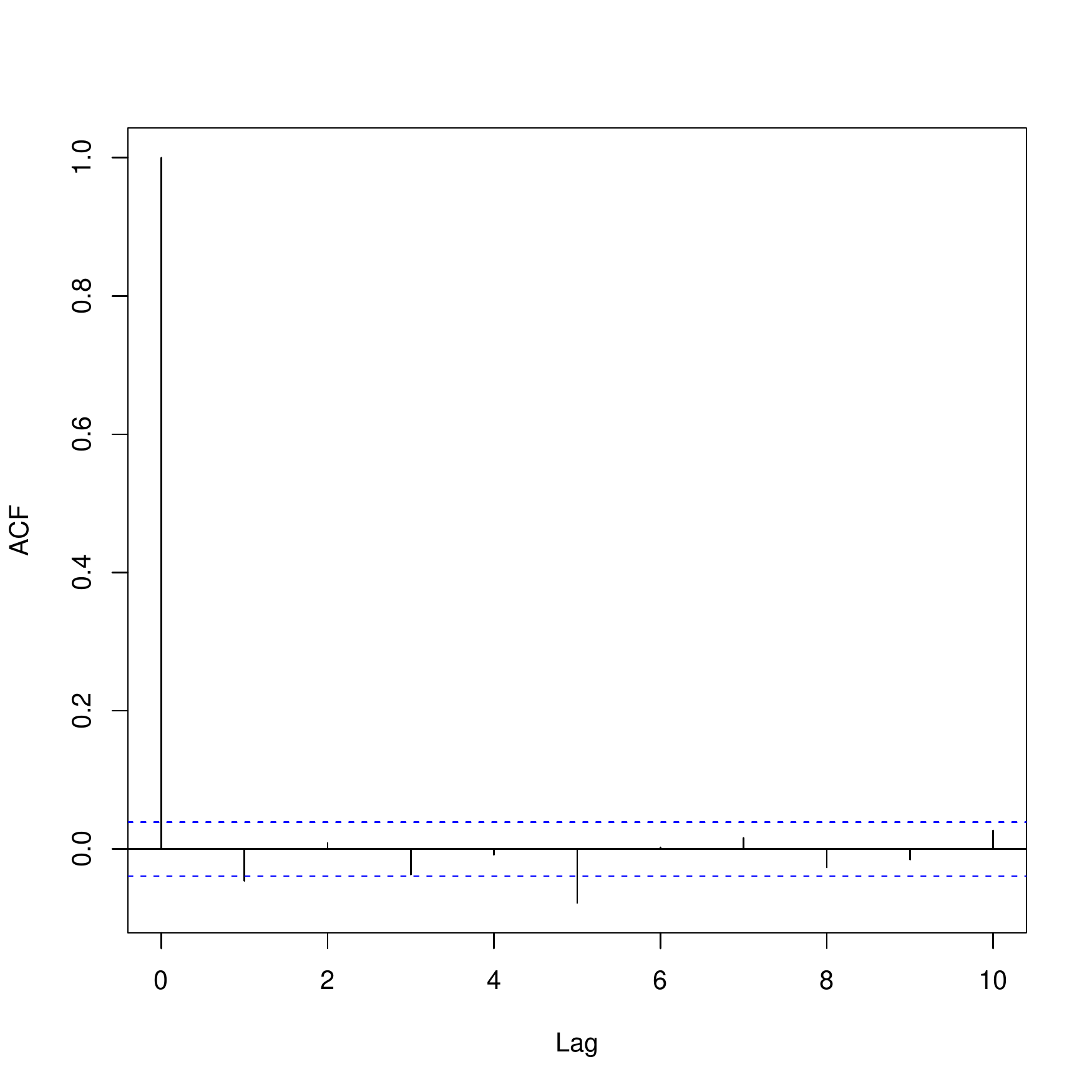}}} \\
{\mbox{\epsfxsize=95mm\epsfbox{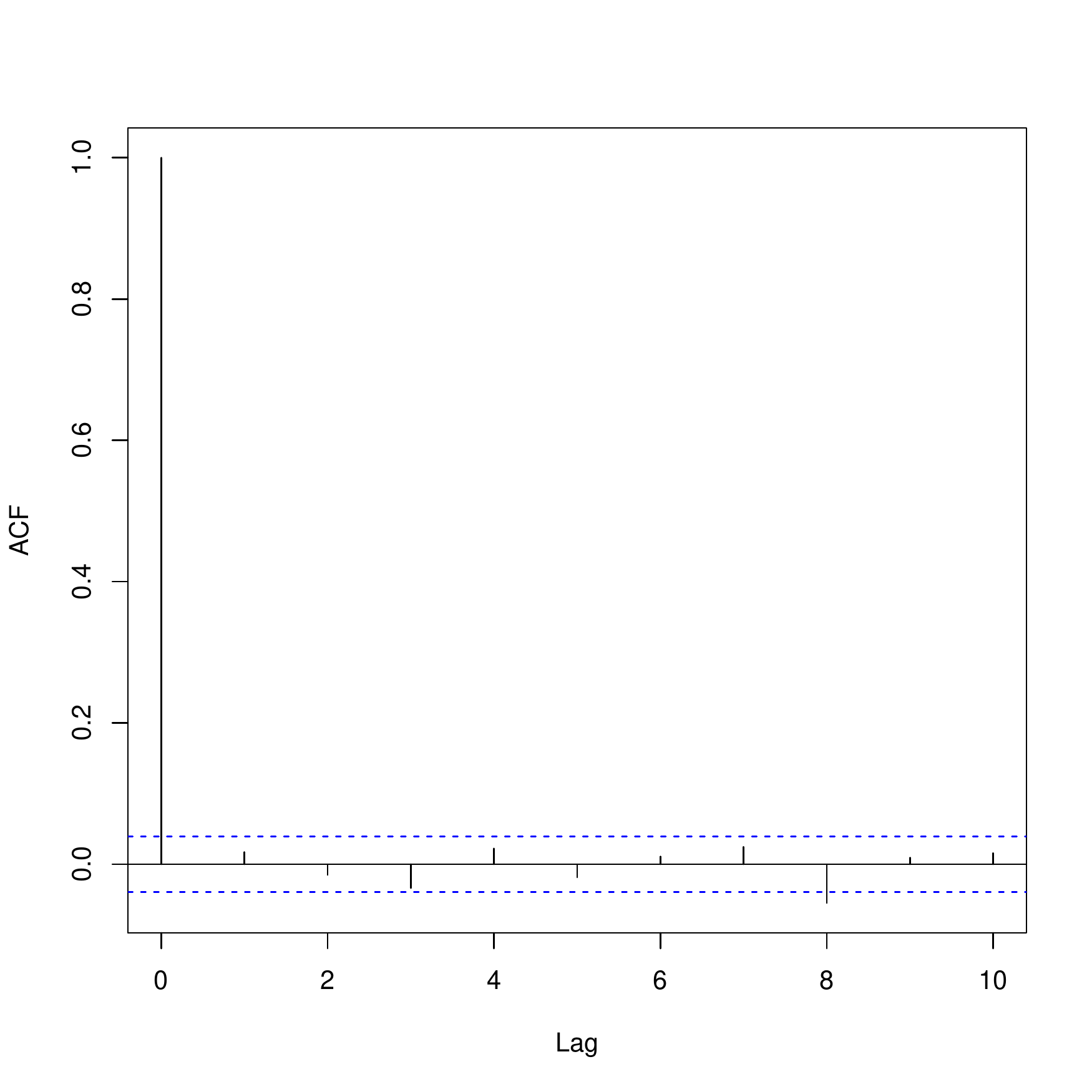}}}
\end{tabular}
\end{center}
\caption[a]{Figure showing plots of the sample autocorrelations $\hat{\rho}_k$ against time. The top plot shows sample autocorrelations for the S\&P 500 index, and the bottom plot shows sample autocorrelations for Apple stock. The dotted blue lines show 95\% confidence intervals under the assumption that the sequence is a Gaussian white noise process.
}
\label{fig.acf}
\end{figure}

We illustrate such portmanteau tests performed on daily closing prices for the S\&P 500 index (SPX) and Apple (AAPL) stock. In particular, the hypothesis 

$$
H_r : \rho(1) = \dots = \rho(r) = 0 \, \, ,
$$

\noindent for $r = 10$, was tested using the studentized permutation tests described in Section \ref{sec.alpha}, with a Bonferroni correction. The results of these tests were compared to the corresponding results of the Ljung-Box test when applied to the data. The test was performed using closing price data, obtained from Yahoo! Finance, between the dates of 01/01/2010 and 12/31/2019. In both cases, days for which data were unavailable, such as weekends and holidays, were omitted. Plots of the log-returns are shown in Figure \ref{fig.ret}, and plots of the sample autocorrelations are shown in figure \ref{fig.acf}.

In both cases, the permutation distribution was approximated using 2,000 random permutations. In addition, the summation parameter used in the studentization term $\hat{\gamma}_n$ used was $b_n = [n^{1/3}]+1$. The $p$-values obtained are shown in Table \ref{tab.ret}.

\begin{table}[]
\centering 
\begin{tabular}{cc rrrrrrrrrr} 
\hline\hline 
$k$&1  &2&    3  &  4 &  5   &6 & 7 &8&9&10\\ [0.5ex]
\hline 
				SPX&0.1675&0.8665&0.4035&0.8305&0.0340&0.9340&0.6025&0.4100&0.4975&0.3610\\
				\hline
				AAPL&0.3375&0.5145&0.1255&0.3435&0.4715&0.6310&0.3845&0.0120&0.6470&0.4895\\
\hline \hline
\end{tabular}
\caption{Table showing the marginal $p$-values obtained using the permutation test for the S\&P 500 index and Apple stock data.} \label{tab.ret}
\end{table}

We observe that, marginally, in the case of the S\&P 500 index, the $p$-value for $k = 5$ was significant, and, in the case of Apple stock, the marginal $p$-value for $k = 8$ was significant. However, in both cases, none of the $p$-values is significant at the $\alpha = 5\%$ level when adjusted using a Bonferroni correction, and so we may conclude that there is no significant evidence in the data for the daily log-returns to indicate deviation from the Efficient Market Hypothesis.

By contrast, the $p$-values obtained using the Ljung-Box test, for all 10 lags simultaneously, were 0.0010 (in the case of the S\&P 500 index), and 0.0827, in the case of Apple stock. While the results in the case of Apple stock are consistent with those of the permutation test, in the case of the S\&P 500 index, we observe that the Ljung-Box test rejects the null hypothesis at the $\alpha = 5\%$ level. However, in light of the results of the permutation test and the simulation results in Section \ref{sec.simulations}, this should not cast doubt upon our conclusion that no significant deviation from the Efficient Market Hypothesis is observed.

\section{Testing autocovariance}
\label{sec.cov}

In this paper, we have discussed testing the null hypothesis 

\beq
H^{(1)} \! \! : \rho_1 = 0 \, \, ,
\eeq 

\noindent where $\rho_1$ is the first-order autocorrelation, using permutation tests with test statistic based on the sample first-order autocorrelation $\hat{\rho}_n$. However, the hypothesis $H^{(1)}$ is equivalent to the null hypothesis

\beq
\tilde{H}^{(1)} \! \!: c_1 = 0 \,\,, 
\eeq

\noindent where $c_1$ is the first-order autocovariance. Since the sample variance $\hat{\sigma}_n ^2$ is permutation invariant, it is clear that an analogous result to that of Theorem \ref{alpha.perm} holds for the permutation distribution based on the test statistic $\sqrt{n}\hat{c}_n$, where $\hat{c}_n$ is the sample first-order autocovariance. In order to obtain a result similar to that of Theorem \ref{thm.unperm.rho.lim}, we can apply Ibragimov's central limit theorem. Then, by results analogous to those of Lemmas \ref{lem.stud.unperm.pap} and \ref{pap.var.perm.stud}, the following holds.

\begin{thm} Let $\{X_n,\, n \in \nn\}$ be a strictly stationary, $\alpha$-mixing sequence, with variance $\sigma^2$ and  first-order autocorrelation $\rho_1$, such that, for some $\delta > 0$,

\beq
\ee \left[ \left| X_1 \right| ^{8 + 4 \delta}\right] < \infty \, \, ,
\label{1d.moment.cond.cov}
\eeq

\noindent and 

\beq
\sum_{n \geq 1} \alpha_X (n) ^\frac{\delta}{2 + \delta} < \infty \, \, .
\label{1d.mixing.cond.cov}
\eeq

\noindent Let $\tau_1^2$ be as in Theorem \ref{thm.unperm.rho.lim}. Suppose that $\tau_1^2 \in (0, \, \infty)$. Let $\hat{T}_n ^2$ be as defined in (\ref{mdep.var.est}). 

\noindent \begin{enumerate}

\item[i)] We have that, as $n \to \infty$,

\beq
\frac{\sqrt{n}\left( \hat{c}_n - c_1 \right)}{\hat{T}_n} \overset{d} \to N \left( 0, \, 1 \right) \, \, .
\label{main.cov.stat.1d}
\eeq

\item[ii)] Let $\hat{R}_n$ be the permutation distribution, with associated group of transformations $S_n$, the symmetric group of order $n$, based on the test statistic $\sqrt{n} \hat{c}_n/\hat{T}_n$. As $n \to \infty$,

\beq
\sup_{t \in \rr} \left| \hat{R}_n (t)-\Phi(t)\right| \overset{p}{\to} 0 \, \, .
\label{stud.cov.test.1d}
\eeq

\end{enumerate}

\label{perm.test.cov.thm}
\end{thm}

In practice, the permutation test for autocovariance produces numerically similar rejection probabillities to the permutation test based on the sample autocorrelation. Table \ref{tab:cov} provides Monte Carlo simulation results for rejection probabilities for the permutation test based on the sample autocovariance, in the same settings as those discussed in Section \ref{sec.simulations}.

\begin{table}[h]
\hspace{-1.4cm}
\begin{tabular}{cc rrrrrrr} 
\hline\hline 
&$n$&10  &20&    50  &  80 &  100   &500 & 1000 \\ [0.5ex]
\hline 
\multirow{4}{*}{$m$-dependent}	&$m=0$&0.0515&0.0483&0.0463&0.0450&0.0498&0.0488&0.0525\\
				&$m=1$&0.0761&0.0571&0.0587&0.0543&0.0551&0.0511&0.0511\\
   				 &$m=2$&0.0830&0.0677&0.0640&0.0546&0.0604&0.0520&0.0511\\
    				 &$m=3$&0.0843&0.0779&0.0665&0.0689&0.0639&0.0507&0.0451\\
				 \hline
\multirow{4}{*}{$\alpha$-mixing}&AR(2), $N(0, \, 1)$ innov.& 0.0463&0.0200&0.0368&0.0388&0.0423&0.0453&0.0530\\
				&{AR(2) Prod., $N(0, \, 1)$ innov.} &0.0467&0.0434&0.0372&0.0411&0.0387&0.0335&0.0345 \\
   				 &{AR(2), $U[-1, \, 1]$ innov.}&0.0529&0.0251&0.0338&0.0390&0.0396&0.0501&0.0469\\
    				 &{AR(2), $t_{9.5}$ innov.}&0.0474&0.0225&0.0365&0.0402&0.0381&0.0493&0.0483\\
\hline \hline
\end{tabular}
\caption{Monte Carlo simulation results for rejection probabilities for tests of $c_1= 0$, in multiple $m$-dependent and $\alpha$-mixing settings.} \label{tab:cov}
\end{table}

\section{Conclusions}

When the fundamental assumption of exchangeability does not necessarily hold, permutation tests are invalid unless strict conditions on underlying parameters of the problem are satisfied. For instance, the permutation test of $\rho(1) =0$ based on the sample first-order autocorrelation is asymptotically valid only when $\sigma^2$, the marginal variance of the distribution, and $\gamma_1 ^2$, where $\gamma_1^2$ is as defined in (\ref{eq.gam.defn}), are equal. Hence rejecting the null must be interpreted correctly, since rejection of the null with this permutation test does not necessarily imply that the true first-order autocorrelation of the sequence is nonzero. We provide a testing procedure that allows one to obtain asymptotic rejection probability $\alpha$ in a permutation test setting. A significant advantage of this test is that it has the exactness property, absent from the Ljung-Box and Box-Pierce tests, under the assumption of independent and identically distributed sequences, as well as achieving asymptotic level $\alpha$ in a much wider range of settings than the aforementioned tests. An analogous testing procedure permits for asymptotically valid inference in a test of the $k$th order autocorrelation. 

As described in the Introduction, correct implementation of a permutation test is crucial if one is interested in confirmatory inference via hypothesis testing; indeed, proper error control of Type 1, 2 and 3 errors can be obtained for tests of autocorrelations by basing inference on test statistics which are asymptotically pivotal. A framework has been provided for a test of serial lack of correlation in time series data, where tests for $\rho(j) = 0$ are conducted simultaneously for a large number of values of $j$, while maintaining error control with respect to the familywise error rate. 


\addcontentsline{toc}{section}{Bibliography}

\bibliographystyle{apalike}    
\bibliography{perm_test_time_series}

\end{document}


\title{Supplement to\\
Permutation Testing for Dependence in Time Series}

\author{Joseph P. Romano 
  \\
Departments of  Statistics and Economics \\
Stanford University \\
\href{mailto:romano@stanford.edu}{romano@stanford.edu}
\and
Marius Tirlea \\
Department of Statistics \\
Stanford University \\
\href{mailto:mtirlea@stanford.edu}{mtirlea@stanford.edu}
 }

\date{\today}

\maketitle

\section{Introduction}

In this supplement, we provide proofs of the results introduced in the main paper, in addition to statements and proofs of auxiliary results used in the proofs of the lemmas and theorems in the main paper. 

\section{Permutation distribution for $m$-dependent sequences} \label{sec.supp.mdep}

\begin{proof} [{\sc Proof of Theorem \ref{thm.clt}}]

By Theorem 1.1 in \cite{rinott}, itself from \cite{stein:book} page 110, for $Y_i = X_i/\sigma_n$, 

\beq
\begin{aligned}
\left| \pp \left( W_n \leq t \right) - \Phi (t) \right| &\leq 2\ee \left[ \left( \sum_{i=1} ^n \sum_{j \in S_i} \left( Y_i Y_j - \ee Y_i Y_j \right) \right)^2\right]^{1/2} + \\
&+\sqrt{\frac{\pi}{2}} \ee \sum_{i=1} ^n \left| \ee \left[ Y_i \, | \, Y_j : j \notin S_i \right] \right| + \\
&+ 2^{3/4} \pi^{-1/4} \ee \left[ \sum_{i=1} ^n \left| Y_i \right| \left( \sum_{j \in S_i} Y_j \right) ^2 \right]^{1/2} \, \,  .
\end{aligned}
\label{stein.thm1}
\eeq

\noindent We now bound the expectations on the right hand side of (\ref{stein.thm1}). Begin by bounding the second expectation. By definition of the $S_i$,
\beq 
\begin{aligned}
\ee \sum_{i=1} ^n \left| \ee \left[ Y_i \, | \, Y_j : j \notin S_i \right] \right| &= \ee \sum_{i=1} ^n \left| \ee Y_i \right| \\
&= 0 \, \, ,
\end{aligned}
\label{piece.1.clt}
\eeq
since $\ee X_i = 0$ for all $i$.

\noindent We turn our attention to the third expectation. Let $i \in [n]$.

\beq
\begin{aligned}
\ee \left[ \left( \sum_{j \in S_i} Y_j \right)^4 \right] &\leq \ee \left[ \left( \sum_{j \in S_i} \left| Y_j \right| \right)^4 \right] \\
&\leq D^4 \ee \left[ \max_{j \in S_i} Y_j ^4 \right] \, \, \, (x \mapsto x^4 \text{ is increasing on } \rr_{\geq 0}) \\
&\leq D^4 \ee \left[ \sum_{j \in S_i} Y_j ^4 \right] \, \, \, (\text{union bound}) \\
&\leq \frac{D^5 M_4}{\sigma_n ^4} \, \, .
\end{aligned}
\label{piece.2.inter.clt}
\eeq

\noindent Hence

\beq
\begin{aligned}
\ee \left[ \sum_{i=1} ^n \left| Y_i \right| \left(\sum_{j \in S_i} Y_j \right) ^2 \right] &= \sum_{i=1} ^n\ee \left[  \left| Y_i \right| \left(\sum_{j \in S_i} Y_j \right) ^2 \right] \\
&\leq \sum_{i=1} ^n \ee \left[ Y_i ^2\right]^{1/2} \ee \left[ \left( \sum_{j \in S_i} Y_j \right)^4 \right]^{1/2} \, \, \, (\text{by Cauchy-Schwarz}) \\
&\leq \sum_{i=1} ^n \left(\frac{M_2}{\sigma_n ^2} \cdot \frac{D^5 M_4}{\sigma_n ^4} \right)^{1/2} \, \, \, (\text{by  (\ref{piece.2.inter.clt})})\\
&= \frac{n\left(D^5 M_2 M_4\right)^{1/2} }{\sigma_n ^3} \, \, .
\end{aligned}
\label{piece.2.clt}
\eeq

\noindent Lastly, we bound the first term in the expectation. Let $U_{ij} = Y_i Y_j - \ee Y_i Y_j, \, i, j \in [n]$.

\beq
\begin{aligned}
\ee \left[ \left( \sum_{i=1}^n \sum_{j\in S_i} (Y_i Y_j - \ee Y_i Y_j ) \right)^2 \right] &= \ee \left[ \sum_{i=1} ^n \sum_{j \in S_i } U_{ij} \sum_{k=1} ^n \sum_{l \in S_k} U_{kl} \right] \\
&= \sum_{i=1} ^n \sum_{j \in S_i} \sum_{k=1} ^n \sum_{l \in S_k} \ee U_{ij} U_{kl} \, \, .
\end{aligned}
\label{piece.3.inter.clt}
\eeq

\noindent By independence, and the fact that $\ee U_{ij} = 0$, observe that $\ee U_{ij} U_{kl} = 0$ unless either $k$ or $l$ is a member of $S_i \cup S_j$. Note that $\left| S_i \cup S_j \right| \leq 2D$.  It follows that we have at most $4nD^3$ nonzero terms on the right hand side of (\ref{piece.3.inter.clt}). This is because there are $n$ choices for $i$, at most $D$ choices for $j$, at most $2D$ choices for $k$, under the assumption that $k \in S_i \cup S_j$, and at most $D$ choices for $l$, and the same upper bound on nonzero terms if we assume that $l \in S_i \cup S_j$ (noting that $l \in S_k \iff k \in S_l$). Also, by repeated application of Cauchy-Schwarz,

\beq
\begin{aligned}
\left| \ee \left[ (Y_i Y_j - \ee Y_i Y_j)  (Y_k Y_l - \ee Y_kY_l) \right]\right| &\leq
\ee \left[ (Y_i Y_j - \ee Y_i Y_j)^2 \right]^{1/2}  \ee \left[ (Y_k Y_l - \ee Y_k Y_l)^2\right]^{1/2} \\
&\leq \ee \left[ Y_i ^2 Y_j ^2 \right]^{1/2}  \ee \left[ Y_k ^2Y_l ^2\right]^{1/2} \\
&\leq \ee \left[ Y_i ^4 \right]^{1/4} \ee \left[ Y_j ^4 \right]^{1/4} \ee \left[ Y_k ^4 \right]^{1/4} \ee \left[ Y_l ^4 \right]^{1/4} \\
&\leq \frac{M_4}{\sigma_n ^4} \, \, .
\end{aligned}
\label{mbound.inter.3.clt}
\eeq

\noindent Substituting the result of (\ref{mbound.inter.3.clt}) into (\ref{piece.3.inter.clt}), we obtain that 

\beq
\ee \left[ \left( \sum_{i=1}^n \sum_{j\in S_i} (Y_i Y_j - \ee Y_i Y_j ) \right)^2 \right] \leq \frac{4nD^3 M_4}{\sigma_n ^4} \, \, .
\label{piece.3.clt}
\eeq

\noindent Combining (\ref{piece.1.clt}), (\ref{piece.2.clt}), and (\ref{piece.3.clt}), and substituting into (\ref{stein.thm1}):

$$
\begin{aligned}
\left| \pp \left( W_n \leq t \right) - \Phi (t) \right| &\leq 2\left( \frac{4nD^3 M_4}{\sigma_n ^4} \right)^{1/2} + 2^{3/4} \pi^{-1/4} \left( \frac{n\left(D^5 M_2 M_4\right)^{1/2} }{\sigma_n ^3} \right)^{1/2} \\
&= \frac{1}{\sigma_n} \left( 4 \left( \frac{nD^3 M_4}{\sigma_n ^2} \right)^{1/2} + 2^{3/4} \pi^{-1/4} \left( \frac{n}{\sigma_n} \right)^{1/2} \left( D^5 M_2 M_4 \right)^{1/4} \right) \, \, ,
\end{aligned}
$$

\noindent as required.
\end{proof}

\begin{proof} [{\sc Proof of Theorem \ref{rand.dist.mdep}, in the i.i.d. setting}]

Let $\Pi_n, \, \Pi' _n \overset{i.i.d.} \sim \text{Unif}\{S_n \}$, independent of $X_1, \, \dots, \, X_n$. Let $\Pi_n \bf{X}$, $\Pi_n ' \bf{X}$ denote the action of $\Pi_n$ and $\Pi_n '$ on $\mathbf{X} = (X_1, \, \dots, \, X_n)$, respectively. By \cite{tsh} Theorem 15.2.3, it suffices to show that Hoeffding's condition holds, i.e. that $\sqrt{n} \hat{\rho} \left( \Pi_n \mathbf{X} \right)$ and $\sqrt{n} \hat{\rho} \left( \Pi_n '\mathbf{X} \right)$ are asymptotically independent, with common distribution $\Phi$.

Note that, for all $m \in \rr$, 

$$
\left(\sqrt{n} \hat{\rho} \left( \Pi_n \mathbf{X} \right), \, \sqrt{n} \hat{\rho} \left( \Pi_n '\mathbf{X} \right) \right) \overset{d} = \left(\sqrt{n} \hat{\rho} \left( \Pi_n (\mathbf{X}  +m )\right), \, \sqrt{n} \hat{\rho} \left( \Pi_n ' (\mathbf{X} +m)\right) \right) \,   ,
$$

\noindent and so, taking $m = -\mu$, we may assume, without loss of generality, that $\ee X_i = 0$. Similarly, since $\hat{\rho}$ is invariant under the transformation $\mathbf{X} \mapsto \tau \mathbf{X}$ for all $\tau > 0$, we may also assume, without loss of generality, that $\sigma^2 = 1$.

Note also that, for all $\pi_n \in S_n$, 

$$
\begin{aligned}
\bar{X}_{\pi_n} &:= \frac{1}{n} \sum_{i=1} ^n X_{\pi_n ^ {-1} (i)} = \bar{X}_n \\
\hat{\sigma}^2 _{\pi_n} &:= \frac{1}{n} \sum_{i=1} ^n \left( X_{\pi^{-1} _n (i)} - \bar{X}_{\pi_n} \right)^2 = \hat{\sigma}_n ^2 \, \, .
\end{aligned}
$$

\noindent Since $\hat{\sigma }_n ^2 \to 1$ almost surely by the strong law of large numbers (under the assumption of finite fourth moments), it follows, by Slutsky's theorem for randomization distributions (see Theorem 5.2 in \cite{chung.perm}), that $ \left(\sqrt{n} \hat{\rho} \left( \Pi_n \mathbf{X} \right), \, \sqrt{n} \hat{\rho} \left( \Pi_n '\mathbf{X} \right) \right) \overset{d} \to N(0, \, I_2)$  if and only if 

\beq
\frac{1}{\sqrt{n}}\sum_{i=1} ^{n-1}
\begin{pmatrix}
(X_{\Pi ^{-1}_n(i)} - \bar{X}_n)(X_{\Pi ^{-1}_n(i+1)} - \bar{X}_n)\\
(X_{\Pi '^{-1}_n(i)} - \bar{X}_n)(X_{\Pi '^{-1}_n(i+1)} - \bar{X}_n )
 \end{pmatrix}
\overset{d} \to N(0, \, I_2) \,\, .
\label{slutsky.iid}
\eeq

\noindent Also, by the moment assumptions, $\bar{X}_n = O_p (1/\sqrt{n})$ by the Central Limit Theorem, and $\sum_{i=1} ^{n-1} X_{\Pi_n ^{-1}(i)} = n\bar{X}_n + O_p (1)$. Hence

$$
\begin{aligned}
\frac{1}{\sqrt{n}}\sum_{i=1} ^{n-1} (X_{\Pi ^{-1}_n(i)} - \bar{X}_n)(X_{\Pi ^{-1}_n(i+1)} - \bar{X}_n) &= \frac{1}{\sqrt{n}} \sum_{i=1} ^{n-1} X_{\Pi ^{-1}_n(i)}X_{\Pi ^{-1}_n(i+1)} - \sqrt{n} \bar{X}_n ^2 + o_p (1) \\
&= \frac{1}{\sqrt{n}} \sum_{i=1} ^{n-1} X_{\Pi ^{-1}_n(i)}X_{\Pi ^{-1}_n(i+1)} + o_p(1) \, \, ,
\end{aligned}
$$

\noindent and similarly for $\Pi_n ' \mathbf{X}$. Hence the convergence in (\ref{slutsky.iid}) holds if and only if 

\beq
\left( \frac{1}{\sqrt{n}}\sum_{i=1} ^{n-1} X_{\Pi ^{-1}_n(i)} X_{\Pi ^{-1}_n(i+1)} , \, \frac{1}{\sqrt{n}}\sum_{i=1} ^{n-1} X_{\Pi '^{-1}_n(i)} X_{\Pi '^{-1}_n(i+1)} \right) \overset{d} \to N(0, \, I_2) \, \, .
\label{clt.iid.2d.inter}
\eeq

\noindent Since the $X_i$ are i.i.d.,we observe that $\Pi_n \mathbf{X}$ and $\mathbf{X}$ are identically distributed. Hence it is easy to see, by relabelling, that $(\Pi_n \mathbf{X}, \, \Pi_n ' \mathbf{X}) \overset{d} =   (\mathbf{X}, \, \Pi_n '\Pi_n^{-1} \mathbf{X})$. Also, $\Pi_n ^{-1} \overset{d} = \Pi_n \overset{d} = \Pi_n ' \overset{d} = \Pi_n ' \Pi_n ^{-1}$, and so the convergence in (\ref{clt.iid.2d.inter}) holds if and only if 

\beq
\left( \frac{1}{\sqrt{n}}\sum_{i=1} ^{n-1} X_i X_{i+1} , \, \frac{1}{\sqrt{n}}\sum_{i=1} ^{n-1} X_{\Pi _n(i)} X_{\Pi_n(i+1)} \right) \overset{d} \to N(0, \, I_2) \, \, .
\label{clt.iid.2d}
\eeq

\noindent The  sequence $\{X_i X_{i+1}, \, i \in [n-1] \}$ has uniformly bounded fourth moment, by the moment condition and independence of the $X_i$. Observe that the sequence $\{X_i X_{i+1}, \, i \in [n-1] \}$ is 1-dependent, with mean 0, variance 1, finite 4th moment, and first-order autocorrelation 0. Therefore, by the $m$-dependent central limit theorem and Slutsky's theorem (rescaling from $\sqrt{n-1}$ to $\sqrt{n}$), as $n \to \infty$, 

\beq
\frac{1}{\sqrt{n}} \sum_{i=1} ^{n-1} X_i X_{i+1} \overset{d} \to N(0,1) \, \, .
\label{clt.iid.1d}
\eeq

\noindent A similar result holds for the corresponding sum of the sequence $\Pi_n \mathbf{X}$. By the Cram\'er-Wold device, (\ref{clt.iid.2d}) holds if and only if, for all $a, \, b \in \rr$, 

\beq
\frac{1}{\sqrt{n}} \sum_{i=1} ^{n-1} a X_i X_{i+1} + \frac{1}{\sqrt{n}} \sum_{i=1} ^{n-1} bX_{\Pi_n (i)} X_{\Pi_n (i+1)} \overset{d} \to N(0, \, a^2 + b^2) \, \, .
\label{cramer.rao.iid}
\eeq

\noindent Let $Y_i = \frac{1}{\sqrt{n}} a X_i X_{i+1} + \frac{1}{\sqrt{n}} b X_{\Pi_n (i)} X_{\Pi_n (i+1) }$,  $i \in [n-1]$. Let $S_i = \{ j: Y_i \text{ and } Y_j \text{ are not independent}\}$. We observe that the $Y_i$ are identically distributed (by symmetry), and, uniformly in $i$,

$$
\begin{aligned}
M_4 := \ee \left[ Y_i ^4 \right] &= \ee \left[ \left( \frac{1}{\sqrt{n}} \sum_{i=1} ^{n-1} a X_i X_{i+1} + \frac{1}{\sqrt{n}} \sum_{i=1} ^{n-1} bX_{\Pi_n (i)} X_{\Pi_n (i+1)} \right)^4 \right] \\
&\leq 2^3 \ee \left[ \frac{a^4}{n^2} X_i ^4 X_{i+1} ^4 + \frac{b^4}{n^2} X_{\Pi_n (i)} ^4 X_{\Pi_n(i+1)} ^4 \right] \\
&\leq \frac{C}{n^2} \, \, ,
\end{aligned}
$$

\noindent for some $C$ depending only on $a$ and $b$. Similarly, for $M_2 := \ee Y_i ^2$, $M_2 \leq C/n$, for some $C$ depending only on $a$ and $b$ (we may choose $C$ to be the same for both $M_2$ and $M_4$). Also, $\ee Y_i = 0$.

Now, condition on $\Pi_n = \pi$, where we leave the dependence of $\pi$ on $n$ implicit. In the following, expectations are taken conditional on $\Pi_n = \pi$. Consider the vectors

$$
\begin{aligned}
v_1 = \begin{pmatrix}
X_i X_{i+1} \\
X_{\pi(i) } X_{\pi(i+1)}
\end{pmatrix}, \, \, \, \, 
v_2 = \begin{pmatrix}
X_j X_{j+1} \\
X_{\pi(j) } X_{\pi(j+1)}
\end{pmatrix} \, \, .
\end{aligned} 
$$

\noindent $v_1$ and $v_2$ are independent unless, for 

$$
\begin{aligned}
A_i &:=  \{ i-1, \, i, \, i+1, \, i+2, \, \pi(i) - 1, \, \pi(i), \, \pi(i) +1, \, \pi(i+1) - 1, \, \pi(i+1), \, \pi(i+1) +1\} \\
B_j &:= \{j, \, j+1, \, \pi(j), \, \pi(j+1) \} \, \,  ,
\end{aligned}
$$

\noindent $A_i \cap B_j \neq \emptyset$. There are at most $\left| A_i \right|$ choices for each element of $B_j$, and so there are at most $\left| A_i \right| ^4 \leq 10^4$ vectors of the form $v_2$ which are not independent of $v_1$. It follows that, for all $i$, $\left| S_i \right| \leq 10^4$, which is, in particular, finite.

We may now apply Theorem \ref{thm.clt}, noting that, by Slutsky's theorem, taking the scaling to be $\sqrt{n}$ or $\sqrt{n-1}$ is equivalent. By Remarks \ref{rem.unif.clt} and \ref{rem.clt.clt}, it now suffices to show that $\sigma_n ^2$, as defined in Theorem \ref{thm.clt}, satisfies $\sigma_n ^2 = a^2 + b^2 + o(1)$. By (\ref{clt.iid.1d}), and the corresponding result for $\pi \mathbf{X}$, 

\beq
\begin{aligned}
\sigma_n ^2 = a^2 + b^2 + o(1) + \frac{2ab}{n} \sum_{i=1} ^{n-1} \sum_{j=1} ^{n-1} \cov \left(X_i X_{i+1}, \, X_{\pi(j)} X_{\pi(j+1)} \right) \, \, .
\end{aligned}
\label{cov.iid.1d}
\eeq

\noindent We claim that the sum of covariances in (\ref{cov.iid.1d}) is $o(1)$. To this end, we have that 

$$
\cov \left(X_i X_{i+1}, \, X_{\pi(j)} X_{\pi(j+1)} \right) = \begin{cases}
1, \, \, \, \text{if } \{\pi(j), \, \pi(j+1) \} = \{i, \, i+1 \} \, \, , \\
0, \, \, \, \text{otherwise} \, \, .
\end{cases}
$$

\noindent Let $M_n$ be the number of pairs of indices $(i, \, j)$ such that $\{\pi(j), \, \pi(j+1) \} = \{i, \, i+1 \}$, $N_n$ be the number of pairs of indices $(i, \, j)$ such that $i = \pi(j)$ and $i+1 = \pi(j+1)$, and let $N_n '$ be the number of pairs of indices $(i, \, j)$ such that $i = \pi(j+1)$ and $i+1 = \pi(j)$. Note that $M_n = N_n + N_n '$. The following argument is taken from \cite{ritzwoller}. 

If $M_n$ were uniformly bounded, it would be the case that $\sigma_n ^2 = a^2 +b ^2 + o(1)$. However, $M_n$ depends on the values of the realizations $\pi$ as $n \to \infty$. We argue as follows: for all $n$, 

$$
\begin{aligned}
\ee N_n &= \sum_{i = 1} ^n \sum_{j=1} ^n \pp (\Pi_n (j) = i, \, \Pi_n (j+1) = i+1) \\
&= \frac{n^2}{n(n-1)} \leq 2 \, \,,
\end{aligned}
$$

\noindent and similarly for $N_n '$. Therefore $\{M_n, \, n \in \nn\}$ are nonnegative random variables with uniformly bounded expectation, and so are tight. Hence, for any subsequence $\{n_j\} \subset \mathbb{N}$, there exists a further subsequence $\{n_{j_k}\}$ such that $M_{n_{j_k}}$ converges in distribution. By Skorokhod's almost sure representation theorem, there exist random variables $\tilde{M}_{n_{j_k}}$ with the same distribution as the $M_{n_{j_k}}$ which converge almost surely, to $\tilde{M}$, say.

Based on $\tilde{M} = m$, construct $\tilde{\Pi}_n$ according to the conditional distribution of $\Pi_n \, | \, M_n = m$, such that, unconditionally, $\tilde{\Pi}_n$ is uniform on $S_n$. Then, along a subsequence, we may use the argument above as if $M$ were finite. In summary, given a subsequence $\{n_j\} \in \mathbb{N}$, there exists a further subsequence $n' = \{n_{j_k}\}$ such that, uniformly in $t$,

$$
\pp\left( \sum_{i=1} ^{n_{j_k}-1} Y_i \leq t \, | \, \Pi_{n_{j_k}} \right) \overset{{a.s.}} \to \Phi \left(\frac{t}{\sqrt{a^2 + b^2}} \right) \, \, .
$$

\noindent Therefore, unconditionally, by dominated convergence,

\beq
\pp\left( \sum_{i=1} ^{n_{j_k}-1} Y_i \leq t \right)  \to \Phi \left(\frac{t}{\sqrt{a^2 + b^2}} \right) \, \, .
\label{dct.iid}
\eeq

\noindent Since, given any subsequence, the result holds along a further subsequence, the limit in (\ref{dct.iid}) holds along the original subsequence. It follows that Hoeffding's condition holds, and so the result follows.
\end{proof}

The following Lemmas are required in the proof of Theorem \ref{rand.dist.mdep}, in the $m$-dependent setting.

\begin{lem} Let $X_1, \, \dots, \, X_n$ be random variables such that $\ee X_i = \nu_i $, $\var(X_i) = 1$, and, the $X_i$ have uniformly bounded $4$th moment.

Let $S_i = \left\{ j \in [n] : X_i \text{ and } X_j \text{ are not independent}\right\}$. Suppose $\left| S_i \right| \leq D < \infty$ for all $i$. Let $\sigma_n ^2 = \ee \left[ \sum_{i=1} ^n X_i \sum_{j \in S_i} X_j \right] = \var\left( \sum_{i=1} ^n X_i \right)$, $\mu_n = \sum_{i=1} ^n \nu_i$, and $W_n = \sum_{i=1} ^n X_i /\sigma_n$. If $\sigma_n ^2 \to 1$ and $\mu_n \to 0$, then, as $n \to \infty$,

$$
\sup_{t \in \rr} \left| \pp (W_n \leq t ) -\Phi (t) \right| \to 0 \, \, .
$$

\label{lem.littleo}
\end{lem}

\begin{proof} Let $Y_i = X_i - \nu_i$. We may now apply the result of Theorem \ref{thm.clt} and Remarks \ref{rem.unif.clt} and \ref{rem.clt.clt}:

$$
\sup_{t \in \rr} \left| \pp \left(\frac{W_n - \mu_n}{\sigma_n} \leq t \right) -\Phi (t) \right| = O\left( n^{-1/4} \right) \, \, .
$$

\noindent Reparametrizing, it follows that 

$$
\sup_{t \in \rr} \left| \pp \left(W_n \leq t \right) -\Phi \left(t \right) \right| = \sup_{t \in \rr} \left| \pp \left(W_n \leq t \right) -\Phi \left(\frac{t-\mu_n}{\sigma_n} \right) \right| \, \, .
$$

\noindent Now, since $\lVert \cdot \rVert_\infty$ is a norm on $\rr$, we have that, by the triangle inequality,

\beq
\begin{aligned}
&\sup_{t \in \rr} \left| \pp \left(W_n \leq t \right) -\Phi \left(t \right) \right| \\
\leq &\sup_{t \in \rr} \left| \pp \left(W_n \leq t \right) -\Phi \left(\frac{t-\mu_n}{\sigma_n} \right) \right| + \sup_{t \in \rr} \left| \Phi (t) -\Phi \left(\frac{t-\mu_n}{\sigma_n} \right) \right| \, \, .
\end{aligned}
\label{cdf.bound}
\eeq

\noindent We now bound the second term in (\ref{cdf.bound}). 

$$
\begin{aligned}
\sup_{t \in \rr} \left| \Phi (t) -\Phi \left(\frac{t-\mu_n}{\sigma_n} \right) \right| &= \sup_{t \in \rr} \left| \int_{(t-\mu_n)/\sigma_n} ^ t \phi(x) \d x\right| \\
&\leq \sup_{t \in \rr} \left| t - \frac{t - \mu_n}{\sigma_n} \right| \exp\left( -\frac{1}{2} \left( t^2 \land \left(\frac{t- \mu_n}{\sigma_n} \right)^2 \right) \right)\\
&\to 0 \, \, .
\end{aligned}
$$

\noindent The result follows.
\end{proof}

\begin{lem} Let $X_1, \, \dots, \, X_n$ be a stationary, $m$-dependent time series, with mean 0, $k$th-order autocorrelation $\rho_k$, $k \in [m]$, variance $\sigma^2 = 1$, and finite 8th moment. Let $\Pi_n \sim \text{Unif}\{S_n\}$.
Let

$$
W_n = \frac{1}{\sqrt{n}} \sum_{i=1} ^{n-1} X_{\Pi_n (i)} X_{\Pi_n (i+1)}  \, \, .
$$

\noindent Then, as $n \to \infty$, 

$$
\sup_{t \in \rr} \left| \pp( W_n \leq t) - \Phi (t) \right| \to 0 \, \, .
$$

\label{clt.mdep}
\end{lem}

\begin{proof} By application of Cauchy-Schwarz and the given properties of $\{X_i, \, i\in [n]\}$, note that $\ee X_{\Pi_n (i) } ^4 X_{\Pi_n (i+1) } ^4$ is uniformly bounded. In addition, conditional on $\Pi_n = \pi$ (where we leave the dependence of $\pi$ on $n$ implicit),

$$
\begin{aligned}
\ee X_{\pi (i) }  X_{\pi (i+1) } &=\sum_{j = 1}^ n \sum_{k \neq j} 1\{\pi(i) = j, \, \pi(i+1) = k\} \ee X_j X_k \\
&= \sum_{j = 1}^ n \sum_{k \neq j} \sum_{l = 1} ^m \rho_l  1\{\pi(i) = j, \, \pi(i+1) = k, \, \left| j - k \right| = l\} \, \, .
\end{aligned}
$$

\noindent Let $Y_i = \frac{1}{\sqrt{n}} \left(X_{\pi(i)} X_{\pi(i+1)} - \ee X_{\pi(i)} X_{\pi(i+1)} \right)$. Observe that, in the setting of Theorem \ref{thm.clt}, $\left| S_i \right| \leq 4m^2$ (since at most $2m$ of the $X_j$ are not independent of $X_{\pi(i)}$, and similarly for $X_{\pi(i + 1)}$), so we may apply Theorem \ref{thm.clt} and Remarks \ref{rem.unif.clt} and \ref{rem.clt.clt}: letting $Z_n = \sum_{i=1} ^{n-1} Y_i/\sigma_n$, where $\sigma_n ^2 = \var \left(\sum_{i=1} ^{n-1} Y_i \right)$,  

$$
\sup_{t \in \rr} \left| \pp( Z_n \leq t) - \Phi (t) \right| = O(n^{-1/4})
$$

\noindent Let $\mu_n = \frac{1}{\sqrt{n}}\sum_{i=1} ^{n-1} \ee X_{\pi(i)} X_{\pi(i+1)}$.

$$
\begin{aligned}
\sigma_n ^2 &=\frac{1}{n} \sum_{i=1} ^{n-1} \ee X_{\pi (i) } ^2 X_{\pi(i+1)} ^2 + \frac{1}{n} \sum_{i = 1}^{n-1} \sum_{j\neq i} \cov \left(X_{\pi(i)} X_{\pi (i+1) }, \, X_{\pi(j)} X_{\pi (j +1) }\right)
\end{aligned}
$$

\noindent We wish to show that $\mu_n \to 0$, and $\sigma_n ^2 \to 1$. Unfortunately, this does not occur for an arbitrary $\pi \in S_n$. However, we may argue as in the proof of the i.i.d. setting of Theorem \ref{rand.dist.mdep}. Let $K$ be such that, for all $i$, $j$, $k$ and $l$, 

$$
\begin{aligned}
\left| \ee X_i X_j \right| &\leq K \\
\left| \ee X_i ^2  X_j ^2 \right| &\leq K \\
\left| \cov \left(X_i X_j, \, X_k X_l \right) \right| &\leq K \, \, .
\end{aligned}
$$ 

\noindent Such a $K$ exists by the moment assumptions on the $X_i$. Note that, for $i < j < k < l$, $\cov( X_i X_j, \, X_k X_l) = 0$ unless $j \leq i + m, \, l \leq k + m$ (since otherwise one of the terms in the product is independent of the others). This also holds for any permutation of $i$, $j$, $k$, and $l$, by the same argument. Let 

\beq
\begin{aligned}
A_n &= \left\{i : \left| \Pi_n (i) - \Pi_n (i+1) \right| \leq m \right\} \\
\Gamma_n &= \left\{ \{i, \, j, \, k, \, l \} : i< j <k<l, \, j \leq i + m, \, l \leq k + m \right\} \\
B_n &= \left\{ \{i, \, j \} : \{\Pi_n (i), \, \Pi_n (i+1), \, \Pi_n (j), \, \Pi_n (j+1) \} \in \Gamma_n \right\} \, \, .
\end{aligned}
\label{set.bounds}
\eeq

\noindent Note that (still conditional on $\Pi_n = \pi$), $\ee X_{\pi(i)} X_{\pi(i+1)} = 0$ unless $i \in A_n$, $\ee X ^2 _{\pi(i)} X^2 _{\pi(i+1)} = 1$ unless $i \in A_n$, and $\cov \left(X_{\pi(i)} X_{\pi (i+1) }, \, X_{\pi(j)} X_{\pi (j +1) }\right)= 0$ unless $\{i, \, j\} \in B_n$. Note that all expectations and covariances are uniformly bounded. Unconditionally, we have that 

\beq
\begin{aligned}
\ee \left|A_n\right| &= \sum_{i=1} ^n \sum_{j=1} ^n \pp \left( \Pi_n (i) = j ,\, \left| \Pi_n (i+1) - j \right| \leq m \right) \\
&\leq \sum_{i=1} ^n \sum_{j=1} ^n \frac{1}{n} \cdot \frac{2m}{n} \\
&= 2m \, \, .
\end{aligned}
\label{an.bound}
\eeq

\noindent Also, by assuming without loss of generality that $\Pi_n (i) < \Pi_n (i+1) < \Pi_n (j) < \Pi_n (j+1)$ and then multiplying by the appropriate number of possible permutations of these values,

\beq
\begin{aligned}
\ee \left| B_n\right| &\leq 24 \sum_{i=1} ^n \sum_{j=1} ^n \sum_{k=1} ^n \sum_{l=1} ^n \pp\left( \Pi_n (i) = j, \, j < \Pi_n (i+1) \leq j+m, \, \Pi_n (k) = l, \, l < \Pi_n (k +1) \leq l + m\right) \\
&\leq \frac{24}{n^2} \sum_{i=1} ^n \sum_{j=1} ^n \sum_{k=1} ^n \sum_{l=1} ^n \frac{m}{n} \cdot \frac{m}{n} \\
&= 24m^2 \, \, .
\end{aligned}
\label{bn.bound}
\eeq

\noindent In particular, letting $M_n = \left|A_n \right| + \left| B_n\right|$, we note that $\{M_n\}$ are nonnegative random variables with uniformly bounded expectation. Hence, proceeding as in the proof of the i.i.d. setting of Theorem \ref{rand.dist.mdep}, the result follows, by Lemma \ref{lem.littleo} (noting that scaling by $\sqrt{n}$ or $\sqrt{n-1}$ is equivalent, by Slutsky's theorem). 
\end{proof}

\begin{proof} [{\sc Proof of Theorem \ref{rand.dist.mdep}, in the $m$-dependent setting}]

By the same argument as in the proof of the i.i.d. setting of Theorem \ref{rand.dist.mdep}, we may assume, without loss of generality, that $\mu = 0$ and $\sigma^2 = 1$.

Note that, both $\bar{X}_n$ and $\hat{\sigma}^2 _n$ are invariant under permutation. Also, by the $m$-dependent CLT, since $\{X_i X_{i+1},\, i\in \nn\}$ is a stationary $(m+1)$-dependent time series with finite mean and variance,  

$$
\bar{X}_n = O_p \left(\frac{1}{\sqrt{n}} \right) \, \, .
$$

\noindent Now, note that, by definition of $m$-dependence, for all $j \in \{1, \, \dots, \, m +1 \}$, $\{X_{j + k(m+1)}, \, {k \geq 0}\}$ are independent. Hence, letting $\hat{\sigma}_{n, \, j} ^2 = \frac{1}{n}\sum_{j + k(m+1) \leq n} X_{j + k(m+1)}^2 $, we have that, by the strong law of large numbers, for all $j$,

$$
\hat{\sigma}_{n, \, j} ^2 \overset{a.s.} \to \frac{1}{m+1} \, \, .
$$

\noindent Summing over $j$, and combining with the asymptotic distribution of $\bar{X}_n$, it follows that 

$$
\hat{\sigma}_n ^2 \overset{p} \to 1 \, \, .
$$	

\noindent Hence, by the same argument as in the proof of the i.i.d. setting of Theorem \ref{rand.dist.mdep}, it suffices to show that, for $\Pi_n, \, \Pi_n ' \sim \text{Unif}\{S_n \}$, independent of each other and of the $X_i$ (noting that $(\Pi_n, \, \Pi_n ' ) \overset{d} = (\Pi_n ^{-1}, \, (\Pi_n ')^{-1})$),

\beq
\left( \frac{1}{\sqrt{n}}\sum_{i=1} ^{n-1} X_{\Pi_n(i)} X_{\Pi_n(i+1)} , \, \frac{1}{\sqrt{n}}\sum_{i=1} ^{n-1} X_{\Pi '_n(i)} X_{\Pi'_n(i+1)} \right) \overset{d} \to N(0, \, I_2) \, \, .
\label{mdep.2d.clt}
\eeq

\noindent By the Cram\'er-Wold device, the convergence in (\ref{mdep.2d.clt}) holds if and only if, for all $a, \, b \in \rr$, 

\beq
\frac{1}{\sqrt{n}}\sum_{i=1} ^{n-1}a X_{\Pi_n(i)} X_{\Pi_n(i+1)} + \frac{1}{\sqrt{n}}\sum_{i=1} ^{n-1} bX_{\Pi '_n(i)} X_{\Pi'_n(i+1)}  \overset{d} \to N(0, \, a^2 + b^2) \, \, .
\label{mdep.1d.clt}
\eeq

Let $Y_i = \frac{1}{\sqrt{n}} a X_{\Pi_n (i)} X_{\Pi_n(i+1)} + \frac{1}{\sqrt{n}} b X_{\Pi'_n (i)} X_{\Pi'_n (i+1) }$,  $i \in [n-1]$. Let 

$$
S_i = \{ j: Y_i \text{ and } Y_j \text{ are not independent}\} \, \, .
$$

\noindent Observe that the $Y_i$ are identically distributed (by symmetry), and, uniformly in $i$,

$$
\begin{aligned}
M_4 := \ee \left[ Y_i ^4 \right] &= \ee \left[ \left( \frac{1}{\sqrt{n}} \sum_{i=1} ^{n-1} a X_{\Pi_n (i)} X_{\Pi_n (i+1)} + \frac{1}{\sqrt{n}} \sum_{i=1} ^{n-1} bX_{\Pi_n '(i)} X_{\Pi_n '(i+1)} \right)^4 \right] \\
&\leq 2^3 \ee \left[ \frac{a^4}{n^2} X_{\Pi_n (i)} ^4X_{\Pi_n (i+1)}^4 + \frac{b^4}{n^2} X_{\Pi'_n (i)} ^4 X_{\Pi'_n(i+1)} ^4 \right] \\
&\leq \frac{C}{n^2} \, \, ,
\end{aligned}
$$

\noindent for some $C$ depending only on $a$ and $b$. Similarly, for $M_2 := \ee Y_i ^2$, $M_2 \leq C/n$, for some $C$ depending only on $a$ and $b$ ($C$ may be chosen to be the same for both $M_2$ and $M_4$).

Now, condition on $\Pi_n = \pi, \, \Pi'_n = \pi'$, where we leave the dependence of $\pi$ and $\pi'$ on $n$ implicit. In the following, expectations are taken conditional on $\Pi_n = \pi, \, \Pi_n ' = \pi'$. Consider the vectors

$$
\begin{aligned}
v_1 = \begin{pmatrix}
X_{\pi(i)} X_{\pi(i+1)} \\
X_{\pi'(i) } X_{\pi'(i+1)}
\end{pmatrix}, \, \, \, \, 
v_2 = \begin{pmatrix}
X_{\pi(j) }X_{\pi(j+1)} \\
X_{\pi'(j) } X_{\pi'(j+1)} \, \, 
\end{pmatrix} \, \, .
\end{aligned}
$$

\noindent $v_1$ and $v_2$ are independent unless, for 

$$
\begin{aligned}
A_i &:=  \bigcup_{r = 0} ^1\{k: \left| k - \pi(i+ r) \right| \leq m \} \cup \bigcup_{r = 0} ^1\{k: \left| k - \pi'(i+ r) \right|   \leq m \}\\
B_j &:= \{\pi(j), \, \pi(j+1), \, \pi'(j), \, \pi'(j+1) \} \, \, ,
\end{aligned}
$$

\noindent $A_i \cap B_j \neq \emptyset$. There are at most $\left| A_i \right|$ choices for each element of $B_j$, and so there are at most $\left| A_i \right| ^4 \leq (4(2m+1))^4$ vectors of the form $v_2$ which are not independent of $v_1$. It follows that, for all $i$, $\left| S_i \right| \leq 256(2m+1)^4$, which is, in particular, finite.

We may now apply Theorem \ref{thm.clt} to $Z_i = Y_i - \ee Y_i$, noting that, by Slutsky's theorem, taking the scaling to be $\sqrt{n}$ or $\sqrt{n-1}$ is equivalent. By Remarks \ref{rem.unif.clt} and \ref{rem.clt.clt}, it now suffices to show that $\sigma_n ^2$, as defined in Theorem \ref{thm.clt}, satisfies $\sigma_n ^2 = a^2 + b^2 + o(1)$. By Lemma \ref{clt.mdep}, and the corresponding result for $\pi' \mathbf{X}$, 

\beq
\begin{aligned}
\sigma_n ^2 = a^2 + b^2 + o(1) + \frac{2ab}{n} \sum_{i=1} ^{n-1} \sum_{j=1} ^{n-1} \cov \left(X_{\pi(i)} X_{\pi(i+1)}, \, X_{\pi'(j)} X_{\pi'(j+1)} \right) \, \, .
\end{aligned}
\label{cov.mdep.1d}
\eeq

\noindent By the moment assumptions, we have that $\cov \left(X_{\pi(i)} X_{\pi(i+1)}, \, X_{\pi'(j)} X_{\pi'(j+1)} \right)$ is uniformly bounded (by $K$, as defined in the proof of Lemma \ref{clt.mdep}).

We claim that the sum of covariances in (\ref{cov.mdep.1d}) is $o(1)$. Note that, unless $\{\pi(i), \, \pi(i+1), \, \pi'(j), \, \pi'(j+1) \} \in \Delta_n$, where

$$
\Delta_n = \{ \{i, \, j, \, k, \, l\} : i \leq j \leq k \leq l, \, j \leq i +m, \, l \leq k + m\} \, \, ,
$$

\noindent we will have $\cov \left(X_{\pi(i)} X_{\pi(i+1)}, \, X_{\pi'(j)} X_{\pi'(j+1)} \right) = 0$, since at least one of the terms will be independent of all the others. Let 

$$
C_n =  \left\{ \{i, \, j \} : \{\Pi_n (i), \, \Pi_n (i+1), \, \Pi_n (j), \, \Pi_n (j+1) \} \in \Delta_n \right\} \, \, .
$$

\noindent Conditional on $\Pi_n, \Pi_n '$, it is the case that $C_n$ is not necessarily finite. However, by the same argument as for $B_n$ in the proof of Lemma \ref{clt.mdep},

$$
\begin{aligned}
&\ee \left|C_n\right|  \\
&\leq 24 \sum_{i, \, j, \, k, \, l \in [n]} \pp\left( \Pi_n (i) = k, \, k \leq \Pi_n (i+1) \leq k+m, \, \Pi_n ' (j) = l, \, l \leq \Pi_n ' (j +1) \leq l + m\right) \\
&\leq \frac{24}{n^2}  \sum_{i, \, j, \, k, \, l \in [n]}  \frac{m}{n} \cdot \frac{m}{n} \\
&= 24m^2 \, \, ,
\end{aligned}
$$

\noindent and so the random variables $\{\left|C_n\right|\}$ are nonnegative and have uniformly bounded expectation. By the same argument as in the proof of the i.i.d. setting of Theorem \ref{rand.dist.mdep}, (\ref{mdep.1d.clt}) holds, and so the result follows.
\end{proof}

\begin{proof} [{\sc Proof of Theorem \ref{thm.unperm.rho.lim}}] We may assume, without loss of generality, that $\ee X_1 = 0$. For $i \in [n-1]$, let 

$$
Y_i = a \left(X_i -\bar{X}_n \right)\left(X_{i+1} -\bar{X}_n \right) +b \left(X_i -\bar{X}_n \right)^2 \, \,.
$$

\noindent Let 

$$
Z_i = a X_i X_{i+1} + b X_i ^2 \, \, .
$$

\noindent Note that 

$$
\begin{aligned}
Y_i &= Z_i + \left[ -\bar{X}_n (aX_i + 2b X_i + aX_{i+1} ) + (a+b) \bar{X}_n ^2 \right]\,\, .
\end{aligned}
$$

\noindent It follows that 

$$
\begin{aligned}
\sum_{i=1} ^{n-1} Y_i &= \sum_{i=1} ^{n-1} Z_i + \left[ - \bar{X}_n \cdot 2(a+b)\left(n\bar{X}_n + O_p(1)\right) + (a+b) n\bar{X}_n ^2 \right] \\
&=\sum_{i=1} ^{n-1} Z_i -{n}(a+b) \bar{X}_n ^2 + o_p \left(\sqrt{n}\right) \,\, .
\end{aligned}
$$

\noindent By the $m$-dependent central limit theorem (in (i)), and Ibragimov's central limit theorem (in (ii); see \cite{ibragimov}) we have that 

$$
\bar{X}_n = O_p \left(\frac{1}{\sqrt{n}} \right) \, \, ,
$$

\noindent and so

$$
\sum_{i=1} ^{n-1} Y_i = \sum_{i=1} ^{n-1}  Z_i + o_p (\sqrt{n}) \,\,.
$$

\noindent Note that the sequence $\{Z_i, \, i \in [n-1]\} $ is, in case (i), $(m+1)$-dependent, with finite second moments, by the triangle inequality and the Cauchy-Schwarz inequality, and in case (ii), is $\alpha$-mixing, with $\alpha$-mixing coefficients satisfying 

$$
\alpha_Z (i) = \alpha_X (i-1) \, \, ,
$$

\noindent and, by the same argument as in case i), has finite $(2+ \delta)$ moments. Also, in both cases, the sequence $\{Z_i\}$ is strictly stationary, and 

$$
\ee[Z_i] = {a}\rho_1 \sigma^2 + {b} \sigma^2 \, \, .
$$

\noindent Also, note that 

$$
\frac{1}{\sqrt{n}}\sum_{i=1} ^{n-1} Y_i = a \sqrt{n} \hat{\rho}_n \hat\sigma_{n}^2 + b \sqrt{n} \hat{\sigma}_n ^2 + o_p (1) \, \, .
$$

\noindent Hence, by the $m$-dependent central limit theorem and Ibragimov's central limit theorem and Slutsky's theorem (to rescale from $\sqrt{n-1}$ to $\sqrt{n}$), we have that

$$
\frac{1}{\sqrt{n}}\sum_{i=1} ^{n-1}  Y_i - \sqrt{n} \left(a\rho_1 \sigma^2 + b \sigma^2\right) \overset{d} \to N\left(0, \, \tilde{\sigma} ^2 \right) \, \, ,
$$

\noindent where

$$
\begin{aligned}
\tilde{\sigma} ^2 &= \var\left(a X_1 X_{2} + bX_1 ^2 \right) + 2\sum_{i \geq 2} \cov \left( a X_1 X_{2} + bX_1 ^2, \, a X_i X_{i+1} + bX_i ^2 \right) \\
&= a^2 \tau_1 ^2 + b^2 \kappa^2 + 2ab \nu_1 \, \, .
\end{aligned}
$$

\noindent It follows that, by the Cram\'er-Wold device,

$$
\sqrt{n} \begin{pmatrix}
\hat{\rho}_n\sigma_n ^2  - \rho_1 \sigma^2 \\
\hat{\sigma}_n ^2 - \sigma^2 
\end{pmatrix} \overset{d} \to N\left(0, \, \begin{pmatrix}
\tau^2 & \nu_1 \\
\nu_1 & \kappa^2 \end{pmatrix}\right) \, \, .
$$

\noindent We may now apply the delta method, with the function $h(x, \, y) = x/y$, to conclude that 

$$
\sqrt{n} \left(\hat{\rho}_n - \rho \right) \overset{d} \to N(0, \, \gamma_1^2) \, \, ,
$$

\noindent as required.
\end{proof}

\begin{proof} [{\sc Proof of Theorem \ref{rand.dist.mdep.stud}}]

We begin by showing the result in (i). By Theorem \ref{thm.unperm.rho.lim}, we have that 

\beq
\sqrt{n} \left(\hat{\rho}_n - \rho \right)\overset{d} \to N\left(0, \, \gamma_1^2\right) \, \, .
\label{clt.mdep.unst}
\eeq

\noindent By definition of $m$-dependence, it follows that $(X_1, \, X_2)$ and $(X_{i+1}, \, X_{i+2})$ are independent for $i \geq m+2$. Hence

$$
\begin{aligned}
\tau^2 &= \var(X_1 X_2) + 2\sum_{i=1} ^{m+1} \cov(X_1 X_2, \, X_{i+1} X_{i+2} ) \\
\kappa^2 &= \var \left(X_1 ^2 \right) + 2\sum_{k=1} ^{m+1} \cov\left(X_1^2, \, X_k^2 \right) \\
\nu_1 &= \cov \left( X_1 X_2, \, X_1 ^2 \right) + \sum_{k =2 }^{m+1} \cov\left( X_1^2, \, X_{k} X_{k+1} \right) + \sum_{k = 2}^{m+2} \cov \left(X_1 X_2,\, X_k ^2 \right) \, \, .
\end{aligned}
$$

\noindent By the proof of Theorem \ref{rand.dist.mdep}, $\hat{\sigma}_n^2$ is a weakly consistent estimator of $\sigma^2 >0$. Hence, by Slutsky's theorem, the continuous mapping theorem, and (\ref{clt.mdep.unst}), it suffices to show that $\hat{\gamma}_n ^2$, as defined in (\ref{eq.gam.defn}), is a weakly consistent estimator of $\gamma_1^2$. We may assume, without loss of generality, that $\ee \left[ X_1 \right] = 0$. Let $Z_i = X_i X_{i+1}$. By the $m$-dependent central limit theorem, we have that $\bar{X}_n = O_p \left(1/\sqrt{n}\right)$. Hence, for all $i$,

$$
\begin{aligned}
Y_i &= Z_i - \bar{X}_n \left(X_i + X_{i+1} \right) + \bar{X}_n ^2  && = Z_i + O_p\left( \frac{1}{\sqrt{n}} \right)  \\
\bar{Y}_n &= \bar{Z}_n - \bar{X}_n ^2 + o_p \left( \frac{1}{n} \right) &&= \bar{Z}_n + O_p \left( \frac{1}{n} \right) \, \, .
\end{aligned}
$$

\noindent It follows that, for all $j \geq 0$, 

$$
\begin{aligned}
\frac{1}{n} \sum_{i = 1} ^{n-j-1} \left(Y_i - \bar{Y}_n \right) \left(Y_{i+j} -\bar{Y}_n \right) &= \frac{1}{n} \sum_{j=1} ^{b_n}\sum_{i = 1} ^{n-j-1} \left(Z_i - \bar{Z}_n  \right) \left(Z_{i+j} -\bar{Z}_n \right) + O_p \left( \frac{1}{\sqrt{n}}\right) \, \, .
\end{aligned}
$$

\noindent In particular, since $b_n = o(\sqrt{n} )$, we have that 

$$
\hat{T}_n ^2 =\frac{1}{n} \sum_{i=1} ^{n-1} \left(Z_i - \bar{Z}_n  \right)^2 +  \frac{2}{n} \sum_{j=1} ^{b_n} \sum_{i = 1} ^{n-j-1} \left(Z_i - \bar{Z}_n  \right) \left(Z_{i+j} -\bar{Z}_n \right) + o_p(1) \, \, .
$$

\noindent In order to show consistency of $\hat{T}_n ^2 $, we may further assume, without loss of generality, that $\ee[Z_i]=\ee\left[ X_1 X_2\right] = 0$.

Since the sequence $\{Z_i\}$ is $(m+1)$-dependent, stationary, and has finite fourth moments, we may once again apply the $m$-dependent central limit theorem, as above, and conclude that 

$$
\begin{aligned}
\hat{T}_n ^2&=\frac{1}{n} \sum_{i=1} ^{n-1} Z_i^2 + \frac{2}{n}  \sum_{j=1} ^{b_n}\sum_{i = 1} ^{n-j-1}Z_i Z_{i+j} + o_p(1) \, \, \\
&=\frac{1}{n} \sum_{i=1} ^{n-1} X_i^2X_{i+1} ^2 + \frac{2}{n}  \sum_{j=1} ^{b_n}\sum_{i = 1} ^{n-j-1}X_i X_{i+1}X_{i+j}X_{i+j+1} + o_p(1) \, \, .
\end{aligned}
$$

\noindent For each $0 \leq j \leq b_n$, let 

$$
\hat{T}_{n, \, j}  = \frac{1}{n} \sum_{i=1} ^{n-j-1} X_i X_{i+1}X_{i+j}X_{i+j+1} \, \, .
$$

\noindent For each $j \in \nn$, let

$$
t_j = \cov\left(X_1 X_2, \, X_{j+1} X_{j+2}  \right)
$$

\noindent Since $b_n= o\left(\sqrt{n} \right)$, it suffices to show that, for each $j$,
\beq
\hat{T}_{n, \, j}-  t_j = O_p \left( \frac{1}{\sqrt{n}} \right) \, \, .
\label{t.slutsky.piece}
\eeq

\noindent Fix $j$. We split the sum defining $\hat{T}_{n, \, j}$ into $(m+1)$ sums of i.i.d. random variables, i.e. for $1 \leq i \leq m+1$, let $A_i = \{ k\in [n]: k \equiv i \mod (m+1) \} $, and let 

$$
\hat{T}_{n, \, j, \, i} = \frac{1}{n} \sum_{k \in A_i, \, k + j + 1\leq n} X_k X_{k+1} X_{k+j} X_{k+j+1}  \, \, .
$$

\noindent Each term in the sum is i.i.d., with mean $t_j$ and finite variance, so, by Chebyshev's inequality, for each $i \leq m +1 $,

$$
\hat{T}_{n, \, j, \, i} - \frac{1}{m+1} t_j = O_p \left( \frac{1}{\sqrt{n}} \right) \, \, .
$$

\noindent Hence, summing over $i$, (\ref{t.slutsky.piece}) holds, and thus $\hat{T}_n ^2$ is a consistent estimator of $\tau_1 ^2$. Similarly, we have that $\hat{K}_n ^2$ and $\hat{\nu}_n$ are consistent estimators of $\kappa^2$ and $\nu_1$, respectively. By Theorem \ref{thm.unperm.rho.lim}, we have that $\hat{\rho}_n$ is a consistent estimator of $\rho_1$, and so, by several applications of Slutsky's theorem, $\hat{\gamma}_n ^2$ is a consistent estimator of $\gamma_1 ^2$. Thus concludes the proof of part (i). 

We turn our attention to part (ii) of the Theorem. Let $\Pi_n \sim \text{Unif} \{S_n\}$. By Slutsky's theorem for randomization distributions, Theorem \ref{rand.dist.mdep}, and the consistency of $\hat{\sigma}_n ^2$ under permutation, it suffices to show that 

\beq
\hat{\gamma}_n ^2 \left( \Pi_n \mathbf{X} \right) \overset{p} \to \sigma^4 \, \, .
\label{need.gamma.conv.mdep}
\eeq

\noindent By Theorem \ref{rand.dist.mdep} and Hoeffding's condition (see \cite{chung.perm}, Theorem 5.1), we have that 

$$
\hat{\rho}_n \left( \Pi_n \mathbf{X} \right) \overset{p} \to 0 \, \, .
$$

\noindent Therefore, by Slutsky's theorem, it suffices to show the following:

$$
\begin{aligned}
\hat{T}_n ^2 (\Pi_n \mathbf{X}) &\overset{p} \to \sigma^4\\
\hat{K}_n ^2 (\Pi_n \mathbf{X}) &\overset{p} \to \var\left( X_1 ^2 \right) \\
\hat{\nu}_n (\Pi_n \mathbf{X}) &\overset{p} \to 0 \, \, .
\end{aligned}
$$

\noindent We begin by showing that $\hat{T}_n ^2 (\Pi_n \mathbf{X}) \overset{p} \to \sigma^4$. Once again, we may assume that $\ee X_1 = 0$. Let $\zeta_i = X_{\Pi_n(i)} X_{\Pi_n (i+1)}$. Noting that the sample mean is invariant under permutation, we may apply the same argument as in the proof of part i), and so

$$
\hat{T}_n ^2 (\Pi_n \mathbf{X}) = \frac{1}{n} \sum_{i=1}^{n-1} \left(  \zeta_i -\bar{\zeta}_n \right)^2 +
\frac{2}{n} \sum_{j=1}^{b_n}\sum_{i=1}^{n-j-1}\left(  \zeta_i -\bar{\zeta}_n \right)\left(  \zeta_{i+j} -\bar{\zeta}_n \right) + o_p (1) \, \, .
$$

\noindent Let 

$$
A_n = \{ i \in [n-1]: \left| \Pi_n (i) - \Pi_n(i+1) \right| \leq m \} \, \, .
$$

\noindent Note that this $A_n$ is the same as the one defined in (\ref{set.bounds}). We have that

\beq
\hat{T}_{n, \, 0} (\Pi_n \mathbf{X}) \overset{p} \to \sigma^4 \, \, ,
\label{rand.t0}
\eeq

\noindent and, for $1 \leq j \leq b_n$,

\beq
\hat{T}_{n, \, j} (\Pi_n \mathbf{X}) = o_p \left( \frac{1}{\sqrt{n}} \right) \, \, .
\label{rand.tj}
\eeq

\noindent We begin by proving (\ref{rand.t0}). Condition on $\Pi_n = \pi$. Unless $i \in A_n$, we have that $\ee X_{\pi(i)} ^2 X_{\pi (i+1)} ^2 = \sigma^4$, and, if $i \in A_n$, $\ee X_{\pi(i)} ^2 X_{\pi (i+1)} ^2  \leq K$, for some finite $K$ (as defined in the proof of Lemma \ref{clt.mdep}). Hence

\beq
\ee \left[ \left| \frac{1}{n} \sum_{i=1} ^{n-1} X_{\Pi_n(i)} ^2 X_{\Pi_n(i+1)}^2 - \sigma^4 \right| \, | \, \Pi_n \right] \leq \frac{1 + \left(K+\sigma^4\right) \left| A_n \right|}{n} \, \, .
\label{bound.1.stud}
\eeq

\noindent Let $\epsilon > 0$. By Markov's inequality, 

\beq
\begin{aligned} 
\pp \left(\left| \frac{1}{n} \sum_{i=1} ^{n-1} X_{\Pi_n(i)} ^2 X_{\Pi_n(i+1)}^2 - 1 \right| > \epsilon \right) &\leq \frac{\ee\left| \frac{1}{n} \sum_{i=1} ^{n-1} X_{\Pi_n(i)} ^2 X_{\Pi_n(i+1)}^2 - 1 \right| }{ \epsilon} \\
&\leq \frac{1 + (K+\sigma^4) \ee\left| A_n \right|}{n \epsilon} \, \, \, (\text{by } (\ref{bound.1.stud})) \\
&\leq \frac{1 + 2m(K+\sigma^4)}{n\epsilon} \, \, \, (\text{by } (\ref{an.bound})) \\
&\to 0 \, \, \, \text{as } n \to \infty \, \, .
\end{aligned}
\label{eq.allconv.mdep}
\eeq

\noindent Hence (\ref{rand.t0}) holds. It remains to prove (\ref{rand.tj}). Let $\Gamma_n$ be as defined in (\ref{set.bounds}). Fix $j \in [m+1]$. Let 

$$
B_n ' = \{ \{i, \, j  \} : \{\Pi_n (i), \, \Pi_n (i+1), \, \Pi_n (i+j), \, \Pi_n (i+j +1) \} \in \Gamma_n \} \, \, .
$$

\noindent Condition once more on $\Pi_n = \pi$. We observe that, for $B_n$ as defined in (\ref{set.bounds}), we have that $\ee \left| B_n ' \right| \leq \ee \left| B_n \right|$. Also note that, unless $\{i, \,  j \} \in B_n '$, we have that 

$$
\ee X_{\Pi_n (i)} X_{\Pi_n (i+1) } X_{\Pi_n (i + j)} X_{\Pi_n (i + j + 1)} = 0\, \, ,
$$

\noindent and, if $\{i, \, j \} \in B_n '$, the expectation is, once again, bounded uniformly by $K$. Therefore, unconditionally,

\beq
\ee \left| \sum_{j=1} ^{n-2}\frac{1}{n}  \sum_{i=1} ^{n-j -1}  X_{\Pi_n (i)} X_{\Pi_n (i+1) } X_{\Pi_n (i + j)} X_{\Pi_n (i + j + 1)} \right| \leq \frac{K \ee \left| B_n ' \right| }{n} \, \, ,
\label{bound.rem.stud}
\eeq

\noindent and, in particular, 

\beq
\ee \left|\frac{1}{n} \sum_{i=1} ^{n-j -1}  X_{\Pi_n (i)} X_{\Pi_n (i+1) } X_{\Pi_n (i + j)} X_{\Pi_n (i + j + 1)} \right| \leq \frac{K \ee\left| B_n ' \right| }{n} \, \, .
\label{bound.2.stud}
\eeq

\noindent Let $\epsilon > 0$. By Markov's inequality, 

$$
\begin{aligned} 
&\pp \left(\left|\frac{1}{n} \sum_{i=1} ^{n-j -1}  X_{\Pi_n (i)} X_{\Pi_n (i+1) } X_{\Pi_n (i + j)} X_{\Pi_n (i + j + 1)} \right| > \epsilon \right) \\ 
&\leq \frac{\ee\left|\frac{1}{n} \sum_{i=1} ^{n-j -1}  X_{\Pi_n (i)} X_{\Pi_n (i+1) } X_{\Pi_n (i + j)} X_{\Pi_n (i + j + 1)} \right| }{ \epsilon} \\
&\leq \frac{K \ee\left| B_n ' \right|}{n \epsilon} \, \, \, (\text{by } (\ref{bound.2.stud})) \\
&\leq \frac{24Km^2}{n\epsilon} \, \, \, (\text{by } (\ref{bn.bound})) \, \, .
\end{aligned}
$$

\noindent Hence the first convergence in (\ref{eq.allconv.mdep}) holds. By an identical argument, the remaining two convergences in (\ref{eq.allconv.mdep}) also hold, and so the result follows. \end{proof}

\section{Permutation distribution for $\alpha$-mixing sequences}

\begin{lem} Let $\{X_n, \, n \in \zz\}$ be a stationary, $\alpha$-mixing sequence. If the sequence $\{X_n\}$ is exchangeable, then the $X_n$ are independent and identically distributed.

\label{lem.exch.indep}
\end{lem}

\begin{proof} It suffices to show that, for all $n \in \nn$, for events $A_i \in \sigma(X_i)$, 

$$
\pp\left( \bigcap_{i=1} ^n A_i \right) = \prod_{i=1} ^n \pp(A_i) \, \, .
$$

\noindent We show this by induction on $n$. As a base case, take $n = 2$. By the exchangeability condition, we have that, for all $N > 2$, 

$$
\begin{aligned}
&\{Y_n, \, n \in \zz \}  
\\&:=\left\{\dots,X_{-2}, \, X_{-1} , \, X_2, \, X_0, \, X_N, \, X_{N-1}, \, \dots,\, X_1, \, X_{N+1}, \, X_{N+2}, \, \dots\right\}\overset{d} = \{X_n, \, n \in \zz \} \, \, .
\end{aligned}
$$

\noindent In particular, considering the $\alpha$-mixing coefficients of $Y$, we have that 

$$
\begin{aligned}
\left | \pp\left(A_1 \cap A_2 \right) - \pp(A_1) \pp(A_2) \right| &\leq \alpha_Y ( N+1) \\
&= \alpha_X ( N+1) \\
&\to 0 \, \, , \, \, \, \text{as $N \to \infty$ }.
\end{aligned}
$$

\noindent Hence the result holds for $n =2$. We now show the inductive step, assuming the inductive hypothesis holds for all $m \leq n$. For all $1 \leq i \leq n+1$, let $A_i \in \sigma(X_i)$. Let 

$$
B = \cap_{i=1}^n A_i \, \, .
$$

\noindent By exchangeability, we have that, for all $N > n + 1$, 

$$
\begin{aligned}
&\{Y_n, \, n \in \zz\}
\\&:=\left\{\dots,X_{-2}, \, X_{-1} , \, X_{n+1}, \, X_0, \, X_N, \, X_{N-1}, \, \dots,\, X_{n+2}, \, X_{1}, \, X_{2}, \, \dots, X_n, \, X_{N+1}, \, X_{N+2}, \, \dots\right\}\\ &\overset{d} = \{X_n, \, n \in \zz\}\, \, .
\end{aligned}
$$

\noindent In particular, it follows that 

$$
\begin{aligned}
\left | \pp\left(A_{n+1} \cap B \right) - \pp(A_{n+1}) \pp(B) \right| &\leq \alpha_Y ( N-n + 1) \\
&= \alpha_X ( N-n +1) \\
&\to 0 \, \, , \, \, \, \text{as $N \to \infty$ }.
\end{aligned}
$$

\noindent Hence, by the inductive hypothesis,

$$
\begin{aligned}
\pp\left( \bigcap_{i=1} ^{n+1} A_i \right) &= \pp( A_{n+1} ) \pp(B) \\
&=\prod_{i=1} ^{n+1} \pp(A_i) \, \, ,
\end{aligned}
$$

\noindent as claimed.
\end{proof}

\begin{rem} \rm Since $m$-dependent sequences are $\alpha$-mixing, the result of Lemma \ref{lem.exch.indep} also holds for $\{X_n\}$ an $m$-dependent sequence. \qed

\end{rem}

The following Lemma, and its proof, will be used extensively to show further results throughout this Section.

\begin{lem} Let $\{X_n\}$ be a stationary, $\alpha$-mixing sequence, with mean $0$ and variance 1, such that 

$$
\lVert X_n \rVert _{8 + 4\delta} < \infty \, \, ,
$$

\noindent for some $\delta > 0$. Suppose the $\alpha$-mixing coefficients of $\{X_n\}$ satisfy 

$$
\sum_{n \geq 1} \alpha_X (n) ^{\frac{\delta}{2 + \delta}} < \infty \, \, .
$$

\noindent Let $\Pi_n \sim \text{Unif}(S_n)$, independent of $\{X_n\}$. We have that, as $n \to \infty$, 

\beq
 \ee \left[ \frac{1}{\sqrt{n}}\sum_{i =1} ^{n-1} X_{\Pi_n (i)} X_{\Pi_n (i +1)} \right] \to 0 \, \, ,
\label{alpha.mean}
\eeq

\noindent and 

\beq
\var\left( \frac{1}{\sqrt{n}} \sum_{i =1} ^{n-1} X_{\Pi_n (i)} X_{\Pi_n (i +1)} \right) \to 1 \, \, .
\label{alpha.var}
\eeq

\label{mean.var.alpha.conv}
\end{lem}

\begin{rem} 
For $\{X_n\}$ a stationary sequence, and for fixed $n \in \nn$, the sequence $\{X_{\Pi_n(i)} X_{\Pi_n (i +1)}, \, i \in [n-1] \}$ is also stationary.\qed
\label{rem.stat.perm}
\end{rem}

\begin{proof} We begin by showing (\ref{alpha.mean}). By Remark \ref{rem.stat.perm}, it suffices to show that, as $n \to \infty$, 

$$
\sqrt{n} \ee \left[ X_{\Pi_n (1)} X_{\Pi_n (2)} \right] \to 0 \, \, .
$$

\noindent For ease of notation, let $\alpha (\cdot)$ denote $\alpha_X (\cdot)$.

$$
\begin{aligned}
\left| \ee \left[ X_{\Pi_n (1)} X_{\Pi_n (2)} \right] \right| &= \frac{2}{n(n-1)} \left |\sum_{i < j} \ee \left[ X_ i X_j \right] \right| \\
&\leq \frac{2}{n(n-1)}\sum_{i < j}  \left | \ee[ X_i X_j ] \right| \\
&=  \frac{2}{n(n-1)}\sum_{i < j}  \left | \cov(X_i, \, X_j )\right| \, \, .
\end{aligned}
$$

\noindent By \cite{doukhan}, Section 1.2.2, Theorem 3, we have that, for $1 \leq i < j \leq n$, 

$$ 
\begin{aligned}
\left| \cov(X_i, \, X_j)  \right| &\leq 8\alpha (j - i) ^{\frac{\delta}{2 + \delta}} \Vert X_i \rVert_{2 + \delta} \lVert X_j \rVert_{2 + \delta} \\
&\leq C_1 \alpha (j - i) ^{\frac{\delta}{2 + \delta}} \, \, ,
\end{aligned}
$$

\noindent for some constant $C_1 \in \rr_+$. Hence, for some constant $C_2$, not depending on $n$,  

\beq
\begin{aligned}
\left| \ee \left[ X_{\Pi_n (1)} X_{\Pi_n (2)} \right] \right| &\leq \frac{C_2}{n^2} \sum_{i < j} \alpha (j - i) ^{\frac{\delta}{2 + \delta}} \\
&= \frac{C_2}{n} \sum_{i=1} ^{n-1} \frac{n - i}{n} \alpha (i) ^{\frac{\delta}{2 + \delta}} \\
&\leq \frac{C_2}{n} \sum_{n \geq 1} \alpha (n)^{\frac{\delta}{2 + \delta}} \\
&= O\left(\frac{1}{n} \right) \, \, ,
\end{aligned}
\label{smallo.1.alpha.mean}
\eeq

\noindent and so (\ref{alpha.mean}) holds. We turn our attention to (\ref{alpha.var}). By (\ref{alpha.mean}), it suffices to show that, as $n \to \infty$, 

$$
\begin{aligned}
\frac{1}{n} \ee\left[ \left(\sum_{i=1} ^{n-1}X_{\Pi_n (i)} X_{\Pi_n (i +1)} \right)^2 \right] \to 1 \, \, .
\end{aligned}
$$

\noindent By repeated application of Remark \ref{rem.stat.perm},

\beq
\begin{aligned}
&\frac{1}{n} \ee\left[ \left(\sum_{i=1} ^{n-1}X_{\Pi_n (i)} X_{\Pi_n (i +1)} \right)^2 \right] \\
&=  \frac{1}{n} \sum_{i=1}^{n-1} \ee\left[ X^2_{\Pi_n (i)} X^2_{\Pi_n (i +1)}\right] + \frac{1}{n} \sum_{i=1} ^{n-2} \ee \left[ X_{\Pi_n (i)} X^2_{\Pi_n (i +1) X_{\Pi_n (i+2)}} \right] + \\
&+ \frac{1}{n} \sum_{ \left| i -j \right| > 1}  \ee \left[ X_{\Pi_n (i)} X_{\Pi_n (i +1)} X_{\Pi_n (j)}X_{\Pi_n (j+1)} \right] \\
&= \frac{n-1}{n}\ee\left[ X^2_{\Pi_n (1)} X^2_{\Pi_n (2)}\right] + 
\frac{n-2}{n}\ee \left[ X_{\Pi_n (1)} X^2_{\Pi_n (2) }X_{\Pi_n (3)} \right] + \\
&+\frac{(n-2)^2}{n}\ee \left[ X_{\Pi_n (1)} X_{\Pi_n (2)} X_{\Pi_n (3)}X_{\Pi_n (4)} \right] \, \, .
\end{aligned}
\label{var.alpha.terms}
\eeq

\noindent We now find the limit of each term on the right hand side of (\ref{var.alpha.terms}). We begin by considering the first term. 

$$
\begin{aligned}
\ee\left[ X^2_{\Pi_n (1)} X^2_{\Pi_n (2)}\right] &= \frac{2}{n (n-1)} \sum_{i < j} \left(\ee\left[X^2_{i}\right] \ee\left[X^2_{j} \right] + \cov\left(X^2_{i}, \,X^2_{j}\right) \right)\\
&= 1 +  \frac{2}{n (n-1)} \sum_{i < j} \cov\left(X^2_{i}, \,X^2_{j}\right) \, \, .
\end{aligned}
$$

\noindent By \cite{doukhan}, Section 1.2.2, Theorem 3, the bounded moment condition on the $\{X_n\}$, and the Cauchy-Schwarz inequality, there exists a constant $C \in \rr_+$, independent of $n$, such that 

$$
\begin{aligned}
\frac{2}{n(n-1)}\left|\sum_{i <j} \cov(X_i ^2, \, X_j ^2)  \right| &\leq \frac{2}{n(n-1)} \sum_{i <j} 8\alpha (j - i) ^{\frac{\delta}{2 + \delta}} \Vert X_i ^2 \rVert_{2+ \delta} \lVert X_j^2 \rVert_{2 + \delta} \\
&\leq  \frac{2C}{n(n-1)} \sum_{i <j}  \alpha (j - i) ^{\frac{\delta}{2 + \delta}} \\
&= o(1) \, \, ,
\end{aligned}
$$

\noindent by the same argument as in (\ref{smallo.1.alpha.mean}). Hence

$$
 \frac{n-1}{n}\ee\left[ X^2_{\Pi_n (1)} X^2_{\Pi_n (2)}\right]  = 1 + o(1) \, \, .
$$

\noindent We now consider the second term in (\ref{var.alpha.terms}).

$$
\begin{aligned}
\left| \ee \left[ X_{\Pi_n (1)} X^2_{\Pi_n (2) }X_{\Pi_n (3)} \right]  \right| & \leq \frac{2}{n(n-1)(n-2)} \sum_{i < k, \, j \neq i, \, k} \left| \ee \left[ X_i X_j ^2 X_k \right] \right| \, \, .
\end{aligned}
$$

\noindent Fix $i$ and $k$. Suppose $j < i$ or $j >k$; without loss of generality, suppose $j < i$. By \cite{doukhan}, Section 1.2.2, Theorem 3, 

$$
\begin{aligned}
\left| \ee \left[ X_i X_j ^2 X_k \right] \right| &= \left| \cov\left( X_i X_j^2, \, X_k \right) \right|  \\
&\leq \lVert X_k \rVert_{2 + \delta} \lVert X_i X_j ^2 \rVert_{2 + \delta} \cdot 8\alpha(k - i)^{\frac{\delta}{2 + \delta}} \\
&\leq8 \lVert X_k \rVert_{2 + \delta} \lVert X_i \rVert_{4 + 2\delta} ^2 \lVert X_j \rVert_{8 + 4\delta} ^4 \alpha(k -i)^{\frac{\delta}{2 + \delta}} \, \, \, \text{(Cauchy-Schwarz)} \\
&\leq C \alpha(k-i) ^{\frac{\delta}{2 + \delta}}\, \, ,
\end{aligned}
$$

\noindent for some universal constant $C$, by the moment condition on the $(X_n)$. Suppose $i < j < k$. Without loss of generality, suppose $j - i \leq k-j \iff j \leq \lfloor \frac{1}{2} (k-i) \rfloor$. By a similar argument,

$$
\begin{aligned}
\left| \ee \left[ X_i X_j ^2 X_k \right] \right| &= \left| \cov\left( X_i X_j^2, \, X_k \right) \right| \\
&\leq \lVert X_k \rVert_{2 + \delta} \lVert X_i \rVert_{4 + \delta} ^2 \lVert X_j \rVert_{8 + 4\delta} ^4 \alpha(k -j) ^{\frac{\delta}{2 + \delta}} \\
&\leq C \alpha(k -j ) ^{\frac{\delta}{2 + \delta}} \, \, ,
\end{aligned}
$$

\noindent by the moment conditions on the $\{X_n\}$. It follows that, for some $C, \, K \in \rr_+$, independent of $n$, 

\beq
\begin{aligned}
\left| \ee \left[ X_{\Pi_n (1)} X^2_{\Pi_n (2) }X_{\Pi_n (3)} \right]  \right| &\leq \frac{C}{n^3}\sum_{i < k} \left[( n- k +i -2) \alpha(k -i)^{\frac{\delta}{2 + \delta}} + 2\sum_{j = i+1} ^{i + \lfloor \frac{1}{2} (k-i) \rfloor} \alpha( k -j)^{\frac{\delta}{2 + \delta}} \right] \\
&\leq \frac{K}{n^3} \sum_{r =1 }^{n-1} \left[ (n-r)^2 + (n-r)(n-r - 1)\right] \alpha(r)^{\frac{\delta}{2 + \delta}} \\
&\leq \frac{K}{n}\sum_{r=1} ^{n-1} \alpha(r)^{\frac{\delta}{2 + \delta}} \\
&= o(1) \, \, .
\end{aligned}
\label{j1.var.term}
\eeq

\noindent Hence the second term in (\ref{var.alpha.terms}) is $o(1)$. We conclude by examining the third term in (\ref{var.alpha.terms}). 

$$
\begin{aligned}
\left|\ee \left[ X_{\Pi_n (1)} X_{\Pi_n (2)} X_{\Pi_n (3)}X_{\Pi_n (4)} \right] \right| &= \frac{24}{n(n-1)(n-2)(n-3)} \left| \sum_{i < j < k < l} \ee \left[ X_i X_j X_k X_l \right] \right| \\
&\leq \frac{24}{n(n-1)(n-2)(n-3)}  \sum_{i < j < k < l} \left| \ee \left[ X_i X_j X_k X_l \right] \right|  \, \, .
\end{aligned}
$$

\noindent Fix $l - i = r \geq 3 $. There are $(n - r)$ such choices for $(i, \, l)$. Without loss of generality, by stationarity, we may subsequently assume $i = 1$, $l = r + 1$. Note that, by the same argument as for the previous two terms, for $1 < j < k < r +1 $, for some constant $C \in \rr_+$ independent of $n$, 

$$
\begin{aligned}
\left| \ee \left[ X_1 X_j X_k X_{r+1} \right] \right| &\leq C \alpha\left( \max \left\{ j-1, \, r+1 - k \right\} \right)^{\frac{\delta}{2 + \delta}} \, \, .
\end{aligned}
$$

\noindent Fix $\max \left\{ j-1, \, r+1 - k \right\} = s$. There are at most $2s$ choices for $(j, \, k)$. Therefore we have that, for some constants $K, \, \tilde{K} \in \rr_+$ independent of $n$, 

\beq
\begin{aligned}
\left|\ee \left[ X_{\Pi_n (1)} X_{\Pi_n (2)} X_{\Pi_n (3)}X_{\Pi_n (4)} \right] \right| &\leq \frac{K}{n^4} \sum_{r = 3} ^{n-1} (n -r) \sum_{s=1} ^{r-2}  2s \cdot \alpha(s)^{\frac{\delta}{2 + \delta}} \\
&\leq \frac{K}{n^4} \sum_{s= 1}^{n-1} 2s(n - s)^2 \alpha(s)^{\frac{\delta}{2 + \delta}} \\
&\leq \frac{\tilde{K}}{n} \sum_{s=1}^{n-1} \left(\frac{s}{n} \right) \left(\frac{n-s}{n}\right)^2 \alpha(s)^{\frac{\delta}{2 + \delta}} \\
&= O \left( \frac{1}{n} \right) \, \, ,
\end{aligned}
\label{final.sum.alpha.var}
\eeq

\noindent by the summability condition on the $\alpha$-mixing coefficients.
\end{proof}

The following Theorem is a slightly more restrictive version of Theorem \ref{alpha.bd.perm}, and its restrictions will be relaxed later.

\begin{thm}

Let $\{X_n\}$ be a stationary, bounded, $\alpha$-mixing sequence, with mean $0$ and variance 1. Suppose that 

$$
\sum_{n \geq 1} \alpha_X (n) < \infty \, \,.
$$

\noindent The permutation distribution of $\sqrt{n} \hat{\rho}_n$ based on the test statistic $\hat{\rho}_n = \hat{\rho}(X_1, \, \dots, \, X_n)$, with associated group of transformations $S_n$, the symmetric group of order $n$, satisfies

$$
\sup_{t \in \rr} \left| \hat{R}_n (t)-\Phi(t)\right| \overset{p}{\to} 0 \, \, ,
$$

\noindent as n $\to \infty$, where $\Phi(t)$ is the distribution of a standard Gaussian random variable.

\label{alpha.bd.perm.known.mean}
\end{thm}

\begin{proof} 
Suppose the sequence $\{X_n\}$ is defined on the probability space $\left( \Omega, \, \mathcal{T}, \, \pp \right)$. Let $\mathcal{F} = \sigma( X_i: i \in \nn)$. Let 

\beq
B(m) := \frac{1}{m}\sum_{i=1} ^m X_i ^2 - \frac{1}{m^2} \left(\sum_{j=1} ^m X_j \right)^2 \, \, .
\label{b.1d.bd.perm}
\eeq

\noindent Note that $B(m) \geq 0$ for all $m \in \nn$. We have that

\beq
\begin{aligned}
\ee\left[\frac{1}{m} \sum_{j=1} ^m X_j \right] &= 0 \, \, ,\\
\ee\left[\frac{1}{m} \sum_{i=1} ^m X_i ^2\right] &= 1 \, \, .
\end{aligned}
\label{bm.means}
\eeq

\noindent Also, by \cite{doukhan}, Section 1.2.2, Theorem 3, by the boundedness of the sequence $\{X_n\}$, there exists a constant $C$, independent of $m$, such that 

$$
\begin{aligned}
0 \leq \frac{1}{m^2} \ee \left[ \left(\sum_{j=1} ^m X_j  \right)^2 \right] &= \frac{1}{m^2} \sum_{j=1} ^m \ee X_j^2 + \frac{2}{m^2} \sum_{i < j} \ee X_i X_j \\
&\leq \frac{1}{m} + \frac{2}{m^2} \sum_{i < j} \lvert \cov\left(X_i, \, X_j \right)\rvert \\
&\leq \frac{1}{m} + \frac{C}{m^2} \sum_{i < j} \alpha_X (j-i)\\
&=\frac{1}{m} + \frac{C}{m^2} \sum_{j=1}^{m-1} (m-j) \alpha_X (j) \\
&= o(1) \, \, ,
\end{aligned}
$$

\noindent since $\alpha_X (j) \to 0$ as $j \to \infty$. By Chebyshev's inequality, it follows that 

\beq
\frac{1}{m} \sum_{j=1} ^m X_j \overset{p} \to 0 \, \, .
\label{bm.xbar.0}
\eeq

\noindent Similarly, there exists a constant $K$, independent of $m$, such that 

$$
\begin{aligned}
0 \leq \var \left(\frac{1}{m} \sum_{i=1} ^m X_i ^2\right) &= \frac{1}{m}\var\left(X_1 ^2 \right) + 2\sum_{i < j} \cov\left( X_i^2, \, X_j^2 \right) \\
&\leq \frac{1}{m} \var\left(X_1 ^2 \right) + \frac{2}{m^2}\sum_{i < j} \lvert \cov\left( X_i^2, \, X_j^2 \right)\rvert \\
&\leq \frac{1}{m}\var\left(X_1 ^2 \right) + \frac{K}{m^2} \sum_{i < j} \alpha_X (j-i)\\
&= o(1) \, \, .
\end{aligned}
$$

\noindent Hence, by Chebyshev's inequality,

\beq
\frac{1}{m} \sum_{i=1} ^m X^2_i \overset{p} \to 1 \, \, .
\label{bm.x2bar.1}
\eeq

\noindent By (\ref{bm.xbar.0}) and (\ref{bm.x2bar.1}), the continuous mapping theorem, and Slutsky's theorem, we have that, as $m \to \infty$, 

$$
B(m) \overset{p} \to 1 \, \, .
$$

\noindent Let $\left\{n_j, \, j \in \nn\right\}\subset \nn$. By an alternative characterization of convergence in probability (see Theorem 20.5 in \cite{billingsley}), we have that there exists a subsequence $\left\{n_{j_k} \right\} \subset \left\{n_j\right\}$ such that  

$$
B\left(n_{j_k} \right)\overset{a.s.} \to 1 \, \, .
$$

\noindent Let 

\beq
A = \left\{ \omega: \lim_{k \to \infty} \inf_{l \geq k} B(n_{j_l}) >0 \right\} \, \, .
\label{set.a.perm.bd}
\eeq

\noindent We have that $\pp (A) = 1$. Also, the sequence $\{X_n\}$ is bounded, and, for all $\omega \in A$, 

$$
\lim_{k \to \infty} \inf_{l \geq k} \left\{\frac{1}{n_{j_l}}\sum_{i=1} ^{n_{j_l}}X_i (\omega)^2 - \frac{1}{n^2_{j_l}} \left(\sum_{i=1} ^{n_{j_l}}X_i (\omega) \right)^2 \right\} > 0 \, \, . 
$$

\noindent Let 

$$
Z_{\Pi_n} (\omega) := \sum_{i=1} ^{n-1} X_{\Pi_n (i)} (\omega)X_{\Pi_n (i+1)}(\omega) +X_{\Pi_n (n)}(\omega) X_{\Pi_n (1)} (\omega) \,\, .
$$

\noindent We may now apply the theorem of \cite{wald1943} (see page 383)\footnote{We do not, strictly speaking, apply the theorem itself. We repeat the proof of the theorem on the subsequence $\left(n_{j_k}\right)$.} to conclude that, for each $\omega \in A$, 

$$
\frac{Z_{\Pi_{n_{j_k} }} (\omega) - \ee_{\Pi_{n_{j_k} }} \left[ Z_{\Pi_{n_{j_k} }} (\omega)\right] }{\var_{\Pi_{n_{j_k} }}\left(Z_{\Pi_{n_{j_k} }}(\omega)\right)^{1/2} } \overset{d} \to N(0, \, 1) \, \, .
$$



\noindent Let $t \in \rr$. We have that, for all $\omega \in A$, 

$$
\bar{R}_{n_{j_k}} (t) (\omega)  \to \Phi(t) \, \, ,
$$

\noindent where 

\beq
\bar{R}_n (t) = \frac{1}{n!} \sum_{\pi_n \in S_n} 1\left\{ Z_{\pi_n} \leq t\var\left(Z_{\Pi_n}  \, | \, \mathcal{F}\right)^{1/2} + \ee \left[ Z_{\Pi_n}  \, | \, \mathcal{F}\right] \right\} \, \, .
\label{rbar.perm}
\eeq 

\noindent In particular, it follows that, since $\pp(A) =1$, 

\beq
\bar{R}_{n_{j_k}} (t) \overset{a.s.} \to\Phi(t) \, \, .
\label{perm.bd.subseq.as}
\eeq

\noindent In conclusion, for every sequence $\left\{n_j \right\} \subset \nn$, there exists a subsequence $\left\{n_{j_k} \right\}$ such that (\ref{perm.bd.subseq.as}) holds. It follows that, for all $t \in \rr$, as $n \to \infty$,

$$
\bar{R}_{n} (t) \overset{p} \to\Phi(t) \, \, .
$$

\noindent Let 

\beq
\begin{aligned}
M_n &=\frac{1}{\sqrt{n}}\ee_{\Pi_n} \left[ \sum_{i=1} ^{n-1}X_{\Pi_n (i)}X_{\Pi_n (i+1)} + \frac{1}{\sqrt{n}} X_{\Pi_n (n)}X_{\Pi_n (1)}\, | \, \mathcal{F}\right] - \frac{1}{\sqrt{n}} X_{\Pi_n (n)} X_{\Pi_n (1)}\, \, , \\
S_n &= \frac{1}{\sqrt{n}} \var_{\Pi_n} \left( \sum_{i=1} ^{n-1}X_{\Pi_n (i)}X_{\Pi_n (i+1)} + \frac{1}{\sqrt{n}} X_{\Pi_n (n)} X_{\Pi_n (1)} \, | \, \mathcal{F}\right)^{1/2} \, \, . \\
\end{aligned}
\label{mn.sn.def}
\eeq

\noindent Note that the expectations and variances in (\ref{mn.sn.def}) are conditional on the data, and the last term in the definition of $M_n$ accounts for the difference between the circular and non-circular definitions of autocovariance. Let

$$
Y_n = \frac{1}{\sqrt{n}} \sum_{i=1} ^{n-1} X_{\Pi_n (i)} X_{\Pi_n (i+1)}  \, \, .
$$

\noindent Note that

$$
\frac{Z_{\Pi_{n } }(\omega) - \ee_{\Pi_n } \left[ Z_{\Pi_n} (\omega)\right] }{\var_{\Pi_n}\left(Z_{\Pi_n}(\omega)\right)^{1/2} } = \frac{Y_n (\omega)- M_n(\omega)}{S_n(\omega)} \, \, .
$$

\noindent In order to show that 

\beq
\hat{R}_n (t) \overset{p} \to \Phi(t) \, \, ,
\label{tilde.perm.bd}
\eeq

\noindent by dividing both sides of the inequality in the indicator function in (\ref{rbar.perm}) by $\sqrt{n}$, and by Slutsky's theorem for randomization distributions (see Theorem 5.1 and Theorem 5.2 of \cite{chung.perm}), it suffices to show that 

\beq
M_n \overset{p} \to 0 \, \,,
\label{mn.bd}
\eeq

\noindent and

\beq
S_n \overset{p} \to 1 \, \, .
\label{sn.bd}
\eeq

\noindent We begin by showing (\ref{mn.bd}). Note that, since the ${X}_n$ are bounded,

\beq
\frac{1}{\sqrt{n}} X_{\Pi_n (n)} X_{\Pi_n (1)} = o_p (1) \, \, .
\label{mn.bd.inter.1}
\eeq

\noindent Also,

\beq
\begin{aligned}
&T_n:= \frac{1}{\sqrt{n}}\ee_{\Pi_n} \left[ \sum_{i=1} ^{n-1}{X}_{\Pi_n (i)}{X}_{\Pi_n (i+1)} + \frac{1}{\sqrt{n}} {X}_{\Pi_n (n)} {X}_{\Pi_n (1)}\, | \, \mathcal{F}\right]\\
&=\sqrt{n} \ee_{\Pi_n} \left[X_{\Pi_n (1)} X_{\Pi_n(2)} \, | \, \mathcal{F} \right] \\
&= \sqrt{n} \cdot \frac{2}{n(n-1)} \sum_{i <j} X_i X_j \, \, .
\end{aligned}
\label{mn.bd.inter.2}
\eeq

\noindent Applying the result of \cite{doukhan}, Section 1.2.2, Theorem 3, and noting that the $X_n$ are bounded, we observe that there exists a constant $C$, independent of $n$, such that 

\beq
\begin{aligned}
\lvert \ee T_n \rvert &\leq \frac{2}{n(n-1)} \sum_{i < j} \lvert \ee \left[ X_i X_j \right] \rvert \\
&\leq \frac{2C}{n(n-1)} \sum_{i < j} \alpha_X (j-i) \\
&= \frac{2C}{n(n-1)} \sum_{i=1}^{n-1} (n - i) \alpha_X (i) \\
&\leq \frac{2C}{n} \sum_{i \geq 1} \alpha_X (i) \\
&= O\left( \frac{1}{n} \right) \, \, ,
\end{aligned}
\label{mn.bd.mean.0}
\eeq

\noindent and so $\ee T_n \to 0$ as $n \to \infty$. We now compute the variance of $T_n$. Note that 

$$
\begin{aligned}
T_n^2 &= \frac{1}{n(n-1)^2} \sum_{i\neq j, \, k \neq l} X_i X_j X_k X_l \\
&=  \frac{1}{n(n-1)^2} \left[ 2\sum_{i \neq j} X_i ^2 X_j ^2 + 4\sum_{i \neq j ,\, k\neq i, \, j} X_i ^2 X_j X_k + 24\sum_{i < j < k < l} X_i X_k X_k X_l \right] \, \,.
\end{aligned}
$$

\noindent By a similar argument to the one showing convergence of (\ref{var.alpha.terms}) in the proof of Lemma \ref{mean.var.alpha.conv}, we observe that 

\beq
\ee T_n ^2 \to 0 \, \, .
\label{mn.bd.var.0}
\eeq

\noindent By (\ref{mn.bd.mean.0}) and (\ref{mn.bd.var.0}), and a straightforward application of Chebyshev's inequality, we have that 

\beq
T_n \overset{p} \to 0 \, \, .
\label{mn.bd.inter.3}
\eeq

\noindent By Slutsky's theorem, and combining (\ref{mn.bd.inter.1}) and (\ref{mn.bd.inter.3}), we see that (\ref{mn.bd}) holds. We turn our attention to (\ref{sn.bd}). Since $M_n \overset{p}\to 0$, it suffices to show that 

\beq
\frac{1}{n}\ee_{\Pi_n} \left[ \left(\sum_{i=1} ^{n-1} X_{\Pi_n (i)} X_{\Pi_n (i+1)} + X_{\Pi_n (n)} X_{\Pi_n (1)} \right)^2 \, | \, \mathcal{F}\right] \overset{p} \to 1 \, \, .
\label{sn.bd.inter.1}
\eeq

\noindent By Lemma 6 of \cite{wald1943} and boundedness of the $\{X_n\}$, we have that 

$$
\frac{1}{n}\ee_{\Pi_n} \left[ \left(\sum_{i=1} ^{n-1} X_{\Pi_n (i)} X_{\Pi_n (i+1)} + X_{\Pi_n (n)} X_{\Pi_n (1)} \right)^2 \, | \, \mathcal{F}\right] =  \ee_{\Pi_n} \left[ X^2_{\Pi_n (1)}  X^2_{\Pi_n (2)} \, | \, \mathcal{F} \right] + o(1) \, \, .
$$

$$
\begin{aligned}
\ee_{\Pi_n} \left[ X_{\Pi_n (1)}  X_{\Pi_n (2)} \, | \, \mathcal{F} \right] &= \frac{2}{n(n-1)} \sum_{i < j} X_i ^2 X_j ^2 \\
&= 1 + \frac{2}{n(n-1)}\sum_{i <j} \left(X_i ^2 X_j ^2 - 1 \right) \, \, .
\end{aligned}
$$

\noindent By a similar argument to the one used to show (\ref{mn.bd}), we have that 

$$
\frac{2}{n(n-1)}\sum_{i <j} \left(X_i ^2 X_j ^2 - 1 \right) \overset{p} \to 0 \, \,,
$$

\noindent and so (\ref{sn.bd}) holds. Hence we have that 

$$
\hat{R}_n (t) \overset{p} \to \Phi(t) \, \, ,
$$

\noindent as claimed.
\end{proof}

The restrictions of Theorem \ref{alpha.bd.perm.known.mean} may now be relaxed in order to show the result of Theorem \ref{alpha.bd.perm}.

\begin{proof}[{\sc Proof of Theorem \ref{alpha.bd.perm}}] Let $\Pi_n, \, \Pi_n ' \sim \text{Unif}(S_n)$ be independent, and also independent of the sequence $(X_i)$.Let 

$$
\hat{\rho}_{\Pi_n} = \frac{\sum_{i=1}^{n-1} \sum_{i=1} \left(X_{\Pi_n (i)} - \bar{X}_n \right) \left(X_{\Pi_n (i+1)} - \bar{X}_n \right)}{\hat{\sigma}_n ^2} \, \, ,
$$

\noindent and define $\hat{\rho}_{\Pi_n '}$ analogously. Without loss of generality, we may assume that $\ee[X_1] = 0$, and $\var\left( X_1 \right) = 1$. By Hoeffding's condition (see \cite{chung.perm}, Theorem 5.1), it suffices to show that 

$$
\sqrt{n} \begin{pmatrix}
\hat{\rho}_{\Pi_n '}\\\hat{\rho}_{\Pi_n '}
\end{pmatrix} \overset{d}\to N(0, \, I_2) \, \, ,
$$

\noindent where $I_2$ is the $2 \times 2$ identity matrix. Note that  

$$
\hat{\rho}_{\Pi_n} = \frac{1}{\hat{\sigma}_n ^2} \left[\frac{1}{\sqrt{n}} \sum_{i=1}^{n-1}X_{\Pi_n (i)}X_{\Pi_n (i+1)}  - 2\sqrt{n} \bar{X}_n ^2 + o(1) \right] \, \, .
$$

\noindent By assumption, we have that 

$$
\ee \left[ \bar{X}_n \right] = 0 \,\,,
$$

\noindent and, by \cite{doukhan}, Section 1.2.2, Theorem 3, we have that, for some universal constant $K$, independent of $n$, 

$$ 
\begin{aligned}
\var\left(\bar{X}_n \right) &= \frac{1}{n} \var(X_1) + \frac{2}{n^2} \sum_{i < j} \cov\left(X_i, \, X_j \right) \\
&\leq O\left(\frac{1}{n} \right) + \frac{K}{n^2} \sum_{i < j} \alpha_X (j - i) \\
&= O\left(\frac{1}{n} \right) \, \, .
\end{aligned}
$$

\noindent It follows that $\bar{X}_n = O_p(1/\sqrt{n})$, and so 

$$
\hat{\rho}_{\Pi_n} = \frac{1}{\hat{\sigma}_n ^2} \left[\frac{1}{\sqrt{n}} \sum_{i=1}^{n-1}X_{\Pi_n (i)}X_{\Pi_n (i+1)}  +o_p (1) \right] \, \, .
$$

\noindent Obtaining an analogous result for $\hat{\rho}_{\Pi_n'}$, and noting that, by the proof of Theorem \ref{alpha.bd.perm.known.mean},

$$
\hat{\sigma}_n^2 \overset{p} \to 1 \, \, ,
$$

\noindent it follows that  

$$
\sqrt{n} \begin{pmatrix}
\hat{\rho}_{\Pi_n '}\\\hat{\rho}_{\Pi_n '}
\end{pmatrix} = \frac{1}{\sigma_n ^2} \begin{pmatrix}
\frac{1}{\sqrt{n}}  \sum_{i=1}^{n-1}X_{\Pi_n (i)}X_{\Pi_n (i+1)} \\
\frac{1}{\sqrt{n}}  \sum_{i=1}^{n-1}X_{\Pi_n (i)}X_{\Pi_n (i+1)} \
\end{pmatrix} + o_p(1) \, \, .
$$

\noindent By Theorem \ref{alpha.bd.perm.known.mean}, and the converse of \cite{chung.perm}, Theorem 5.1, the result follows.
\end{proof}

The proof of Lemma \ref{trunc.diag.lem} follows, and this will be used in extending the result of Theorem \ref{alpha.bd.perm} to the setting of Theorem \ref{alpha.perm}.

\begin{proof} [{\sc Proof of Lemma \ref{trunc.diag.lem}}]

Enumerate $\mathbb{Q}$ as $\mathbb{Q} = \left\{q_n : n \in \nn \right\}$. For all $N \in \nn$, there exists $\tilde{f}(N) \in \nn$ such that, for all $n \geq \tilde{f}(N)$ and all $j \leq N$, 

$$
\pp \left( \left| G_{N ,\, n} \left( q _j\right)- g_N(q_j) \right| \geq \frac{1}{N} \right) \leq \frac{1}{N} \, \, .
$$

\noindent By letting $f(1) = \tilde{f} (1)$, and recursively defining, for $N \geq 2$,

$$
f(N) = \max \left\{\tilde{f} (N), \, f(N-1) +1 \right\} \, \, ,
$$

\noindent we have that $f$ satisfies the same property as $\tilde{f}$, and is strictly increasing. For all $n \geq f(1)$, let 

$$
N_n = \sup\left\{N: \, f(N)  \leq n \right\} \, \, .
$$

\noindent Note that $N_n \to \infty$ as $n \to \infty$, and that $N_n$ is nondecreasing. Fix $j \in \nn$. Let $r, \, s \in \nn$. We wish to show that, for all $n$ sufficiently large,

$$
\pp \left(\lvert G_{N_n ,\, n}(q_j) - g_{N_n}(q_j) \rvert \geq \frac{1}{r} \right) \leq \frac{1}{s}\, \, .
$$

\noindent It is a fact that, for all $N \geq \max \{ j, \, r, \, s \}$, for all $n \geq f(N)$,

$$
\pp\left(\lvert G_{N, \, n} (q_j)- g_{N_n}(q_j) \rvert \geq \frac{1}{r} \right) \leq \frac{1}{s} \, \, .
$$

\noindent Also, there exists $M \in \nn$ such that, for all $n \geq M$, 

$$
N_n \geq \max\{j, \, r, \, s\}  \, \, .
$$

\noindent Since $f\left(N_n \right) \leq n$, it follows that, for all $n \geq \max\{M, \, f(1) \}$, 

$$
\pp\left(\lvert G_{N_n, \, n} (q_j)- g_{N_n}(q_j) \rvert \geq \frac{1}{r} \right) \leq \frac{1}{s}\, \, .
$$

\noindent Hence we may conclude that, for all $q \in \mathbb{Q}$, as $n \to \infty$, 

$$
G_{N_n, \, n} (q) - g_{N_n}(q)\overset{p}\to 0 \, \, .
$$

\noindent By Slutsky's theorem and convergence of $g_{N_n} (t)$ for all $t \in \rr$, it follows that 

$$
G_{N_n, \, n} (q) \overset{p} \to g(q) \, \, .
$$

\noindent Since $\mathbb{Q}$ is dense in $\rr$, $g$ is continuous, and $G_{N_n, \, n}$ is nondecreasing, it follows that, for all $t \in \rr$, as $n \to \infty$, 

$$
G_{N_n, \, n} (t) \overset{p} \to g(t) \, \, ,
$$

\noindent as desired.
\end{proof}

We are now in a position to prove Theorem \ref{alpha.perm}.

\begin{proof} [{\sc Proof of Theorem \ref{alpha.perm}}]

For each $N \in \nn$, $n \in \nn$, let 

$$
\begin{aligned}
Y^{(N)}_n &= X_n 1\left\{ \lvert X_n \rvert \leq N \right\} \\
Z^{(N)}_n &= X_n 1\left\{ \lvert X_n \rvert > N \right\} \, \, .
\end{aligned}
$$

\noindent We have that 

\beq
\begin{aligned}
&\frac{1}{\sqrt{n}} \sum_{i=1} ^{n-1} X_{\Pi_n (i)} X_{\Pi_n (i+1)} \\
&= \frac{1}{\sqrt{n}} \sum_{i=1} ^{n-1} \left(Y^{(N)}_{\Pi_n (i)}  + Z^{(N)}_{\Pi_n (i)} \right)\left(Y^{(N)}_{\Pi_n (i+1)}  + Z^{(N)}_{\Pi_n (i+1)} \right) \\
&=\frac{1}{\sqrt{n}} \sum_{i=1} ^{n-1} Y^{(N)}_{\Pi_n (i)} Y^{(N)}_{\Pi_n (i+1)}  + 
\frac{1}{\sqrt{n}} \sum_{i=1} ^{n-1} Z^{(N)}_{\Pi_n (i)} Y^{(N)}_{\Pi_n (i+1)} + \\
&+\frac{1}{\sqrt{n}} \sum_{i=1} ^{n-1} Y^{(N)}_{\Pi_n (i)} Z^{(N)}_{\Pi_n (i+1)} +
\frac{1}{\sqrt{n}} \sum_{i=1} ^{n-1} Z^{(N)}_{\Pi_n (i)} Z^{(N)}_{\Pi_n (i+1)} \,\,.
\end{aligned}
\label{perm.4.terms}
\eeq

\noindent Note that all of these terms are random with respect to both the sequence $\{X_n\}$ and the permutation $\Pi_n$. Let 

$$
\begin{aligned}
\mu_N &= \ee\left[ X_1 1\left\{ \lvert X_1 \rvert \leq N \right\} \right] \, \, , \\
\sigma_N ^2 &= \var \left( X_1 1\left\{ \lvert X_1 \rvert \leq N \right\} \right) \, \, ,\\
\nu_N  &= \ee\left[ X_1 1\left\{ \lvert X_1 \rvert > N \right\} \right] \, \, ,\\
\tau_N ^2 &=  \var \left( X_1 1\left\{ \lvert X_1 \rvert > N \right\} \right) \, \, .
\end{aligned}
$$

\noindent For each $N \in \nn$, let $\hat{R}_{N, \, n}$ be the permutation distribution of $\sqrt{n} \hat{\rho}_n$, based on the test statistic 

$$
\hat{\rho}_n = \hat{\rho} _n\left( \frac{Y^{(N)} _1 - \mu_N}{\sigma_N}, \, \dots, \, \frac{Y^{(N)} _n - \mu_N}{\sigma_N} \right) \,\,.
$$

\noindent Note that $\Phi$ is continuous, and, for each $N$, $\hat{R}_{N, \, n}$ is nondecreasing in $t$.  Also, by Theorem \ref{alpha.bd.perm}, we have that, for all $N \in \nn$, $t \in \rr$, as $n \to \infty$, 

\beq
\hat{R}_{N ,\, n} (t)\overset{p} \to \Phi (t) \, \, .
\label{conv.bd.perm}
\eeq

\noindent Hence, by Lemma \ref{trunc.diag.lem}, there exists a sequence $\{N_n, \, n\in \nn\}$, with $N_n \to \infty$, such that 

\beq
\hat{R}_{N_n, \, n} (t) \overset{p} \to \Phi(t) \, \, .
\label{perm.n.inc}
\eeq

\noindent We now show that 

\beq
\frac{\mu_{N_n}}{ \sqrt{n}}\sum_{i=1} ^{n-1} \left(Y^{(N_n)}_{\Pi_n(i)} - \mu_{N_n} \right) +
\frac{\mu_{N_n}}{ \sqrt{n}}\sum_{i=1} ^{n-1} \left(Y^{(N_n)}_{\Pi_n(i+1)} - \mu_{N_n} \right) = o_p(1) \, \, .
\label{2.yn.conv}
\eeq

\noindent Note that, by dominated convergence, as $n \to \infty$, 

\beq
\begin{aligned}
\mu_{N_n} \to 0 \, \, ,\\
\sigma_{N_n} ^2 \to 1 \,\,.
\end{aligned}
\label{mu.sig.yn}
\eeq

\noindent Note also that, since the second moments of the sequence $\{X_n\}$ are finite, we have that, by Chebyshev's inequality,

$$
\begin{aligned}
&\frac{\mu_{N_n}}{\sigma^2_{N_n} \sqrt{n}}\sum_{i=1} ^{n-1} \left(Y^{(N_n)}_{\Pi_n(i)} - \mu_{N_n} \right) +
\frac{\mu_{N_n}}{\sigma^2_{N_n} \sqrt{n}}\sum_{i=1} ^{n-1} \left(Y^{(N_n)}_{\Pi_n(i+1)} - \mu_{N_n} \right)   \\&=\frac{2\mu_{N_n}}{\sigma^2_{N_n} \sqrt{n}}\sum_{i=1} ^{n} \left(Y^{(N_n)}_{i} - \mu_{N_n} \right) +
o_p(1) \, \, .
\end{aligned} 
$$

\noindent Hence, in order for (\ref{2.yn.conv}) to hold, by Slutsky's theorem, it suffices to show that 

\beq
\frac{1}{\sigma^2_{N_n} \sqrt{n}}\sum_{i=1} ^{n} \left(Y^{(N_n)}_{i} - \mu_{N_n} \right) = O_p (1) \, \, .
\label{use.neumann}
\eeq

\noindent For $i \in [n]$, let

$$
\xi_{n, \, i} = \frac{1}{\sigma_{N_n} \sqrt{n}}\left(Y^{(N_n)}_{i} - \mu_{N_n} \right) \, \, .
$$

\noindent For all $n \in \nn$, $i \in [n]$, we have that 

$$
\begin{aligned}
\ee \xi_{n, \, i} &= 0\, \, , \\
\sum_{i=1} ^n \ee \xi_{n, \, i} ^2 &= 1 \, \, .
\end{aligned}
$$

\noindent Also, by dominated convergence, as $n \to \infty$, 

$$
\var \left( \sum_{i=1} ^n \xi_{n, \, i} \right) \to \sigma^2 :=\ee\left[X_1 ^2 \right] + 2\sum_{k \geq 2} \ee \left[X_1 X_k \right] \geq 0\, \, .
$$

\noindent Since the $\xi_{n, \, i}$ satisfy the Lyapunov condition 

$$
\begin{aligned}
\ee\left[ \lvert \xi_{n, \, i} \rvert ^{2 + \delta} \right] &\leq \frac{\left(\lVert X_i \rVert_{2 + \delta} + \lvert \mu_{N_n} \rvert \right)^{2+\delta} }{\sigma_{N_n} ^{2 + \delta} n^{1 + \delta/2}}\\
&= O \left(n ^{-1 - \delta/2} \right) \, \, ,
\end{aligned}
$$

\noindent we have that, for all $\epsilon > 0$, as $n \to \infty$, 

$$
\sum_{i=1} ^n \ee \left[ \xi_{n, \, i } 1\left\{ \xi_{n,\, i} > \epsilon \right\} \right] \to 0 \, \, .
$$

\noindent Since the $\xi_{n, \, i}$ have the same $\alpha$-mixing coefficients as the $X_i$, it follows that conditions (2.3) and (2.4) of \cite{neumann} Theorem 2.1 are satisfied, and so 

$$
\sum_{i=1} ^n \xi_{n, \, i} \overset{d} \to N(0, \, \sigma^2) \, \,.
$$

\noindent In particular, by Slutsky's theorem, it follows that (\ref{use.neumann}) holds, and so (\ref{2.yn.conv}) also holds. 

We now show that 

\beq
\frac{1}{\sqrt{n}} \sum_{i=1} ^{n-1} Z^{(N_n)}_{\Pi_n (i)} Y^{(N_n)}_{\Pi_n (i+1)} + \sqrt{n} \mu_{N_n} ^2= o_p(1) \, \, .
\label{yn.zn.conv}
\eeq

\noindent Since $\mu_{N_n} + \nu_{N_n} = 0$ for all $n \in \nn$, we have that 

$$
\begin{aligned}
&\frac{1}{\sqrt{n}} \sum_{i=1} ^{n-1} Z^{(N_n)}_{\Pi_n (i)} Y^{(N_n)}_{\Pi_n (i+1)} + \sqrt{n} \mu_{N_n} ^2 \\
&= \frac{1}{\sqrt{n}} \sum_{i=1} ^{n-1} \left( Z^{(N_n)}_{\Pi_n (i)} -\nu_{N_n} \right)\left( Y^{(N_n)}_{\Pi_n (i+1)} -\mu_{N_n} \right) +\frac{\nu_{N_n}}{\sqrt{n}} \sum_{i=1} ^{n-1} \left( Y^{(N_n)}_{\Pi_n (i+1)} -\mu_{N_n} \right)  +\\
&+ \frac{\mu_{N_n}}{\sqrt{n}}  \sum_{i=1} ^{n-1}\left( Z^{(N_n)}_{\Pi_n (i)} -\nu_{N_n} \right) \, \, .
\end{aligned}
$$

\noindent By (\ref{mu.sig.yn}), $\nu_{N_n} \to 0$ as $n \to \infty$, so by (\ref{use.neumann}), we have that 

\beq
\frac{\nu_{N_n}}{\sqrt{n}} \sum_{i=1} ^{n-1} \left( Y^{(N_n)}_{\Pi_n (i+1)} -\mu_{N_n} \right) = o_p (1) \, \, .
\label{yn.zn.1}
\eeq

\noindent By a similar argument, it follows that 

\beq
\frac{\mu_{N_n}}{\sqrt{n}}  \sum_{i=1} ^{n-1}\left( Z^{(N_n)}_{\Pi_n (i)} -\nu_{N_n} \right) = o_p(1) \, \, .
\label{yn.zn.2}
\eeq

\noindent We now fully utilize the moment and mixing conditions on the sequence $\{X_n \}$. By an similar proof to that of Lemma \ref{mean.var.alpha.conv}\footnote{The only difference in the proof is that we keep track of a multiplicative factor of $\lVert Z^{(N_n)}_1 - \nu_{N_n}\rVert^k_r$ in front of each term, for some $k > 1, \, r \in [1, \, 8+ 4\delta]$. }, noting that, for all $r \in [1, \, 8 + 4\delta]$, as $n \to \infty$, 

$$
\lVert Z^{(N_n)}_{i} -\nu_{N_n} \rVert_r \to 0 \, \, ,
$$

\noindent uniformly in $r$, we have that 

$$
\ee \left[  \frac{1}{\sqrt{n}} \sum_{i=1} ^{n-1} \left( Z^{(N_n)}_{\Pi_n (i)} -\nu_{N_n} \right)\left( Y^{(N_n)}_{\Pi_n (i+1)} -\mu_{N_n} \right)  \right] = o(1) \,\, ,
$$

\noindent and 

$$
\var \left( \frac{1}{\sqrt{n}} \sum_{i=1} ^{n-1} \left( Z^{(N_n)}_{\Pi_n (i)} -\nu_{N_n} \right)\left( Y^{(N_n)}_{\Pi_n (i+1)} -\mu_{N_n} \right)  \right) = o(1) \,\, .
$$

\noindent By Chebyshev's inequality, it follows that

\beq
 \frac{1}{\sqrt{n}} \sum_{i=1} ^{n-1} \left( Z^{(N_n)}_{\Pi_n (i)} -\nu_{N_n} \right)\left( Y^{(N_n)}_{\Pi_n (i+1)} -\mu_{N_n} \right) = o_p (1) \, \, .
\label{yn.zn.3}
\eeq

\noindent By Slutsky's theorem, combining (\ref{yn.zn.1}), (\ref{yn.zn.2}), and (\ref{yn.zn.3}), we see that (\ref{yn.zn.conv}) holds. 

Since, for all $n, \, N \in \nn$, 

$$
\frac{1}{\sqrt{n}} \sum_{i=1} ^{n-1} Z^{(N)}_{\Pi_n (i)} Y^{(N)}_{\Pi_n (i+1)} \overset{d} =
\frac{1}{\sqrt{n}} \sum_{i=1} ^{n-1} Y^{(N)}_{\Pi_n (i)} Z^{(N)}_{\Pi_n (i+1)} \, \, ,
$$

\noindent we have that 

\beq
\frac{1}{\sqrt{n}} \sum_{i=1} ^{n-1} Y^{(N_n)}_{\Pi_n (i)} Z^{(N_n)}_{\Pi_n (i+1)} + \sqrt{n} \mu_{N_n} ^2 = o_p(1) \, \, .
\label{yn.zn.conv.sym}
\eeq

\noindent Now, we consider the fourth term in (\ref{perm.4.terms}). Note that 

$$
\begin{aligned}
&\frac{1}{\sqrt{n}} \sum_{i=1} ^{n-1} Z^{(N_n)}_{\Pi_n (i)} Z^{(N_n)}_{\Pi_n (i+1)} - \sqrt{n} \mu_{N_n} ^2 \\
&=\frac{1}{\sqrt{n}} \sum_{i=1} ^{n-1} \left( Z^{(N_n)}_{\Pi_n (i)} -\nu_{N_n} \right)\left( Z^{(N_n)}_{\Pi_n (i+1)} -\nu_{N_n} \right) + \frac{\nu_{N_n}}{\sqrt{n}} \sum_{i=1} ^{n-1} \left( Z^{(N_n)}_{\Pi_n (i)} -\nu_{N_n} \right) + \\
&+ \frac{\nu_{N_n}}{\sqrt{n}} \sum_{i=1} ^{n-1} \left( Z^{(N_n)}_{\Pi_n (i+1)} -\nu_{N_n} \right) \, \, .
\end{aligned}
$$

\noindent By an identical argument to the one used to show (\ref{yn.zn.conv}), we have that 

\beq
\frac{1}{\sqrt{n}} \sum_{i=1} ^{n-1} Z^{(N_n)}_{\Pi_n (i)} Z^{(N_n)}_{\Pi_n (i+1)} - \sqrt{n} \mu_{N_n} ^2 = o_p(1) \, \, .
\label{zn.conv}
\eeq

\noindent By (\ref{perm.4.terms}), we have that

$$
\begin{aligned}
&\frac{1}{\sqrt{n}} \sum_{i=1} ^{n-1} X_{\Pi_n (i)} X_{\Pi_n (i+1)} \\
&=\frac{1}{\sqrt{n}} \sum_{i=1} ^{n-1} Y^{(N_n)}_{\Pi_n (i)} Y^{(N_n)}_{\Pi_n (i+1)}  + 
\frac{1}{\sqrt{n}} \sum_{i=1} ^{n-1} Z^{(N_n)}_{\Pi_n (i)} Y^{(N_n)}_{\Pi_n (i+1)} + \\
&+\frac{1}{\sqrt{n}} \sum_{i=1} ^{n-1} Y^{(N_n)}_{\Pi_n (i)} Z^{(N_n)}_{\Pi_n (i+1)} +
\frac{1}{\sqrt{n}} \sum_{i=1} ^{n-1} Z^{(N_n)}_{\Pi_n (i)} Z^{(N_n)}_{\Pi_n (i+1)}  \, \, ,
\end{aligned}
$$

\noindent By (\ref{2.yn.conv}),

\beq
\begin{aligned}
&\frac{1}{\sqrt{n}} \sum_{i=1} ^{n-1} Y^{(N_n)}_{\Pi_n (i)} Y^{(N_n)}_{\Pi_n (i+1)} \\
&= \frac{\sigma_{N_n} ^2}{\sqrt{n}} \sum_{i=1} ^{n-1} \left(\frac{Y^{(N_n)}_{\Pi_n (i)} - \mu_{N_n}}{\sigma_{N_n}} \right)\left(\frac{Y^{(N_n)}_{\Pi_n (i+1)} - \mu_{N_n}}{\sigma_{N_n}} \right) + \frac{\mu_{N_n}}{ \sqrt{n}}\sum_{i=1} ^{n-1} \left(Y^{(N_n)}_{\Pi_n(i)} - \mu_{N_n} \right) + \\
&+\frac{\mu_{N_n}}{ \sqrt{n}}\sum_{i=1} ^{n-1} \left(Y^{(N_n)}_{\Pi_n(i+1)} - \mu_{N_n} \right) + \sqrt{n} \mu_{N_n} ^2 \\
&=  \frac{\sigma_{N_n} ^2}{\sqrt{n}} \sum_{i=1} ^{n-1} \left(\frac{Y^{(N_n)}_{\Pi_n (i)} - \mu_{N_n}}{\sigma_{N_n}} \right)\left(\frac{Y^{(N_n)}_{\Pi_n (i+1)} - \mu_{N_n}}{\sigma_{N_n}} \right) + \sqrt{n} \mu_{N_n} ^2 + o_p (1) \, \, .
\end{aligned}
\label{fin.a.perm.1}
\eeq

\noindent Hence, substituting the results of (\ref{yn.zn.conv}), (\ref{yn.zn.conv.sym}), and (\ref{zn.conv}) into (\ref{fin.a.perm.1}), 

$$
\frac{1}{\sqrt{n}} \sum_{i=1} ^{n-1} X_{\Pi_n (i)} X_{\Pi_n (i+1)} = \frac{\sigma_{N_n} ^2}{\sqrt{n}} \sum_{i=1} ^{n-1} \left(\frac{Y^{(N_n)}_{\Pi_n (i)} - \mu_{N_n}}{\sigma_{N_n}} \right)\left(\frac{Y^{(N_n)}_{\Pi_n (i+1)} - \mu_{N_n}}{\sigma_{N_n}} \right) + o_p(1) \, \, .
$$

\noindent By (\ref{perm.n.inc}), (\ref{mu.sig.yn}), and Slutsky's theorem for randomization distributions (\cite{chung.perm}, Theorem 5.2), it follows that, for all $t \in \rr$, as $n \to \infty$,

$$
\hat{R}_n (t) \overset{p} \to \Phi(t) \, \, ,
$$

\noindent as claimed.
\end{proof}

The following results show the consistency of the variance estimators of $\gamma_1 ^2$.

\begin{lem} Let $\{X_n, \, n \in \nn\}$ be a stationary, $\alpha$-mixing sequence, with mean $\mu$, such that, for some $\delta >0$,

$$
\lVert X_1 \rVert_{4 + 2\delta} < \infty \, \, ,
$$

\noindent and 

$$
\sum_{n \geq 1} \alpha_X (n) ^{\frac{\delta}{2 + \delta} }< \infty  \, \, .
$$ 

\noindent Let 

$$
\tau ^2 = \var \left(X_1\right) + 2\sum_{k \geq 2} \cov (X_1, \, X_k )\,\,.
$$

\noindent Suppose that $\tau^2 >0$. Let $\{b_n, \, n \in \nn\}$ be such that, for all $n$, $b_n\leq n-1$, $b_n = o \left(\sqrt{n} \right)$, and $b_n \to \infty$ as $n \to \infty$. Let

$$
\begin{aligned}
\hat{\tau}_n ^2 &= \frac{1}{n} \sum_{i=1} ^n \left(X_i- \bar{X}_n\right) ^2 + \frac{2}{n}\sum_{j=1}^{b_n}\sum_{i=1} ^{n-j} \left(X_i - \bar{X}_n \right) \left(X_{i+j} -\bar{X}_n \right)\, \, . \\
\end{aligned}
$$

\noindent We have that, as $n \to \infty$, 

$$
\begin{aligned}
\hat{\tau}_n^2 \overset{p}\to \tau^2 \, \,.
\end{aligned}
$$

\label{var.stud.unperm}
\end{lem}

\begin{proof} Let 

$$
\tilde{\tau}_n ^2 = \frac{1}{n} \sum_{i=1} ^n \left(X_i- \mu \right) ^2 + \frac{2}{n}\sum_{j=1}^{b_n}\sum_{i=1} ^{n-j} \left(X_i - \mu\right) \left(X_{i+j} -\mu \right) \, \,.
 $$

\noindent We have that 

$$
\begin{aligned}
\hat{\tau}_n ^2 &= \tilde{\tau}_n ^2 +\left( \frac{2}{n}\sum_{j=1}^{b_n} (n-j) -1 \right) \left( \bar{X}_n - \mu \right)^2 - \frac{2\left(\bar{X}_n - \mu\right)}{n} \sum_{j =1}^{b_n} \sum_{i=1} ^{n-j} \left[ (X_i - \mu) + (X_{i + j} - \mu) \right] \, \, .
\end{aligned}
$$

\noindent By \cite{doukhan}, Section 1.2.2, Theorem 3, there exists a constant $C \in \rr_+$ such that

\beq
\begin{aligned}
\var\left( \bar{X}_n \right) &= \frac{1}{n} \var(X_1) + \frac{1}{n^2} \sum_{i <j} \cov \left(X_i, \, X_j \right)\\
&\leq O\left(\frac{1}{n} \right) + \frac{C}{n^2} \sum_{i < j} \alpha_X (j - i)^{\frac{\delta}{2+ \delta}} \\
&= O\left( \frac{1}{n} \right) \, \, ,
\end{aligned}
\label{mean.op.rootn}
\eeq

\noindent by the same argument as in (\ref{smallo.1.alpha.mean}). Hence, since $\ee \bar{X}_n = \mu$, 

$$
\bar{X}_n - \mu = O_p \left( \frac{1}{\sqrt{n}} \right) \, \,,
$$

\noindent and so

\beq
\tilde{\tau}_n ^2 = \hat{\tau}_n ^2 - \frac{2\left(\bar{X}_n - \mu\right)}{n} \sum_{j =1}^{b_n} \sum_{i=1} ^{n-j} \left[ (X_i - \mu) + (X_{i + j} - \mu) \right]+ o_p(1) \, \,.
\label{tau.tilde.piece0}
\eeq

\noindent Note that 

\beq
\begin{aligned}
\sum_{j =1}^{b_n} \sum_{i=1} ^{n-j} \left[ (X_i - \mu) + (X_{i + j} - \mu) \right] &= \sum_{j=1}^{b_n} \left( 2n\left( \bar{X}_n-\mu\right) -\sum_{i=1}^{n-j-1} (X_i- \mu) - \sum_{r = n- j + 1}^{n} (X_r- \mu) \right) \\
&= 2nb_n \left( \bar{X}_n-\mu\right) - \\
&-\sum_{i=1} ^n \left[ \left(b_n \wedge (n- i -1) \right)_+ + \left( b_n + i - n\right)_+\right] (X_i - \mu) \, \, .
\end{aligned}
\label{tau.tilde.piece1}
\eeq

\noindent By another application of \cite{doukhan}, Section 1.2.2, Theorem 3 and Chebyshev's inequality, we obtain that 

\beq
\sum_{i=1} ^n \left[ \left(b_n \wedge (n- i -1) \right)_+ + \left( b_n + i - (n-1)\right)_+\right] (X_i - \mu) = O_p\left( b_n \sqrt{n} \right) \, \, .
\label{tau.tilde.piece2} 
\eeq

\noindent Combining (\ref{tau.tilde.piece0}), (\ref{tau.tilde.piece1}), and (\ref{tau.tilde.piece2}), we have that

\beq
\hat{\tau}_n ^2 = \tilde{\tau}_n ^2 + o_p (1) \, \, .
\eeq

\noindent Therefore, assuming, without loss of generality, that $\mu = 0$, it suffices to show that 

$$
\tilde{\tau}_n ^2 \overset{p} \to \tau^2 \, \, .
$$

\noindent Fix $0 \leq j \leq b_n$. By stationarity, we have that 

\beq
\begin{aligned}
\ee \left[ \frac{1}{n} \sum_{i=1}^{n-j} X_i X_{i+j} \right] &= \ee X_1 X_{1 + j}  -\frac{j}{n} \ee X_1 X_{1+j} \\
&= \ee X_1 X_j - O\left( \frac{1}{n} \right) \, \, .
\end{aligned} 
\label{mean.tau.alpha.piece}
\eeq

\noindent By \cite{doukhan}, Section 1.2.2, Theorem 3, we also have that, for some constant $C \in \rr_+$,

$$
\begin{aligned}
\var \left( \frac{1}{n} \sum_{i=1}^{n-j} X_i X_{i+j} \right) &= \frac{n-j}{n}\var \left( X_1 X_{1+j} \right) +\frac{2}{n^2} \sum_{1\leq i < k \leq n - j } \cov\left(X_{i} X_{i+j}, \, X_k X_{k + j} \right) \\
&\leq O\left( \frac{1}{n} \right) + \frac{2}{n^2} \sum_{1 \leq i < k \leq n-j } \left| \cov\left(X_{i} X_{i+j}, \, X_k X_{k + j} \right)\right| \\
&\leq O\left( \frac{1}{n} \right) + \frac{C}{n^2} \sum_{1 \leq i < k \leq n-j} \alpha_X \left( \left( k - j - i \right)_+ \right)^\frac{\delta}{2 + \delta} \\
&= O \left( \frac{1}{n} \right) \, \, ,
\end{aligned}
$$

\noindent as in (\ref{smallo.1.alpha.mean}). In particular, by Chebyshev's inequality, it follows that 

\beq
\frac{1}{n} \sum_{i=1}^{n-j} X_i X_{i+j} = \ee X_1 X_{1 + j} + O_p\left( \frac{1}{\sqrt{n}} \right) \,\, .
\label{tau.alpha.piece.conv}
\eeq

\noindent Hence 

$$ 
\tilde{\tau}_n^2 = \tau^2 + O_p \left( \frac{b_n}{\sqrt{n}} \right) \, \, ,
$$

\noindent and the result follows.
\end{proof}

\begin{cor} Let $\{X_n\}$ be a stationary, $\alpha$-mixing sequence such that, for some $\delta >0$, 

$$
\lVert X_1 \rVert_{8 + 4\delta} < \infty \, \, ,
$$

\noindent and 

$$
\sum_{n \geq 1} \alpha_X (n) ^{\frac{\delta}{2 + \delta} }< \infty  \, \, .
$$

\noindent For $i \in \nn$, let $Y_i = \left(X_i  - \bar{X}_n \right) \left( X_{i+1}  - \bar{X}_n\right)$. Let $b_n= o\left(\sqrt{n}\right)$ be such that, as $n \to \infty$, $b_n \to \infty$. Let

$$
\hat{T}_n ^2 = \frac{1}{n} \sum_{i=1} ^{n-1} \left(Y_i - \bar{Y}_n \right)^2 + \frac{2}{n}\sum_{j=1}^{b_n}\sum_{i=1} ^{n-j-1} \left(Y_i - \bar{Y}_n \right) \left(Y_{i+j} -\bar{Y}_n \right) \, \, .
$$

\noindent Let 

$$
{\tau}_1^2 = \var(X_1 X_2) + 2\sum_{k \geq 2} \cov \left( X_1X_2, \, X_k X_{k+1} \right) \, \, .
$$

\noindent We have that, as $n \to \infty$, 

$$
\hat{T}_n ^2 \overset{p} \to \tau_1^2 \, \, .
$$

\label{cor.stud.unperm}
\end{cor}

\begin{proof} Without loss of generality, we may assume that $\ee X_1 = 0$. Let $Z_i = X_i X_{i+1}$. Let $j \geq 0$. Noting that, by the same argument as in (\ref{mean.op.rootn}), $\bar{X}_n = O_p(1/\sqrt{n} )$, we have that 

$$
\begin{aligned}
Y_i &= Z_i - \bar{X}_n \left(X_i + X_{i+1} \right) + \bar{X}_n ^2  && = Z_i + O_p\left( \frac{1}{\sqrt{n}} \right)  \\
\bar{Y}_n &= \bar{Z}_n - \bar{X}_n ^2 + o_p \left( \frac{1}{n} \right) &&= \bar{Z}_n + O_p \left( \frac{1}{n} \right) \, \, .
\end{aligned}
$$

\noindent Hence

$$
\begin{aligned}
\frac{1}{n} \sum_{i = 1} ^{n-j-1} \left(Y_i - \bar{Y}_n \right) \left(Y_{i+j} -\bar{Y}_n \right) &= \frac{1}{n} \sum_{i = 1} ^{n-j-1} \left(Z_i - \bar{Z}_n  \right) \left(Z_{i+j} -\bar{Z}_n \right) + O_p \left( \frac{1}{\sqrt{n}}\right) \, \, .
\end{aligned}
$$

\noindent In particular, since $b_n = o(\sqrt{n} )$, we have that 

$$
\hat{T}_n ^2 =\frac{1}{n} \sum_{i=1} ^{n-1} \left(Z_i - \bar{Z}_n  \right)^2 +  \frac{2}{n}\sum_{j=1} ^{b_n} \sum_{i = 1} ^{n-j-1} \left(Z_i - \bar{Z}_n  \right) \left(Z_{i+j} -\bar{Z}_n \right) + o_p(1) \, \, .
$$

\noindent Since \noindent $\alpha_Z (i) = \alpha_X (i-1)$, we may apply the result of Lemma \ref{var.stud.unperm}, and the result follows.
\end{proof}

Having shown consistency of the estimator $\hat{T}_n^2$, the consistency of estimators of the other terms in the definition of $\gamma_1^2$ follow similarly.

\begin{proof} [{\sc Proof of Lemma \ref{lem.stud.unperm.pap}}] By an identical argument to that used in the proof of Lemma \ref{var.stud.unperm} and Corollary \ref{cor.stud.unperm}, we have that 

$$
\begin{aligned}
\hat{K}_n ^2 &\overset{p} \to \kappa^2 \\
\hat{\nu}_n &\overset{p} \to \nu_1 \, \, .
\end{aligned}
$$

\noindent As a consequence of Theorem \ref{thm.unperm.rho.lim}, it follows that $\hat{\rho}_n$ is a consistent estimator of $\rho_1$. Hence, by the continuous mapping theorem and Slutsky's theorem, we have that 

$$
\hat{\gamma}_n^2 \overset{p} \to \gamma_1 ^2 \, \, ,
$$

\noindent as required.
\end{proof}

We now proceed to show that, when a random permutation is applied, the studentizing factor converges in probability to 1.

\begin{lem} In the setting of Corollary \ref{cor.stud.unperm}, let $\Pi_n \sim \text{Unif}(S_n)$, independent of the sequence $\{X_i, \, i \in \nn\}$. For $i \in \nn$, let

$$
Z_i = \left(X_{\Pi_n (i)} - \bar{X}_n \right) \left( X_{\Pi_n(i+1)} -\bar{X}_n \right) \, \, .
$$

\noindent Let 

$$
\hat{T}_n ^2 = \frac{1}{n} \sum_{i=1} ^{n-1} \left(Z_i- \bar{Z}_n\right) ^2 + \frac{2}{n}\sum_{j=1}^{b_n}\sum_{i=1} ^{n-j-1} \left(Z_i- \bar{Z}_n \right) \left(Z_{i+j} -\bar{Z}_n \right) \, \, .
$$

\noindent Then, as $n \to \infty$,

$$
\hat{T}_n ^2 \overset{p} \to \var(X_1)^2 \, \, .
$$

\label{var.perm.stud}
\end{lem}

\begin{proof} We may assume, without loss of generality, that $\ee X_1 = 0$. Let

$$
\Gamma_i = X_{\Pi_n (i)}  X_{\Pi_n (i+1)} \, \, .
$$

\noindent By the same argument as in the proof of Corollary \ref{cor.stud.unperm},

$$
\begin{aligned}
\hat{T}_n ^2 = \frac{1}{n} \sum_{i=1} ^n \left(\Gamma_i - \bar{\Gamma}_n \right)^2 + \frac{2}{n} \sum_{j=1} ^{b_n} \sum_{i=1} ^{n-j-1}\left(\Gamma_i - \bar{\Gamma}_n \right)\left(\Gamma_{i+j} - \bar{\Gamma}_n \right) + o_p(1) \, \, .
\end{aligned}
$$

\noindent By Lemma \ref{mean.var.alpha.conv}, we have that $\bar{\Gamma}_n = O_p\left(1/\sqrt{n} \right)$. Let 

$$
\tilde{T}_n^2 = \sum_{i=1} ^n \Gamma_i ^2  +  \frac{2}{n} \sum_{j=1} ^{b_n} \sum_{i=1} ^{n-j-1}\Gamma_i \Gamma_{i+j} \, \, .
$$

\noindent Since $b_n = o \left( \sqrt{n} \right)$, we have that 

$$
\hat{T}_n ^2 = \tilde{T}_n ^2 + o_p(1) \, \, .
$$

\noindent It therefore suffices to show that 

$$
\tilde{T}_n ^2 \overset{p} \to \var \left( X_1\right) ^2\, \, .
$$

\noindent Dividing $\tilde{T}_n ^2$ by $\var \left( X_1\right) ^2$, we may assume, without loss of generality, that $\var \left( X_1\right) = 1$, and subsequently it suffices to show that 

$$
\tilde{T}_n ^2 \overset{p} \to 1 \, \, .
$$

\noindent We have 

\beq
\tilde{T}_n ^2 = \frac{1}{n} \sum_{i=1} ^{n-1} X^2_{\Pi_n (i)} X^2_{\Pi_n (i+1)} + \frac{2}{n} \sum_{j = 1}^ {b_n} \sum_{i = 1} ^{n-j-1} X_{\Pi_n (i)} X_{\Pi_n (i+1)}X_{\Pi_n (i+j)} X_{\Pi_n (i+j+1)} \, \, .
\label{tilde.t.perm}
\eeq

\noindent Note that, applying the result of Lemma \ref{mean.var.alpha.conv} to the sequence $\left\{X_i ^2 - 1, \, i \in \nn\right\}$, and by the same argument as in (\ref{mean.tau.alpha.piece}), we have that

$$
\begin{aligned}
\frac{1}{n} \sum_{i=1} ^{n-1} X^2_{\Pi_n (i)} X^2_{\Pi_n (i+1)} &= -1 + \frac{2}{n} \sum_{i=1}^n X_i ^2 + \frac{1}{n} \sum_{i=1} ^{n-1} \left(X^2_{\Pi_n (i)} -1 \right) \left(X^2_{\Pi_n (i+1)}  -1 \right) + o_p(1) \\ 
&= 1 + o_p(1) \, \, .
\end{aligned}
$$

\noindent In particular, combining with (\ref{tilde.t.perm}), since $b_n = o(\sqrt{n})$, it suffices to show that, for each $1 \leq j \leq b_n$,

$$
\frac{1}{\sqrt{n}} \sum_{i=1} ^{n-j-1} X_{\Pi_n (i)} X_{\Pi_n (i+1)}X_{\Pi_n (i+j)} X_{\Pi_n (i+j+1)} = O_p \left( 1 \right) \, \, .
$$

\noindent For $j = 1$, for all $i$, we have that 

$$
X_{\Pi_n (i)} X_{\Pi_n (i+1)}X_{\Pi_n (i+j)} X_{\Pi_n (i+j+1)} \overset{d} = X_{\Pi_n (1) } X_{\Pi_n (2)} ^2 X_{\Pi_n (3)} \, \, .
$$

\noindent By (\ref{j1.var.term}), we have that, for some constant $K \in \rr_+$,

\beq
\begin{aligned}
\left| \ee X_{\Pi_n (1) } X_{\Pi_n (2)} ^2 X_{\Pi_n (3)}  \right | &\leq \frac{K}{n} \sum_{r = 1} ^{n-1} \alpha_X (r) ^\frac{\delta}{2 + \delta} \\
&= O\left( \frac{1}{n} \right) \, \, .
\end{aligned}
\label{mean.j1.student}
\eeq

\noindent It follows that 

$$
\ee \left[ \frac{1}{\sqrt{n}} \sum_{i=1} ^ {n-2} X_{\Pi_n (i)} X_{\Pi_n (i+1)} ^2 X_{\Pi_n (i+2)} \right] = o(1) \, \, .
$$

\noindent We also have that, due to the moment conditions on the sequence $\{X_n\}$ and the Cauchy-Schwarz inequality, 

$$
\begin{aligned}
&\var\left(\frac{1}{\sqrt{n}} \sum_{i=1} ^ {n-2} X_{\Pi_n (i)} X_{\Pi_n (i+1)} ^2 X_{\Pi_n (i+2)}  \right) \\&= \frac{1}{n} \sum_{i=1} ^{n-2} \var\left( X_{\Pi_n (i)} X_{\Pi_n (i+1)} ^2 X_{\Pi_n (i+2)} \right) + \frac{2}{n} \sum_{i < j} \cov\left( X_{\Pi_n (i)} X_{\Pi_n (i+1)} ^2 X_{\Pi_n (i+2)}, \, X_{\Pi_n (j)} X_{\Pi_n (j+1)} ^2 X_{\Pi_n (j+2)} \right) \\
&= \frac{1}{n} \sum_{i < j - 2} \cov\left( X_{\Pi_n (i)} X_{\Pi_n (i+1)} ^2 X_{\Pi_n (i+2)}, \, X_{\Pi_n (j)} X_{\Pi_n (j+1)} ^2 X_{\Pi_n (j+2)} \right) + O(1) \, \, .
\end{aligned}
$$

\noindent Note that, for all $i < j-2$, the covariances in the above sum are equal. Therefore, in order to conclude the proof, by Chebyshev's inequality it suffices to show that 

\beq
\cov\left( X_{\Pi_n (i)} X_{\Pi_n (i+1)} ^2 X_{\Pi_n (i+2)}, \, X_{\Pi_n (j)} X_{\Pi_n (j+1)} ^2 X_{\Pi_n (j+2)} \right)= O\left(\frac{1}{n}\right) \, \, .
\label{var.j.order.student}
\eeq

\noindent By (\ref{mean.j1.student}), it suffices to show that, for $i < j-2$, 

$$
\ee \left[X_{\Pi_n (i)} X_{\Pi_n (i+1)} ^2 X_{\Pi_n (i+2)}X_{\Pi_n (j)} X_{\Pi_n (j+1)} ^2 X_{\Pi_n (j+2)}\right] = O\left( \frac{1}{n} \right) \, \, .
$$

\noindent There exist constants $C_1, \, C_2, \, C_3, \, C_4, \, C_5 \in \rr_+$ such that 

\beq
\begin{aligned}
&\ee \left[X_{\Pi_n (i)} X_{\Pi_n (i+1)} ^2 X_{\Pi_n (i+2)}X_{\Pi_n (j)} X_{\Pi_n (j+1)} ^2 X_{\Pi_n (j+2)}\right] = \frac{C_1}{n^6} \sum_{i, \, j, \, k, \, l, \, m, \, p\text{\, dist.}}\ee \left[X_{i} X_{j}  X_{k}X_{l} X_{m} ^2 X_{p} ^2\right] \\
&= \frac{C_2}{n^6} \sum_{m < i < j < k < l < p} \ee \left[X_{i} X_{j}  X_{k}X_{l} X_{m} ^2 X_{p} ^2\right] +  \frac{C_3}{n^6} \sum_{ i < j < k < l < m, \,  i < p < m} \ee \left[X_{i} X_{j}  X_{k}X_{l} X_{m} ^2 X_{p} ^2\right] + \\
&+  \frac{C_4}{n^6} \sum_{m < i < j < k < l, \, m < p < l} \ee \left[X_{i} X_{j}  X_{k}X_{l} X_{m} ^2 X_{p} ^2\right]
 + \frac{C_5}{n^6} \sum_{i < j < k < l, \, i < m < l, \,i < p < l}\ee \left[X_{i} X_{j}  X_{k}X_{l} X_{m} ^2 X_{p} ^2\right] \, \, .
 \end{aligned}
\label{var.j1.student}
\eeq

\noindent We begin by considering the fourth sum on the right hand side of (\ref{var.j1.student}). Fix $l- i = r \geq 5$. There are $(n-r)$ such choices for $(i, \, l)$. By stationarity, we may assume, without loss of generality, that $m =1$, $p = r+1$. Assume, without loss of generality, that $m < p$. By \cite{doukhan},  Section 1.2.2, Theorem 3, there exists a constant $K \in \rr_+$ such that

$$
\left|\ee \left[X_{i} X_{j}  X_{k}X_{l} X_{m} ^2 X_{p} ^2\right] \right| \leq K \alpha_X \left( \max\{ \min\{j, \, m \} -1, \, l - \max\{k, \, p \} \} \right) ^\frac{\delta}{2 + \delta} \, \, .
$$

\noindent Fix $s = \max\{ \min\{j, \, m \} -1, \, l - \max\{k, \, p \} \}$. There are, at most, $2s$ choices for the second smallest and second largest index combined, and there are at most $(n-s)^2$ choices for the remaining two indices. Hence

$$
\begin{aligned}
\frac{1}{n^6} \sum_{i < j < k < l, \, i < m < l, \,i < p < l} \left | \ee \left[X_{i} X_{j}  X_{k}X_{l} X_{m} ^2 X_{p} ^2\right] \right | &\leq \frac{1}{n^6} \sum_{r = 5} ^{n-1} (n-r) \sum_{s=1}^r 2s (n-s)^2 \alpha_X (s) ^\frac{\delta}{2 + \delta} \\
&\leq  \frac{2}{n^6} \sum_{s = 1} ^{n-1} s (n-s)^4 \alpha_X (s)^\frac{\delta}{2 + \delta} \\
&= O\left( \frac{1}{n} \right) \, \, .
\end{aligned}
$$

\noindent We now consider the third term on the right hand side of (\ref{var.j1.student}). Fix $l - \max\{p, \, k \}= r$. There are at most $n -r$ choices for the pair $\left( \max\{ p, \, k \} , \, l \right)$. Now, there are at most $C \cdot (n-r)^4$ choices for the remaining indices, for some universal constant $C \in \rr_+$. Hence, once more applying \cite{doukhan}, Section 1.2.2, Theorem 3, we have that 

$$
\begin{aligned}
\frac{1}{n^6} \sum_{m < i < j < k < l, \, m < p < l} \left|\ee \left[X_{i} X_{j}  X_{k}X_{l} X_{m} ^2 X_{p} ^2\right] \right| &\leq \frac{C}{n^6} \sum_{r=1} ^{n-1} \left(n-r \right)^5 \alpha_X (r)^\frac{\delta}{2 + \delta} \\
&= O\left( \frac{1}{n} \right) \, \, .
\end{aligned}
$$

\noindent Similarly, the second term in (\ref{var.j1.student}) is also $O(1/n)$.

We conclude our analysis for $j = 1$ by considering the first sum in (\ref{var.j1.student}). A double application of \cite{doukhan}, Section 1.2.2, Theorem 3, and the same argument as for the previous sums, gives us that, for some universal constant $C \in \rr_+$, 

$$
\begin{aligned}
&\frac{1}{n^6} \sum_{m < i < j < k < l < p} \left | \ee \left[X_{i} X_{j}  X_{k}X_{l} X_{m} ^2 X_{p} ^2\right] \right | \\
&\leq \frac{C}{n^6} \sum_{r = 1} ^{n-1} (n-r)^5 \alpha_X (r) ^\frac{\delta}{2 + \delta} + \frac{1}{n^6} \sum_{m < i < j < k < l < p}  \left | \ee \left[X_{i} X_{j}  X_{k}X_{l} X_{p} ^2\right] \right | \\
&\leq O\left( \frac{1}{n} \right) + \frac{C}{n^5} \sum_{r = 1} ^{n-1} (n-r) ^4 \alpha_X (r) ^\frac{\delta}{2 + \delta} \\
&= O\left( \frac{1}{n} \right) \, \,.
\end{aligned}
$$

\noindent It follows that (\ref{var.j.order.student}) holds, by Chebyshev's inequality. We turn our attention to the corresponding result for $j >1$. Note that, for all $i$, for all $j > 1$, $X_{\Pi_n (i)} X_{\Pi_n (i+1)}X_{\Pi_n (i+j)} X_{\Pi_n (i+j+1)} \overset{d} = X_{\Pi_n (i)} X_{\Pi_n (2)}X_{\Pi_n (3)} X_{\Pi_n (4)}$. By (\ref{final.sum.alpha.var}), we have that 

\beq
\ee \left[ \frac{1}{\sqrt{n}} \sum_{i=1} ^{n-j-1} X_{\Pi_n (i)} X_{\Pi_n (i+1)}X_{\Pi_n (i+j)} X_{\Pi_n (i+j+1)} \right] = o(1) \, \, .
\label{mean.j.big.student}
\eeq

\noindent By the Cauchy-Schwarz inequality,
\beq
\begin{aligned}
&\var\left(\frac{1}{\sqrt{n}} \sum_{i=1} ^{n-j-1} X_{\Pi_n (i)} X_{\Pi_n (i+1)}X_{\Pi_n (i+j)} X_{\Pi_n (i+j+1)}\right) = \var \left(X_{\Pi_n (1)} X_{\Pi_n (2)}X_{\Pi_n (3)} X_{\Pi_n (4)}\right) + \\ &+\frac{1}{n}\sum_{i \neq k} \cov \left( X_{\Pi_n (i)} X_{\Pi_n (i+1)}X_{\Pi_n (i+j)} X_{\Pi_n (i+j+1)}, \, X_{\Pi_n (k)} X_{\Pi_n (k+1)}X_{\Pi_n (k+j)} X_{\Pi_n (k+j+1)}\right) \\
&= \frac{1}{n}\sum_{i \in \{k \pm 1, \, k \pm j, \, k \pm (j+1)  \}^C}  \cov \left( X_{\Pi_n (i)} X_{\Pi_n (i+1)}X_{\Pi_n (i+j)} X_{\Pi_n (i+j+1)}, \, X_{\Pi_n (k)} X_{\Pi_n (k+1)}X_{\Pi_n (k+j)} X_{\Pi_n (k+j+1)}\right) + \\
&+ O(1) \, \, .
\end{aligned}
\label{eq.last.covariance.done}
\eeq

\noindent By stationarity of the sequence $\{X_i\}$ and Chebyshev's inequality, in order to conclude the proof, it suffices to show that

$$
\cov\left( X_{\Pi_n (1)} X_{\Pi_n (2)}X_{\Pi_n (3)} X_{\Pi_n (4)}, \, X_{\Pi_n (5)} X_{\Pi_n (6)}X_{\Pi_n (7)} X_{\Pi_n (8)}\right) = O\left( \frac{1}{n} \right) \, \, ,
$$

\noindent since the covariances in the sum of the last line of (\ref{eq.last.covariance.done}) are equal. By (\ref{mean.j.big.student}), it therefore suffices to show that 

$$
\ee\left[X_{\Pi_n (1)} X_{\Pi_n (2)}X_{\Pi_n (3)} X_{\Pi_n (4)}X_{\Pi_n (5)} X_{\Pi_n (6)}X_{\Pi_n (7)} X_{\Pi_n (8)}\right] = O\left( \frac{1}{n} \right) \, \,.
$$

\noindent There exists a universal constant $C \in \rr_+$ such that 

$$
\begin{aligned}
&\left|\ee\left[X_{\Pi_n (1)} X_{\Pi_n (2)}X_{\Pi_n (3)} X_{\Pi_n (4)}X_{\Pi_n (5)} X_{\Pi_n (6)}X_{\Pi_n (7)} X_{\Pi_n (8)} \right] \right| \\&\leq \frac{C}{n^8} \sum_{i < j <k <l < p < q < r < s} \left| \ee X_ i X_j X_k X_l X_p X_q X_r X_s\right|
\end{aligned}
$$

\noindent Fix $j - i = r$. There are $(n-r)$ choices (at most) for the pair $(i, \, j)$. Subsequently, there are (at most) $(n-r)^6$ choices for the remaining indices. Hence, by \cite{doukhan}, Section 1.2.2, Theorem 3, we have that, for some universal constant $C \in \rr_+$, 

$$
\begin{aligned}
\frac{C}{n^8} \sum_{i < j <k <l < p < q < r < s} \left| \ee X_ i X_j X_k X_l X_p X_q X_r X_s\right| &\leq \frac{C}{n} \sum_{r = 1} ^{n-1} \left( \frac{n-r}{n} \right)^7 \alpha_X (r) ^\frac{\delta}{2 + \delta} \\
&= O \left( \frac{1}{n} \right) \, \, .
\end{aligned}
$$

\noindent The result follows.
\end{proof}

\begin{proof}  [{\sc Proof of Lemma \ref{pap.var.perm.stud}}] 

By an identical argument to the one used in the proof of Lemma \ref{var.perm.stud}, we have that 

$$
\begin{aligned}
\hat{K}_n ^2 &\overset{p} \to \var\left(X_1 ^2 \right) \\
\hat{\nu}_1 &\overset{p} \to 0 \, \, .
\end{aligned} 
$$

\noindent Assuming, without loss of generality, that $\ee X_1 = 0$, we have that, as a consequence of Theorem \ref{alpha.perm}, Hoeffding's condition (see \cite{chung.perm}, Theorem 5.1), and Slutsky's Theorem, the sample autocorrelation $\hat{\rho}_n$ under random permutation satisfies

$$
\hat{\rho}_n \overset{p}\to 0 \, \, .
$$

\noindent Hence, by Slutsky's theorem and Lemma \ref{var.perm.stud}, it follows that 

$$
\hat{\gamma}_n ^2 \overset{p} \to 1 \, \, ,
$$

\noindent as required.
\end{proof}

Using the previous results, we may now prove Theorem \ref{perm.test.big.thm}.

\begin{proof} [{\sc Proof of Theorem \ref{perm.test.big.thm}}] 

By Theorem \ref{alpha.perm}, Lemma \ref{pap.var.perm.stud}, and Slutsky's theorem for randomization distributions (\cite{chung.perm}), (\ref{stud.perm.test.1d}) holds.

By Theorem \ref{thm.unperm.rho.lim}, Lemma \ref{lem.stud.unperm.pap}, and Slutsky's theorem, (\ref{main.test.stat.1d}) holds.
\end{proof}

We illustrate the proof of Theorem \ref{1d.local.perm}, which proceeds almost identically to the proof of Theorem \ref{alpha.bd.perm}.

\begin{proof}[{\sc Proof of Theorem \ref{1d.local.perm}}] 

The proof is exactly the same as that of Theorem \ref{alpha.bd.perm}, since we have that, for each fixed $r$, we have that there exists a constant $C \in \rr_+$ such that

$$
\sum_{i=1} ^{r -1}\alpha_{X^{(r)}} (i)< C \, \, ,
$$

\noindent and the result of \cite{wald1943} still applies. The only difference is that the definition of the random variable $B(m)$ (see \ref{b.1d.bd.perm}) becomes

$$
B(m) = \sum_{i=1} \frac{1}{m}\sum_{i=1} ^m \left(X_i^{(m)}\right)^2 - \frac{1}{m^2} \left(\sum_{j=1} ^m X_j ^{(m)}\right)^2 \, \, .
$$

\noindent The remainder of the proof is identical.
\end{proof}

\section{Testing multiple lags}

\begin{proof} [{\sc Proof of Theorem \ref{thm.multiple}}] We may assume, without loss of generality, that For $i \in \nn$, and $\{a_i, \, i \in  \{0,\, 1, \, \dots, \, k\} \} \subset \rr^{k+1}$, let 

$$
Y_i =  \sum_{k = 0}^r a_i \left(X_i - \bar{X}_n \right) \left(X_{i+k} - \bar{X}_n \right) \, \, .
$$

\noindent Let 

$$
Z_i = \sum_{k = 0}^r a_i X_i X_{i+k} \, \, .
$$

\noindent By the same argument as in the proof of Theorem \ref{thm.unperm.rho.lim}, we have that

$$
Y_i = Z_i + o_p (\sqrt{n}) \, \, .
$$

\noindent Also, note that the sequence $(Z_i)$ is $\alpha$-mixing, with $\alpha$-mixing coefficients given by 

$$
\alpha_Z (i) = \alpha_X (i - r) \, \, .
$$

\noindent Noting also that 

$$
E[Z_i] = a_0 \sigma^2 + \sum_{k=1} ^r a_k \rho_k \sigma^2 \, \, ,
$$

\noindent by applying Ibragimov's central limit theorem as in the proof of Theorem \ref{thm.unperm.rho.lim}, we may conclude that 

$$
\frac{1}{\sqrt{n}} \sum_{i=1} ^{n-1} Y_i - \sqrt{n} \left( a_0 \sigma^2 + \sum_{k=1} ^r a_k  \rho_k \sigma^2 \right) \overset{d} \to N\left(0, \, \tilde{\sigma}^2\right) \, \, ,
$$

\noindent where

$$
\begin{aligned}
\tilde{\sigma}^2 &= \var(Y_1 ) + 2 \sum_{k >1} \cov(Y_1, \, Y_k) \\
&=\sum_{k = 0} ^r  a_k^2 \cdot \var(X_1 X_{1 + k} ) + 2\sum_{j < k} a_j a_k \cov(X_1 X_{1+j}, \, X_1 X_{1+k}) + 2 \sum_{l > 1} \sum_{j = 0}^r \sum_{k=0}^r a_j a_k \cov(X_1 X_{1 + j} , \, X_{l }X_{l + k} ) \, \, .
\end{aligned}
$$

\noindent Since the $a_i$ were arbitrary, it follows that, by the Cram\'er-Wold device,  

$$
\sqrt{n} \left [\hat{\sigma}_n ^2 \begin{pmatrix}
1\\
\hat{\rho}_1\\
\vdots \\
\hat{\rho}_r \end{pmatrix}  -\sigma^2 \begin{pmatrix}
1\\
{\rho}_1\\
\vdots \\
{\rho}_r \end{pmatrix} \right]
\overset{d} \to N(0, \, \Sigma) \, \, .
$$

\noindent We may now apply the delta method, with the function $h(x_0, \, \dots, \, x_r) = \left(x_1 / x_0, \, x_2/x_0, \, \dots, \, x_r/ x_0\right)$, to conclude that

$$
\sqrt{n}  \left [ \begin{pmatrix}
\hat{\rho}_1\\
\vdots \\
\hat{\rho}_r \end{pmatrix}  - \begin{pmatrix}
{\rho}_1\\
\vdots \\
{\rho}_r \end{pmatrix} \right] \overset{d} \to N\left(0, \, A^T \Sigma A\right) \, \, , 
$$

\noindent where 

$$
A = \begin{pmatrix} 
-\frac{\rho_1}{\sigma^4} & \dots & \dots &\dots&- \frac{ \rho_r}{\sigma^4} \\
\frac{1}{\sigma^2} & 0 & \dots&\dots  & 0 \\
0 & \frac{1}{\sigma^2} & 0 &\dots & 0 \\
\vdots &\vdots &\vdots & \vdots & \vdots \\
0 & \dots &\dots &\dots & \frac{1}{\sigma^2} 
\end{pmatrix} \, \, ,
$$

\noindent as required.
\end{proof}

\addcontentsline{toc}{section}{Bibliography}

\bibliographystyle{apalike}

\bibliography{perm_test_time_series}